\documentclass[11pt]{article}

\usepackage{amssymb,amsmath,amsthm,amsopn,amsfonts}
\usepackage[T1]{fontenc}
\usepackage{ae,aecompl}
\usepackage{bm}

\usepackage{algorithmic}
\usepackage{algorithm}

\usepackage{graphicx}

\usepackage[margin=20pt, justification=justified, font=small,labelfont=bf,labelsep=space]{caption}
\usepackage{float}
\usepackage{placeins}
\usepackage[percent]{overpic}  
\usepackage{morefloats}

\usepackage{fancybox,layout,color}
\usepackage{fancyhdr}

\usepackage{scrtime}

\usepackage{setspace}	

\usepackage{url}
\usepackage{placeins}

\usepackage{color}
\usepackage{colortbl}

\usepackage{mathtools}

\DeclarePairedDelimiter\floor{\lfloor}{\rfloor}

\usepackage{ulem}  
\usepackage{cancel}  

\numberwithin{equation}{section}
\definecolor{purple}{RGB}{160,32,40}

\newtheorem{teo}{Theorem}[section]
\newtheorem{nota}[teo]{Remark}
\newtheorem{ex}[teo]{Example}
\newtheorem{coro}[teo]{Corollary}
\newtheorem{defi}[teo]{Definition}
\newtheorem{prop}[teo]{Proposition}

\newcommand{\R}{\ensuremath{\mathbb{R}} }

\newcommand{\dist}{\mathrm{dist}}

\newcommand{\osc}{\mathrm{osc}}
\newcommand{\epi}{\mathrm{epi}}
\newcommand{\Aff}{\mathrm{Aff}}

\DeclareMathOperator{\co}{\mathsf{co}}

\DeclareMathOperator{\diam}{\mathsf{diam}}

\title{Compensated Convexity on Bounded Domains,\\  Mixed Moreau
Envelopes and Computational Methods}

\author{\normalsize  Kewei Zhang\thanks{School  of Mathematical Sciences,
University of Nottingham, University Park, Nottingham, NG7 2RD, UK},\,
Antonio Orlando\thanks{CONICET, Departamento de Bioingenier\`ia,
Universidad Nacional de Tucum\'an, Argentina}
\,and 
Elaine Crooks\thanks{Department of Mathematics, Swansea University,
Singleton Park, Swansea, SA2 8PP, UK}
}

\date{ }

\begin{document}

\maketitle


\singlespacing 


\pagestyle{fancy}
\fancyhead{}
\cfoot{\thepage}


\begin{abstract}
Compensated convex transforms have been introduced for extended real valued functions defined over $\R^n$. In their application to 
image processing, 
interpolation and shape interrogation, where one deals with functions defined over a bounded domain, 
one was making the implicit assumption that the function coincides with its transform at the boundary of the data domain. 
In this paper,
we introduce local compensated convex transforms for functions defined in bounded open convex subsets $\Omega$
of $\R^n$ by making specific extensions of the function to the whole space, and establish their relations to 
globally defined compensated convex transforms via the mixed critical
Moreau envelopes. We find that the compensated convex transforms of such extensions
coincide with the local compensated convex transforms in the closure of $\Omega$.
We also propose a numerical scheme for computing Moreau envelopes, establishing
convergence of the scheme with the rate of convergence depending on the regularity of the original function. 
We give an estimate of the number of iterations needed for computing the discrete Moreau envelope. 
We then apply the local compensated convex transforms to image processing and shape interrogation.
Our  results are compared with those obtained by using schemes based on computing the convex envelope from the 
original definition of compensated convex transforms.
\end{abstract}


\medskip
{\bf 2000 Mathematics Subjects Classification number}:
90C25, 90C26, 49J52, 52A41, 65D17

\medskip
{\bf Email}:  kewei.zhang@nottingham.ac.uk,  aorlando@herrera.unt.edu.ar, e.c.m.crooks@swansea.ac.uk

\normalsize


\setcounter{equation}{0}
\section{Introduction}\label{Sec.Intro}
This paper contains two parts. In the first part, 
we `localize' the global notions of compensated convex transforms 
defined over $\R^n$, which were first introduced in \cite{Zha08a,Zha08b}, 
by defining such transforms over bounded convex closed domains in $\mathbb{R}^n$ so that their values in the domain 
agree with the globally defined transforms applied to some special extensions of the function to $\mathbb{R}^n$. 
The motivation for such local definitions is mainly from applications to digital images and data arrays, 
where we have to consider functions defined on a rectangular box. 
In the second part, we propose a new scheme for the 
computation of Moreau envelopes and we prove its convergence. The scheme can be used to compute the 
Moreau based definitions of the local compensated convex transforms which we introduce in this paper 
or it can be useful on its own for general applications of Moreau envelopes.

The theory of compensated convexity transforms has been applied to, for example, digital image processing and computational
geometry. So far,  applications of the theory include the design of multiscale, parametrized,
geometric singularity extraction of ridges, valleys and edges from graphs of functions and from characteristic functions of 
closed sets in $\mathbb{R}^n$ \cite{ZOC15a}. Several robust methods have been developed to date, namely, for the extraction of 
the set of intersections between two or more smooth compact manifolds \cite{ZOC15b,ZCO16b}, for the extraction of the multiscale 
medial axis from geometric objects \cite{ZCO15c} and for the interpolation and approximation of 
sampled functions \cite{ZCO16a}.
By `robustness' here we mean the Hausdorff stability, that is, the error between the outputs obtained for two data samples  
is controlled by the Hausdorff distance between the two sampled input data sets.

In the applications mentioned above, the data domains are usually represented by boxes in $\R^n$. The  numerical schemes 
used in \cite{ZOC15a, ZOC15b, ZCO15c, ZCO16a, ZCO16b} for the  evaluation of the compensated convex transforms relied on the availability of schemes for computing the 
convex envelope of a given function and on the implicit assumption that the transforms coincide with the function at the boundary of the data
domain. In those works, we were applying such schemes mainly to demonstrate the numerical feasibility of the transforms
rather than \textit{(i)} designing efficient numerical schemes for their computation and/or
\textit{(ii)} analysing the effect of the boundary assumptions. In this paper we 
also 
address the more practical
question of accurately and effectively computing
the compensated convex transforms for functions defined on a bounded convex domain without using  numerical schemes 
that compute the convex envelope. In order to do so,
we will explore the alternative definitions of the compensated convex transforms based on the Moreau envelopes.
One of the advantages of this approach is that when it comes to the numerical
implementation by the scheme we propose, we obtain an estimate of the number of iterations which provides
the exact discrete Moreau envelope. This is different from the application of iterative schemes to
compute the convex envelope which can be shown to converge but for which no convergence rate is known to be available \cite{Obe08}.

Before we  introduce our local versions of compensated convex transforms on a bounded closed
convex domain, we recall from \cite{Zha08a} the notions of quadratic compensated convex
convex transforms in $\R^n$.
For a function $g:\R^n\to\R\cup\{\infty\}$
satisfying the growth condition $g(x)\geq -c_0-c_1|x|^2$, $x\in\mathbb{R}^n$, for some constants
$c_0$, $c_1>0$ with $|x|=\sqrt{x\cdot x}$ the Euclidean norm of $x\in\mathbb{R}^n$
and $x\cdot y$ the inner product on $\mathbb{R}^n$ between $x\in\R^n$ and $y\in\R^n$, 
the quadratic lower compensated convex transform
(lower transform for short) of $g$ for $\lambda>c_1$ is defined for $x\in\R^n$ by \cite{Zha08a},

\begin{equation}\label{Eq.01.Def.LwTr}
	C_{\lambda}^l(g)(x)=\co[g+\lambda|\cdot|^2](x)-\lambda|x|^2\,,
\end{equation}
where $\co[h]$ is the convex envelope of the function $h:\mathbb{R}^n\rightarrow\mathbb{R}\cup\{\infty\}$ bounded below.
Given $g:\mathbb{R}^n\rightarrow\mathbb{R}\cup\{-\infty\}$ such that $g(x)\leq c_0+c_1|x|^2$, $x\in\mathbb{R}^n$,
the quadratic upper compensated convex transform
(upper transform for short) of $g$ of module $\lambda>0$ is defined for $x\in\R^n$ by

\begin{equation}\label{Eq.01.Def.UpTr}
	C_{\lambda}^u(g)(x)=\lambda|x|^2-\co[\lambda|\cdot|^2-g](x)\,.
\end{equation}
The lower and upper transform can be, in turn, characterized in
terms of the critical mixed Moreau envelopes as follows 
(see Appendix \ref{Sec.Proofs})

\begin{equation}\label{Eq.01.MixMor}
	C^l_\lambda(g)(x)=M^\lambda(M_\lambda(g))(x)
	\quad\text{and}\quad
	C^u_\lambda(g)(x)=M_\lambda(M^\lambda(g))(x)\quad\text{for all }x\in\mathbb{R}^n
\end{equation}
where, in our notation,

\begin{equation}\label{Eq.01.Def.UpLwMor}
\begin{split}
	&M_{\lambda}(g)(x)=\inf_{y\in\mathbb{R}^n}\{g(y)+\lambda|y-x|^2 \},\\
	\text{and}\quad&
	M^{\lambda}(g)(x)=\sup_{y\in\mathbb{R}^n}\{g(y)-\lambda|y-x|^2\},
\end{split}
\end{equation}
are the lower and upper Moreau envelope of $g$, respectively \cite{Mor65,Mor66,LL86,AA93,CS04},
defined as the \textit{inf-} and \textit{sup-}convolution of $g$ with quadratic perturbations, respectively.

In mathematical morphology \cite{Ser82,Soi04,Shi09},
the Moreau lower and upper envelopes can be viewed as `greyscale' erosions and dilations  by
quadratic structuring elements \cite{BS92,BDMS96,Jac92}.  
Despite the rich structure of these filters, the use of quadratic structuring 
functions has not however been fully exploited compared, for instance, to the widespread use of the 
flat structuring elements \cite{Soi04,Shi09}. This is usually attributed to the fact that 
the gray--scale dilation and erosion operators based on a flat structuring function 
are easy to implement (given that they reduce to a local maximum and minimum filter) and
result therefore in fast algorithms. By the definition in terms of the compensated convex transforms 
it is possible however to offer not only an alternative
interpretation of the transforms \eqref{Eq.01.MixMor} as `one-step' morphological openings and closings, 
but also to provide an easy evaluation of the geometric properties of such filters \cite{ZOC15a}
and to obtain competetively fast algorithms as we illustrate in this paper.

In variational analysis, by contrast, the lower compensated convex transform \eqref{Eq.01.Def.LwTr} is also known 
as proximity hull \cite[Example 11.4]{RW98} and it has been redefined as quadratic envelope in \cite{Car19}. 
But it is only starting from 
\cite{Zha08a,ZOC15b} that the relation of the proximity hull to the convex envelope, and to the
envelope of quadratic functions with given curvature established first in \cite[Eq. (1.4)]{ZOC15a} and then in \cite{Car19}, 
has been studied in a systematic
manner and applied as such to image processing, shape interrogation, and function approximation in \cite{ZCO15c,ZCO18,OCZ20} 
and to compressive sensing and low rank approximation in \cite{Car19}.

Note that \eqref{Eq.01.Def.LwTr}, \eqref{Eq.01.Def.UpTr} and  their alternative representations \eqref{Eq.01.MixMor}
are given for extended real valued
functions. They can thus be applied to functions $g$, defined in proper subsets $\Omega$
of $\R^n$ by defining, as a common practice in convex analysis \cite{HL01}, 
the following extensions of $g$ in $\R^n\setminus \Omega$,

\begin{equation}\label{Eq:Intro:01}
	g^{\infty}:\,x\in\R^n\mapsto g(x)=\left\{
				\begin{array}{ll}
					g(x),&	x\in \Omega\\[1.5ex]
					+\infty,&	x\in\R^n\setminus\Omega\,,
				\end{array}
	\right.
\end{equation}
for the definition of $C_{\lambda}^l(g^{\infty})$, and

\begin{equation}\label{Eq:Intro:02}
	g^{-\infty}:\,x\in\R^n\mapsto g(x)=\left\{
				\begin{array}{ll}
					g(x),&	x\in \Omega\\[1.5ex]
					-\infty,&	x\in\R^n\setminus\Omega\,,
				\end{array}
	\right.
\end{equation}
for the definition of $C_{\lambda}^u(g^{-\infty})$.
However, such natural and direct definitions
of local compensated convex transforms for functions defined in $\Omega$  using \eqref{Eq.01.MixMor} 
depend on values of the Moreau envelope at points outside the domain $\Omega$, while convex envelope based 
methods using \eqref{Eq.01.Def.LwTr} or \eqref{Eq.01.Def.UpTr} will rely on 
an a-priori assumption about the boundary values of the transforms and on the application of 
convex envelope based schemes which, as far as we know, are neither efficient nor accurate.

These problems lead us to design the following simple local compensated convex transforms based on
the mixed Moreau envelope definitions \eqref{Eq.01.MixMor} of the compensated convex
transforms without the need of calculating values of the Moreau envelopes outside the bounded
closed domain $\Omega$. Before introducing our local transforms, we introduce some notation and recall
some definitions.

Let  $\Omega$ be a non-empty bounded open convex subset of $\mathbb{R}^n$.
We consider bounded functions $f:\Omega\subset\mathbb{R}^n\to \mathbb{R}$
satisfying $m\leq f(x)\leq M$ in $\Omega$ for some constants $-\infty<m\leq M< +\infty$.
Without loss of generality, and if not otherwise specified,  we set $m=\inf_{\Omega}f$ and
$M=\sup_{\Omega}f$ and define the oscillation of $f$ in $\Omega$ by 

\[
	O_f:=M-m\geq 0.
\]
We consider the auxiliary functions, $f^{-}_{\overline \Omega}$ and  $f^{+}_{\overline \Omega}$
that extend $f$ from $\Omega$ to its closure $\overline \Omega$,

\begin{subequations}\label{Eq.01.Ext.Local}
\begin{align}
	& f_{\overline \Omega}^{-}:\,x\in\overline\Omega\mapsto\, f_{\overline \Omega}^{-}(x)=
				\left\{\begin{array}{ll}
					\displaystyle f(x),	& \displaystyle  x\in \Omega,\\[1.5ex]
					\displaystyle \inf_{\Omega}f,  & \displaystyle  x\in \partial\Omega\,,
				\end{array}\right. \label{Eq.01.Ext.Local.LW} \\[1.5ex]
	\text{and}\quad &f_{\overline \Omega}^{+}:\,x\in\overline\Omega\mapsto\, f_{\overline \Omega}^{+}(x)=
				\left\{\begin{array}{ll}
					\displaystyle f(x),	& \displaystyle  x\in \Omega,\\[1.5ex]
					\displaystyle \sup_{\Omega}f,  & \displaystyle  x\in \partial\Omega\,,
				\end{array}\right. \label{Eq.01.Ext.Local.UP}
\end{align}
\end{subequations}
and the auxiliary functions $f_{\R^n}^{-}$ and $f_{\R^n}^{+}$ that extend $f$ from $\Omega$ to the whole space $\mathbb{R}^n$,

\begin{subequations}\label{Eq.01.Ext.Glob}
\begin{align}
	& f_{\R^n}^{-}:\,x\in\mathbb{R}^n\mapsto\, f_{\R^n}^{-}(x)=
				\left\{\begin{array}{ll}
					\displaystyle f(x),	& \displaystyle  x\in \Omega,\\[1.5ex]
					\displaystyle  \inf_{\Omega}f,& \displaystyle  x\in \mathbb{R}^n\setminus \Omega\,,
				\end{array}\right.  \label{Eq.01.Ext.Glob.LW}\\[1.5ex]
\text{and}\quad & f_{\R^n}^{+}:\,x\in\mathbb{R}^n\mapsto\, f_{\R^n}^{+}(x)=
				\left\{\begin{array}{ll}
					\displaystyle f(x),	& \displaystyle  x\in \Omega,\\[1.5ex]
					\displaystyle  \sup_{\Omega}f,& \displaystyle  x\in \mathbb{R}^n\setminus \Omega\,.
				\end{array}\right. \label{Eq.01.Ext.Glob.UP}
\end{align}
\end{subequations}

In practice, the extensions $f_{\overline \Omega}^{-}$, $f_{\overline \Omega}^{+}$ of $f$ to the boundary $\partial\Omega$ 
correspond to adding a frame of one pixel wide layer on the boundary of the data array and defining $f$ at each point 
of the  frame by the maximum value or the minimum value of the function.

For the function $f_{\R^n}^{-}$ and $f_{\R^n}^{+}$ we then consider the following transformations
which are well defined for $x\in\R^n$,

\begin{equation}\label{Eq.Def.LwtrUpTr}
	\begin{split}
		C_{\lambda}^l(f_{\R^n}^{-})(x)
				&=\co[f_{\R^n}^{-}+\lambda|\cdot|^2](x)-\lambda|x|^2,\\[1.5ex]
		C_{\lambda}^u(f_{\R^n}^{+})(x)
				&=\lambda|x|^2-\co[\lambda|\cdot|^2-f_{\R^n}^{+}](x)\,,
	\end{split}
\end{equation}
and their characterization  \eqref{Eq.01.MixMor} in terms of the critical mixed Moreau envelopes:

\begin{equation}\label{Eq.DefCmpTr.Mor}
	C^l_\lambda(f_{\R^n}^{-})(x)=M^\lambda(M_\lambda(f_{\R^n}^{-}))(x),\qquad
	C^u_\lambda(f_{\R^n}^{+})(x)=M_\lambda(M^\lambda(f_{\R^n}^{+}))(x)\,.
\end{equation}

For a bounded function $g: \overline\Omega\to \mathbb{R}$ defined in the closure of a convex bounded open set $\Omega$,
we also introduce for $x\in \overline \Omega$ the notation $M_{\lambda,\Omega}(g)(x)$ and $M^{\lambda}_{\Omega}(g)(x)$
to denote the following \textit{inf-} and  \textit{sup-}convolutions of $g$ with quadratic weights,

\begin{subequations}\label{Eq.01.Def.LocMrEnv}
\begin{align}
	&M_{\lambda,\Omega}(g)(x)=\underset{y\in \overline\Omega}{\inf}\{ g(y)+\lambda|y-x|^2\}\,,
		\label{Eq.01.Def.LocMrEnv.Lw}\\[1.5ex]
		\text{and}\quad
	&M^{\lambda}_{\Omega}(g)(x)=\underset{y\in \overline\Omega}{\sup}\{ g(y)-\lambda|y-x|^2\}\,,
		\label{Eq.01.Def.LocMrEnv.Up}
\end{align}
\end{subequations}
which coincide with the restriction to $\overline\Omega$
of $M_{\lambda}(g^{\infty})$ and $M^{\lambda}(g^{-\infty})$, respectively.

Motivated by the characterization \eqref{Eq.01.MixMor} for compensated convex transforms, we now define \textbf{the local lower compensated convex transform}
of $f$ in $\overline\Omega$, as

\begin{equation}\label{Def.1.LLCT}
	C^l_{\lambda,\Omega}(f_{\overline \Omega}^{-})(x):=M^{\lambda}_{\Omega}(M_{\lambda,\Omega}(f^{-}_{\overline \Omega}))(x)
	\quad\text{for }x\in\overline \Omega\,,
\end{equation}
and \textbf{the local upper compensated convex transform} of $f$ in $\overline\Omega$, as

\begin{equation}\label{Def.1.LUCT}
	C^u_{\lambda,\Omega}(f_{\overline \Omega}^{+})(x):=M_{\lambda,\Omega}(M^{\lambda}_{\Omega}(f^{+}_{\overline \Omega}))(x)
	\quad\text{for }x\in\overline \Omega\,.
\end{equation}

An important feature of the special extensions $f_{\overline \Omega}^{-}$, $f_{\overline \Omega}^{+}$, $f_{\R^n}^{-}$, $f_{\R^n}^{+}$ of
functions to $\overline\Omega$ and to $\mathbb{R}^n$ is that we have, for $x\in\Omega$ that

\begin{equation}\label{Sec1.Eq.Ident}
	C^l_{\lambda,\,\Omega}(f_{\overline\Omega}^{-})(x)=C^l_{\lambda}(f_{\R^n}^-)(x)\quad\text{and}\quad
	C^u_{\lambda,\,\Omega}(f_{\overline\Omega}^{+})(x)=C^l_{\lambda}(f_{\R^n}^+)(x)\,.
\end{equation}

The equalities \eqref{Sec1.Eq.Ident} do not hold, in general, if we had considered the extensions  
$f^{\infty}_{\overline\Omega}$ and $f^{\infty}$  [respect. $f^{-\infty}_{\overline\Omega}$ and $f^{-\infty}$] with 
$f^{\infty}_{\overline\Omega}$ [respect. $f^{-\infty}_{\overline\Omega}$] that extends $f$ over $\overline\Omega$ by 
$\infty$ [respect. $-\infty$] on $\partial \Omega$. 
While our extensions can make the local versions equals to the global ones, the infinity 
versions would not have this property.

After this brief introduction, in the next Section we present background results on convex analysis and the theory of compensated
convex transforms 
and we will also compare the definitions of the compensated convex transforms for 
the different extensions of $f$ over $\R^n$
by showing that for any $x\in \R^n$, there holds (Proposition \ref{Sec2:Prop:CmprCmpTr}) 

\begin{equation*}\label{Sec1:Eq:Comparison}
	C_{\lambda,\,\Omega}^l(f_{\overline\Omega}^{-})(x)\leq C_{\lambda,\,\Omega}^l(f)(x)
	\leq C_{\lambda}^l(f^{\infty})(x)\leq f(x)\leq C_{\lambda}^u(f^{-\infty})(x)\leq 
	C_{\lambda,\,\Omega}^u(f)(x)\leq C_{\lambda,\,\Omega}^u(f_{\overline\Omega}^{+})(x)
	\,,
\end{equation*}
where $C_{\lambda,\Omega}^l(f)$ and $C_{\lambda,\Omega}^u(f)$ are defined, respectively, for $x\in\overline\Omega$ as

\[
	C_{\lambda,\Omega}^l(f)(x):=M_{\Omega}^{\lambda}(M_{\lambda,\Omega}(f))(x)\quad\text{and}
	\quad C_{\lambda,\Omega}^u(f)(x):=M_{\lambda,\Omega}(M_{\Omega}^{\lambda}(f))(x)
\]
by means of \eqref{Eq.01.Def.LocMrEnv.Lw} and \eqref{Eq.01.Def.LocMrEnv.Up}.
The main theoretical results are stated in Section \ref{Sec.MainRslt}. Here, we state that
if $\Omega\subset\R^n$ is a bounded convex open set and $f:\Omega\to\R$ a bounded 
function, then for any $x\in\overline\Omega$ (see Theorem \ref{Theo.LocLwTr})

\begin{equation}
	C^l_{\lambda,\,\Omega}(f_{\overline\Omega}^{-})(x)=C^l_{\lambda}(f_{\R^n}^-)(x)\quad\text{and}\quad
	C^u_{\lambda,\,\Omega}(f_{\overline\Omega}^{+})(x)=C^l_{\lambda}(f_{\R^n}^+)(x)\,.
\end{equation}
Furthermore, we also show that there exists a constant $\kappa(\lambda,\,f)$ which depends on $\lambda$ and  $f$ such that
at points $x\in\Omega$ with $\dist(x,\partial\Omega)>\kappa(\lambda,\,f)$, the values of
$C^l_{\lambda,\,\Omega}(f_{\overline\Omega}^{-})(x)$ and $C^u_{\lambda,\,\Omega}(f_{\overline\Omega}^{+})(x)$
depend only on the values of $f$ on $\Omega$.
For the special case of the characteristic function $\chi_K$ of a compact set $K\subset\Omega$,
which represents `geometric shapes' (binary data),
we use the natural and simple extension $\chi^{\Omega}_K$ defined in \eqref{Sec3.Eq.06}, which  
is the restriction of the characteristic function
$\chi_K$ to $\overline{\Omega}$, and $\chi_K$ itself, rather than the extensions
defined in \eqref{Eq.01.Ext.Local} and \eqref{Eq.01.Ext.Glob}, respectively. 
Under the condition that $\dist^2(x,\partial\Omega)>1/\lambda$, we establish 
$C^u_{\lambda,\Omega}(\chi^{\Omega}_K)(x)=C^u_\lambda(\chi_K)(x)$ 
(see also the remarks about Theorem \ref{Thm.3.ChrctFnct} below).
We present in Section \ref{Sec.NumSchm} an algorithm that allows the numerical realization of the Moreau envelope,
whereas Section \ref{Sec.NumExprm} contains numerical experiments which illustrate how to apply our theoretical findings
to carry out, for instance, image processing and computational geometry tasks. As examples, we discuss the
finding of ridges in the graph of a function (multiscale medial axis map)
and intersections of curves in a plane, image inpainting and salt \& pepper denoising. 
Section \ref{Sec.Conc} concludes the paper with some final remarks whereas the Appendix 
contains the proofs of the main results.


\setcounter{equation}{0}
\section{Notation and Preliminaries}\label{Sec.Notat}

This section presents a brief overview of some basic results in convex analysis and  the theory
of compensated convex transforms that will be used in the sequel for the proof of the main results;
for a comprehensive account of convex analysis, we refer to Refs.~\cite{HL01,Roc70}, and 
to Refs.~\cite{Zha08a,ZOC15a}
for an account of the theory of compensated convex transforms.

\begin{prop}\label{Prp.GlCnvx}
	Let $f:\mathbb{R}^n\to (-\infty,\,+\infty]$ be coercive in the sense that
	$f(x)/|x|\to\infty$ as $|x|\to \infty$, and $x_0\in\R^n$. 
	Denote by $\Aff(\R^n,\R)$ the set 
	of real--valued affine functions on $\R^n$. Then
   \begin{itemize}
	\item[(i)]
		The value $\co\left[f\right](x_0)$ of the convex envelope of $f$ at
		$x_0\in \mathbb{R}^n$ is given by

		\begin{equation}\label{Eq.Prp.GlCnvx.1}
			\begin{array}{ll}
			\displaystyle  \co\left[f\right](x_0)=\underset{x_1,\ldots,x_{n+1}}{\inf}\,\Bigg\{ \sum^{n+1}_{i=1}\lambda_if(x_i):\;  
			        &\displaystyle\sum^{n+1}_{i=1}\lambda_i=1,\; \sum^{n+1}_{i=1}\lambda_ix_i=x_0,\\[1.5ex]
				&\displaystyle 	\lambda_i\geq 0,\; x_i\in \R^n
			\Bigg\}\,.
		\end{array}
		\end{equation}
		If, in addition, $f$ is lower semicontinuous, then 
		the infimum is attained, that is, for some $1\leq p\leq n+1$ there are
		$\lambda_k^\ast>0$, $y_k^\ast\in\mathbb{R}^n$ for $k=1,\dots,p$, satisfying
		$\sum^p_{k=1}\lambda^\ast_k=1$ and $\sum^p_{k=1}\lambda^{\ast}_ky^{\ast}_k=x_0$ such that
		the points $\left(y_k^{\ast}, f(y_k^{\ast})\right)$, $k=1,\ldots,\,p$, lie in the intersection of a supporting plane of
		the epigraph of $f$, $\epi(f)$, and $\epi(f)$, and 
		
		\begin{equation}
			\co [f](x_0)=\sum^p_{k=1}\lambda^{\ast}_kf(y^{\ast}_k)\,.
		\end{equation}
	\item[(ii)]
		The value $\co\left[f\right](x_0)$, for $f$ taking only finite values, can also be obtained as follows:
	
		\begin{equation}\label{Eq.Prp.GlCnvx.2}
			\co\left[f\right](x_0)=\sup\left\{
				\ell(x_0):\; \ell\;\;\text{\rm affine}\quad\text{\rm and}\quad
				\ell(y)\leq f(y)\;\;\text{\rm for all }y\in\R^n
				\right \}
		\end{equation}
		with the $\sup$ attained by an affine function $\ell^\ast\in\Aff(\R^n)$ if $f$ is lower semicontinuous.
     \end{itemize}
\end{prop}

We will also introduce the following local version of convex envelope at a point.

\begin{defi}  \label{Def.LocCnvx}
	Let $r>0$ and $x_0\in\R^n$. Denote by $B(x_0;r)$ the open ball centered at $x_0$ with radius $r$,
	and by $\overline{B}(x_0;r)$ the corresponding closed ball. 
	Suppose $f:\overline B(x_0;r)\to \mathbb{R}$ is a bounded  function in $\overline B(x_0;r)$. Then
	the value $\co_{\bar{B}(x_0;r)}\left[f\right](x_0)$ of the local convex envelope of $f$ at $x_0$ in
	$\overline B(x_0;r)$ is defined by

	\begin{equation*}\label{Eq.Def.LocCnvx.1}
	\begin{array}{ll}
			\displaystyle  \co_{\bar{B}(x_0;r)}\left[f\right](x_0)
			=\underset{x_1,\ldots,x_{n+1}}{\inf}\,\Bigg\{\sum^{n+1}_{i=1}\lambda_i f(x_i): 
				&\displaystyle	\sum^{n+1}_{i=1}\lambda_i=1,\;\sum^{n+1}_{i=1}\lambda_ix_i=x_0,\\[1.5ex]
				&\displaystyle	\lambda_i\geq 0,\;|x_i-x_0|\leq r,\; x_i\in\R^n		
			\Bigg\}.
	\end{array}
	\end{equation*}
\end{defi}

\begin{nota}\label{Rem.LocCnvx}
	If $f$ is lower semicontinuous, then by using the second part of Proposition \ref{Prp.GlCnvx}$(i)$,
	we see that the infimum is attained in $\overline{B}(x_0;r)$. This means that 
	for some $1\leq p\leq n+1$ there are $\lambda_k>0$, $y_k\in\mathbb{R}^n$ such that $|y_k-x_0|\leq r$  for $k=1,\dots,p$, satisfying
	$\sum^p_{k=1}\lambda_k=1$, $\sum^p_{k=1}\lambda_ky_k=x_0$ and
	
	\begin{equation}
		\co_{\bar B(x_0;\,r)}[f](x_0)=\sum^p_{k=1}\lambda_kf(y_k)\,,
	\end{equation}
	thus, in this case, $\co[f](x_0)$ depends only on the values of $f$ in $\overline B(x_0;\,r)$.
\end{nota}

We recall also the following ordering properties for compensated convex transforms which can be found 
in Ref.~\cite{Zha08a},

\begin{equation}\label{Sec2.Eq.Ordering}
	C^l_\lambda(f)(x)\leq f(x)\leq C^u_\lambda(f)(x),\quad x\in\mathbb{R}^n\,,
\end{equation}
whereas for $f\leq g$ in $\mathbb{R}^n$, we have that

\begin{equation}\label{Sec2.Eq.Ineq}
	C^l_\lambda(f)(x)\leq C^l_\lambda(g)(x)\quad\text{and}\quad
	C^u_\lambda(f)(x)\leq C^u_\lambda(g)(x),\quad x\in\mathbb{R}^n\,.
\end{equation}

\begin{prop}\label{Prop.Tra} (Translation invariance property)
	For any $f:\mathbb{R}^n\to (-\infty,\,+\infty]$
	bounded below and for any affine function $\ell:\mathbb{R}^n\to \mathbb{R}$,
	$\co [f+\ell]=\co[f]+\ell$. Consequently,
	both $C^u_\lambda(f)$ and $C^l_\lambda(f)$ are translation invariant against the weight function, that is

	\begin{equation*}\label{Eq.Prop.Tra.1}
	\begin{array}{l}
		\displaystyle C^l_\lambda(f)(x)=\co\left[\lambda |(\cdot)-x_0|^2+f\right](x)-\lambda |x-x_0|^2\,,\\[1.5ex]
		\displaystyle C^u_\lambda(f)(x)=\lambda |x-x_0|^2-\co\left[\lambda |(\cdot)-x_0|^2-f\right](x)
	\end{array}
	\end{equation*}
	for all $x\in \R^n$ and for every fixed $x_0$. In particular, at $x_0$, we have
	
	\begin{equation*}\label{Eq.Prop.Tra.2}
			C^l_\lambda(f)(x_0)=\co[\lambda  |(\cdot)-x_0|^2+f](x_0)\,, \quad
			C^u_\lambda(f)(x_0)=-\co[\lambda |(\cdot)-x_0|^2-f](x_0)\,.
	\end{equation*}
\end{prop}

For both theoretical and numerical developments,
the following property on the locality of the compensated  convex transforms for Lipschitz 
functions and bounded functions plays a fundamental role. 
The result for bounded functions is a slight modification of the locality property stated
in Theorem 3.10, Ref.~\cite{ZOC15a}. For a locally bounded function 
$f:\mathbb{R}^n\to \mathbb{R}$, we define the upper and the lower semicontinuous closure $\overline{f}$ 
and $\underline{f}$\cite{HL01,Roc70}, respectively, by

\begin{equation}\label{Sec2.Eq.SemClos}
	\overline{f}(x)=\underset{y\to x}{\limsup}\,f(y)\quad\text{\rm and}\quad
	\underline{f}(x)=\underset{y\to x}{\liminf}f(y).
\end{equation}
We have  the following result.

\begin{prop}\label{Prop.Loc}
Suppose $f:\mathbb{R}^n\rightarrow \mathbb{R}$ is bounded. Let $\lambda>0$ and $x_0\in\mathbb{R}^n$. 
Then the  following locality properties hold:

\begin{equation}\label{Sec2.Eq.Loc}
	\begin{split}
		C^l_\lambda(f)(x_0)&=\co_{B(x_0;R_{f,\lambda})}[\underline{f}+\lambda|(\cdot)-x_0|^2](x_0)\\[1.5ex] 
		C^u_\lambda(f)(x_0)&=-\co_{B(x_0;R_{f,\lambda})}[\lambda|(\cdot)-x_0|^2-\overline{f}](x_0)\\[1.5ex] 			 
	\end{split}
\end{equation}
with $R_{f,\lambda}=3\sqrt{O_f/\lambda}$. If $f$ is Lipschitz continuous with Lipschitz constant $L\geq 0$, 
then $R_{f,\lambda}=4L/\lambda$. 
\end{prop}

\begin{nota}
The values of $R_{f,\lambda}$ given here have improved upon those obtained in 
Theorem 3.10, Ref.~\cite{ZOC15a}.
\end{nota}

By Remark \ref{Rem.LocCnvx}, the convex envelope over $B(x_0;R_{f,\lambda})$ which enters \eqref{Sec2.Eq.Loc} is therefore given by 

\begin{equation}\label{Sec2.Eq.LocCnvx}
	\co_{B(x_0;R_{f,\lambda})}[\underline{f}+\lambda|(\cdot)-x_0|^2](x_0)=\sum^p_{k=1}\lambda_k(\lambda|y_k-x_0|^2+\underline{f}(y_k))
\end{equation}
for some $1\leq p\leq n+1$, $\lambda_k>0$, $y_k\in\mathbb{R}^n$ for $k=1,\dots,p$, satisfying
$\sum^p_{k=1}\lambda_k=1$, $\sum^p_{k=1}\lambda_ky_k=x_0$ and

\begin{equation}\label{Sec2.Eq.LocEst}
	|y_k-x_0|\leq R_{f,\lambda}
\end{equation}
for all $k=1,\dots,p$. Similar conclusion can be drawn for $\co_{B(x_0;R_{f,\lambda})}[\lambda|(\cdot)-x_0|^2-\overline{f}](x_0)$.

Next we state the locality properties for the Moreau envelopes.

\begin{prop}\label{Prop.3.Mor}
Let $f:\mathbb{R}^n\to \mathbb{R}$.
\begin{itemize}
\item[$(i)$] If $f$ is bounded, then for any fixed $x\in\mathbb{R}^n$, if $y_x,\,z_x\in\mathbb{R}^n$ satisfy

	\begin{equation}\label{Sec2.Eq.LocMor}
		\begin{split}
			M_\lambda(f)(x)&=\underset{y\in\mathbb{R}^n}{\inf}\{\underline{f}+\lambda|y-x|^2\}=
			\underline{f}(y_x)+\lambda|y_x-x|^2,\\
			M^\lambda(f)(x)&=\underset{y\in\mathbb{R}^n}{\sup}\{\overline{f}(y)-\lambda|y-x|^2\}=
			\overline{f}(z_x)-\lambda|z_x-x|^2,
		\end{split}
	\end{equation}
	then 

	\begin{equation}\label{Sec2.Eq.LocMorEst}
		|y_x-x|\leq \sqrt{O_f/\lambda}\text{ and }|z_x-x|\leq \sqrt{O_f/\lambda}.
	\end{equation}

\item[$(ii)$] If $f$ is a Lipschitz function with Lipschitz constant $L\geq 0$.
	 Then for any fixed $x\in\mathbb{R}^n$, if $y_x,\,z_x\in\mathbb{R}^n$ satisfy

	\begin{equation}
		\begin{split}
		M_\lambda(f)(x)&=\inf\{ f(y)+\lambda|y-x|^2,\; y\in\mathbb{R}^n\}=
		 f(y_x)+\lambda|y_x-x|^2,\\
		M^\lambda(f)(x)&=\sup\{ f(y)-\lambda|y-x|^2,\; y\in\mathbb{R}^n\}=
		 f(z_x)-\lambda|z_x-x|^2,
		\end{split}
	 \end{equation}
	then 

	\begin{equation}
		|y_x-x|\leq L/\lambda\text{ and }|z_x-x|\leq L/\lambda\,.
	\end{equation}
\end{itemize}
\end{prop}

\begin{nota}
Proposition \ref{Prop.3.Mor}$(i)$ was also established in \cite[Lemma 3.5.7]{CS04} whereas
\ref{Prop.3.Mor}$(ii)$ is partially contained in \cite[Theorem 4.2]{CLSW98}.
In Appendix \ref{Sec.Proofs} we give a new the proof of these results.
\end{nota}

For completeness, we recall the following relationship between the lower and 
upper Moreau envelope and the lower and upper compensated convex transform given by 

\begin{equation}\label{Eq.RelMor}
	M_{\lambda}(f)=-M^{\lambda}(-f)\quad\text{and}\quad
	C_{\lambda}^l(f)=-C_{\lambda}^u(-f)\,,
\end{equation}
and the ordering properties of the Moreau envelope which we give next for the lower Moreau envelope
in the case of $f\leq g$ with $f,\,g$ bounded from below, 

\begin{equation}\label{Eq:OrderMoreau}
	M_{\lambda}(f)(x)\leq f(x)\quad\text{and}\quad M_{\lambda}(f)(x)\leq M_{\lambda}(g)(x)\quad\text{for }x\in\R^n\,.
\end{equation}

The following result precises the connection among the compensated convex transforms for the different extensions
of $f:\Omega\subseteq\R^n\to \R$ over $\R^n$. We give this relationship only for the lower compensated 
convex transform given that by \eqref{Eq.RelMor} we can establish similar ones
for the upper transform.

\begin{prop}\label{Sec2:Prop:CmprCmpTr}
	Let $\Omega\subset \R^n$ be a convex open set and assume $f:\Omega\subset \R^n\rightarrow \R$ to be bounded. 
	Denote by $\underline{f}:\overline\Omega\subset \R^n\rightarrow \R$ the lower semicontinuous envelope of $f$
	defined by \eqref{Sec2.Eq.SemClos} at the interior points $x$ and by

	\begin{equation}\label{Sec2.Eq.SemClosBnd}
		\underline{f}(x)=\liminf_{y\to x,\, y\in\Omega}f(y)		
	\end{equation}
	at the boundary points $x\in\partial \Omega$. The following inequalities hold

	\begin{equation}
		C_{\lambda,\,\Omega}^l(f_{\overline\Omega}^{-})(x)\leq C_{\lambda,\,\Omega}^l(f)(x)
		\leq C_{\lambda}^l(f^{\infty})(x)\leq \underline{f}(x) \leq f(x)	
	\end{equation}
	for any $x\in \R^n$.
\end{prop}

Next we recall from \cite{DL93} the definition of modulus of continuity of a function along 
with some of its properties which will be used to analyse the numerical scheme that computes 
the Moreau envelope. 

\begin{defi}\label{Sec2.Def.MoC}
Let $f:\R^n \to \R$ be a bounded and uniformly continuous function in $\R^n$. Then,

\begin{equation}\label{Sec2.Eq.Def.MoC}
	\omega_f:t\in[0,\,\infty)\mapsto \omega_{f}(t)=\sup\Big\{|f(x)-f(y)|:\,x,y\in\R^n\text{ and }|x-y|\leq t\Big\}
\end{equation}
is called the modulus of continuity of $f$.
\end{defi}

\begin{prop}\label{Sec2.Pro.MoC}
Let $f:\R^n \to \R$ be a bounded and uniformly continuous function in $\R^n$. Then the modulus of continuity $\omega_f$ of
$f$ satisfies the following properties:

\begin{equation}\label{Sec2.Eq.Pro.MoC}
	\begin{array}{ll}
		$(i)$	&\omega_f(t)\to \omega(0)=0,\text{ as }t\to 0;	\\[1.5ex]
		$(ii)$	&\omega_f \text{ is non-negative and non-decreasing continuous function on }[0,\infty);	\\[1.5ex]
		$(iii)$	&\omega_f \text{ is subadditive: }\omega_f(t_1+t_2)\leq\omega_f(t_1)+\omega_f(t_2)
			\text{ for all }t_1,\,t_2\geq 0\,.
	\end{array}
\end{equation}
\end{prop}
A function $\omega$ defined on $[0,\,\infty)$ satisfying \eqref{Sec2.Eq.Pro.MoC} is called a
modulus of continuity. A modulus of continuity
$\omega$ can be bounded from above by an affine function (see Lemma 6.1 of Ref.~\cite{DL93}),
that is, there exist some constants
$a>0$ and $b\geq 0$ such that

\begin{equation}\label{Sec2.Eq.Pro.MoC.0}
	\omega(t)\leq at+b\quad(\text{for all }t\geq 0).
\end{equation}

We conclude this Section by recalling the following definitions.

Let $C$ be a subset of $\R^n$. We define the distance of $x\in\R^n$ from $C$ as

\begin{equation}
	\dist(x,C):=\inf_{y\in C}|y-x|\,,
\end{equation}
the diameter of the set $C$ as

\begin{equation}
	\mathrm{diam}(C):=\sup_{x,y\in C}|y-x|\,,
\end{equation}
the indicator function $i_{C}$ of $C\subset\R^n$ as the function defined in $\R^n$ such that

\begin{equation}
	i_C(x):=\left\{\begin{array}{ll}
			0,	&\text{if }x\in C;\\[1.5ex]
			+\infty,&\text{otherwise}\,,
		\end{array}\right.
\end{equation}
and the characteristic function $\chi_{C}$ of $C\subset\R^n$ as the function defined in $\R^n$ such that

\begin{equation}
	\chi_C(x):=\left\{\begin{array}{ll}
			1,	&\text{if }x\in C;\\[1.5ex]
			0,&\text{otherwise}\,.
		\end{array}\right.
\end{equation}

It is then not difficult to verify that for any $x\in\R^n$

\begin{equation}\label{Sec2.Eq.EucDistMor}
	\dist^2(x,\,C)=\underset{y\in\R^n}{\inf}\left\{i_C(y)+|y-x|^2 \right\}=\frac{1}{\lambda}M_{\lambda}(\lambda i_C)(x)\,.
\end{equation}
that is, $\dist^2(\cdot,\,C)$ is the \textit{inf}-convolution of $i_C$ with $\|\cdot\|^2$, and is proportional to
the lower Moreau envelope of $\lambda i_C$ with parameter $\lambda$.


\setcounter{equation}{0}
\section{Main Results}\label{Sec.MainRslt}
The main results given in this section consist of two parts. In the first part, we establish the 
relationship between the local Moreau envelopes 
of \eqref{Eq.01.Ext.Local} and the global Moreau envelopes of \eqref{Eq.01.Ext.Glob}, and between   
the corresponding mixed Moreau envelopes. This relationship is a consequence of the type of auxiliary 
functions under consideration. In the second part, we give conditions that ensure that the local Moreau envelopes 
and the corresponding mixed Moreau envelopes depend only on the local values of $f$. The precise meaning of this statement will be
specified for each result. Next we consider the case of a bounded function $f$ defined on a bounded domain $\Omega$.

\begin{teo}\label{Theo.LocLwTr}
Let $\Omega\subset\R^n$  be a bounded open set
and $f:\Omega\subset\mathbb{R}^n\to \mathbb{R}$ a bounded function. Consider the extensions
$f_{\overline\Omega}^{-}$ and $f_{\mathbb{R}^n}^{-}$ given by \eqref{Eq.01.Ext.Local.LW} and \eqref{Eq.01.Ext.Glob.LW}, respectively,
and the extensions $f_{\overline\Omega}^{+}$ and $f_{\mathbb{R}^n}^{+}$ given by \eqref{Eq.01.Ext.Local.UP}
and \eqref{Eq.01.Ext.Glob.UP}, respectively. 
Then, for any $x\in\overline\Omega$,  

\begin{align} 
	M_{\lambda,\Omega}(f_{\overline\Omega}^{-})(x)&=M_\lambda(f_{\R^n}^{-})(x)\,, \label{Sec3.Eq.01}\\[1.5ex]
	M^{\lambda}_{\Omega}(f_{\overline\Omega}^{+})(x)&=M^{\lambda}(f_{\R^n}^{+})(x)\,,
\end{align}
and 

\begin{align} 
	M^{\lambda}_{\Omega}(M_{\lambda,\Omega}(f_{\overline\Omega}^{-}))(x)&=M^{\lambda}(M_{\lambda}(f_{\R^n}^{-}))(x)\,,  \label{Sec3.Eq.02}\\[1.5ex]
	M_{\lambda,\Omega}(M^{\lambda}_{\Omega}(f_{\overline\Omega}^{+}))(x)&=M_{\lambda}(M^{\lambda}(f_{\R^n}^{+}))(x)\,.
\end{align}
Consequently, for any $x\in\overline\Omega$,  

\begin{align}\label{Sec3.Eq.03} 
	C^l_{\lambda,\Omega}(f_{\overline\Omega}^{-})(x)&=C^l_\lambda(f_{\R^n}^{-})(x)\,,\\[1.5ex]
	C^u_{\lambda,\Omega}(f_{\overline\Omega}^{+})(x)&=C^u_\lambda(f_{\R^n}^{+})(x)\,.
	\label{Sec3.Eq.03a} 
\end{align}
Furthermore, we have the following locality results:

\begin{itemize}
	\item[$(i)$] If $x\in\Omega$ is such that $\dist^2(x,\partial\Omega)>O_f/\lambda$ 
		and there is $z_x\in\mathbb{R}^n$ such that 

		\begin{equation*}
		\begin{split}
			&M_{\lambda,\Omega}(f_{\overline\Omega}^-)(x)=f_{\overline\Omega}^-(z_x)+\lambda|z_x-x|^2\\[1.5ex]
			&[\text{resp. }M^{\lambda}_{\Omega}(f_{\overline\Omega}^+)(x)=f_{\overline\Omega}^+(z_x)-\lambda|z_x-x|^2]\,,
		\end{split}
		\end{equation*}
		then  $M_{\lambda,\Omega}(f_{\overline\Omega}^-)(x)$ 
		[resp. $M^{\lambda}_{\Omega}(f_{\overline\Omega}^+)(x)$] 
		is determined by values of $f$ on $\Omega$, in the sense that 
		$z_x\in \Omega$.
	\item[$(ii)$] If $x\in\Omega$ is such that $\dist^2(x,\partial\Omega)>4O_f/\lambda$ 
		and there is a $z_x\in\mathbb{R}^n$ such that 
		
		\begin{equation*}
		\begin{split}
			&M^{\lambda}_{\Omega}(M_{\lambda,\Omega}(f_{\overline\Omega}^-))(x)=
			    M_{\lambda,\Omega}(f_{\overline\Omega}^-)(z_x)-\lambda|z_x-x|^2\\[1.5ex]
			&[\text{resp. }M^{\lambda}_{\Omega}(M_{\lambda,\Omega}(f_{\overline\Omega}^+))(x)=
				M_{\lambda,\Omega}(f_{\overline\Omega}^+)(z_x)+\lambda|z_x-x|^2]
		\end{split}
		\end{equation*}				
		then $M^{\lambda}_{\Omega}(M_{\lambda,\Omega}(f_{\overline\Omega}^-))(x)$ 
		[resp. $M^{\lambda}_{\Omega}(M_{\lambda,\Omega}(f_{\overline\Omega}^+))(x)$] 
		is determined by values of $f$ on $\Omega$, in the sense that 
		$M_{\lambda,\Omega}(f_{\overline\Omega}^-)(z_x) = f_{\overline\Omega}^-(y_x) + \lambda |y_x-x|^2$ where 
		$z_x,\,y_x \in \Omega$.
\end{itemize}
\end{teo}

\begin{nota}
The locality properties of Theorem \ref{Theo.LocLwTr} state that under the conditions $(i)$ and $(ii)$ on $x\in\Omega$, respectively,
the values of $C^l_{\lambda,\Omega}(f_{\overline\Omega}^{-})(x)$ [resp. $C^l_{\lambda,\Omega}(f_{\overline\Omega}^{+})(x)$] 
depend on the values of $f$ on $\Omega$.
This means that the values of $C^l_{\lambda,\Omega}(f_{\overline\Omega}^{-})(x)$ 
[resp. $C^l_{\lambda,\Omega}(f_{\overline\Omega}^{+})(x)$]
are not influenced by the values of $f$ we defined  on $\partial\Omega$
when we define $f_{\overline\Omega}^{-}(x)$ [resp. $f_{\overline\Omega}^{+}(x)$]. 
We express this by saying that $C^l_{\lambda,\Omega}(f_{\overline\Omega}^{-})(x)$ 
[resp. $C^l_{\lambda,\Omega}(f_{\overline\Omega}^{+})(x)$] 
is not affected by boundary values.
\end{nota}

As an application of Theorem \ref{Theo.LocLwTr} we consider the case where 
$f$ is the squared Euclidean distance function to a closed set. 
The following two results are useful, for instance, when we need to compute the multiscale medial axis map \cite{ZCO15c}. 
Let $K$ be a nonempty closed set, the quadratic multiscale medial axis map of $K$ with scale $\lambda>0$ is defined 
in \cite[Definition 3.1]{ZCO15c} for $x\in\R^n$ by

\[
	\mathcal{M}(\lambda;\,K)(x)=(1+\lambda)\Big(\dist^2(x;\,K)-C_{\lambda}^l(\dist^2(\cdot;\,K))(x)
		\Big)\,.
\]
Next we describe how  $\mathcal{M}(\lambda;\,K)(x)$ can be expressed in terms of the local lower transform.
The first result can be applied to find the multiscale medial axis map of the set $\Omega\setminus K$,
where $\Omega$ is an open subset of $\mathbb{R}^n$ and $K\subset \Omega$ a compact set.

\begin{coro}\label{Cor.SetBnd}
 	Let $\Omega\subset \mathbb{R}^n$ be an open set and $K\subset \Omega$ a non-empty compact set.
	Let $f(x):=\dist^2(x, K\cup \Omega^c)$ for $x\in \R^n$ and $f_{\overline\Omega}^{-}(x)$
	be defined by \eqref{Eq.01.Ext.Local.LW}. 
 	Then for $x\in\overline\Omega$,

	\begin{equation}\label{Sec3.Eq.04}
		\mathcal{M}(\lambda;\,K\cup\Omega^c)(x)=(1+\lambda)\Big(f_{\overline\Omega}^{-}(x)-
		C^l_{\lambda,\Omega}(f_{\overline\Omega}^{-})(x)\Big)\,.
	\end{equation}
\end{coro}
 
\begin{nota}\label{Sec3.Rem.MMA}
	Equation \eqref{Sec3.Eq.04} actually gives $\mathcal{M}(\lambda;\,K\cup\Omega^c)(x)$ for any $x\in\R^n$
	given that $\mathcal{M}(\lambda;\,K\cup\Omega^c)(x)=0$ for $x\in\R^n\setminus\overline\Omega$.
\end{nota}

The next result, on the other hand, applies when we need to define the multiscale medial axis map 
of an open set $A\subset \Omega$. 

\begin{coro}\label{Cor.SetCmp}
	Let $\Omega\subset \mathbb{R}^n$  be an open set and $A\subset \Omega$ a non-empty open set.
	Let $f(x):=\dist^2(x, A^c)$ and define  $f_{\overline\Omega}^{-}(x)$ by \eqref{Eq.01.Ext.Local.LW}. 
	Then for any $x\in\overline\Omega$ 
	
	\begin{equation}\label{Sec3.Eq.05}
		\mathcal{M}(\lambda;\,A^c)(x)=(1+\lambda)\Big(f_{\overline\Omega}^{-}(x)-
		C^l_{\lambda,\Omega}(f_{\overline\Omega}^{-})(x)\Big)\,.
	\end{equation}
 \end{coro}

\begin{nota}\label{Sec3.Rem.MMA2}
In this case we also have that
$\mathcal{M}(\lambda;\,A^c)(x)=0$ for $x\in\R^n\setminus\overline\Omega$.
\end{nota}

Next we consider the behavior of the local upper compensated transform of the characteristic function of a compact set 
in view of applications that involve the processing of binary images. 

\begin{teo}\label{Thm.3.ChrctFnct}
Let $\Omega\subset \mathbb{R}^n$ be an open set and $K\subset \Omega$ a non-empty compact set.
Let $\chi_K$ denote the characteristic function of $K$ defined in $\mathbb{R}^n$
and $\chi_K^{\Omega}$ the restriction of $\chi_K$ to $\overline{\Omega}$, that is,
	
\begin{equation}\label{Sec3.Eq.06}
		\chi_K^{\Omega}(x)=\left\{\begin{array}{ll}
				\displaystyle	1& \displaystyle \text{if }x\in K\\[1.5ex]
				\displaystyle	0& \displaystyle \text{if }x\in \overline\Omega\setminus K\,.
			\end{array}\right.
\end{equation}
Then, if $\dist^2(K,\,\partial\Omega)>1/\lambda$, for any $x\in\overline\Omega$ 

\begin{equation}\label{Sec3.Eq.07}
	M_{\Omega}^{\lambda}(\chi_K^{\Omega})(x)=M^{\lambda}(\chi_K)(x)
\end{equation}
and

\begin{equation}\label{Sec3.Eq.08}
	M_{\lambda,\Omega}(M^{\lambda}_{\Omega}(\chi_K^{\Omega}))(x)=M_{\lambda}(M^{\lambda}(\chi_K))(x)\,.
\end{equation}
Consequently, if we define $C_{\lambda,\Omega}^u(\chi_K^{\Omega})(x)=M_{\lambda,\Omega}(M^{\lambda}_{\Omega}(\chi_K^{\Omega}))(x)$, 
it follows that  

\begin{equation}\label{Sec3.Eq.09}
	C_{\lambda,\Omega}^u(\chi_K^{\Omega})(x)=C_{\lambda}^u(\chi_K)(x)\,.
\end{equation}
\end{teo}

\begin{nota}
\begin{itemize}
	\item[$(i)$] 
		Compared to Theorem \ref{Theo.LocLwTr}, Theorem \ref{Thm.3.ChrctFnct} states that in the case of binary functions
		we can establish the equalities \eqref{Sec3.Eq.07} and \eqref{Sec3.Eq.08} using $\chi_K^{\Omega}$
		as defined by \eqref{Sec3.Eq.06}, and $\chi_K$, rather than the corresponding 
		auxiliary functions \eqref{Eq.01.Ext.Local.UP} and \eqref{Eq.01.Ext.Glob.UP}, respectively,
		which are the type of functions that are used in Theorem \ref{Theo.LocLwTr}.
	\item[$(ii)$] It is possible to establish a locality result for $M_{\Omega}^{\lambda}(\chi_K^{\Omega})$ and
		$M_{\lambda,\Omega}(M_{\Omega}^{\lambda}(\chi_K^{\Omega}))$ 
		in the following sense. If $x\in\overline\Omega$ is such that 
		$\dist^2(x,\,\partial\Omega)>4/\lambda$ then $\dist^2(z_x,\,\partial\Omega)>1/\lambda$, where $z_x\in\overline\Omega$
		is such that 

		\begin{equation*}
			M_{\Omega}^{\lambda}(\chi_K^{\Omega})(x)=\chi_K^{\Omega}(z_x)-\lambda|z_x-x|^2\,.
		\end{equation*}
		Also, if $x\in\overline\Omega$, then 
		$\dist^2(z_x,\,\partial\Omega)>1/\lambda$ and $\dist^2(y_x,\,\partial\Omega)>1/\lambda$ where 
		$y_x,\,z_x\in\overline\Omega$ are such that 

		\begin{equation*}
			\begin{split}
				M_{\lambda,\Omega}(M_{\Omega}^{\lambda}(\chi_K^{\Omega}))&=M_{\Omega}^{\lambda}(\chi_K^{\Omega})(z_x)+\lambda|z_x-x|^2\\[1.5ex]
					&=\chi_K^{\Omega}(y_x)-\lambda|y_x-z_x|^2+\lambda|z_x-x|^2\,.
			\end{split}
		\end{equation*}
		As a consequence, if   
		$K=\{y\in\overline\Omega:\,\dist^2(y,\,\partial\Omega)\geq 1/\lambda\}$ and 
		$\dist^2(x,\,\partial\Omega)>9/\lambda$, then both 
		$M_{\Omega}^{\lambda}(\chi_K^{\Omega})(x)$ and
		$M_{\lambda,\Omega}(M_{\Omega}^{\lambda}(\chi_K^{\Omega}))(x)$ are determined only by 
		$K$, i.e. $y_x,\,z_x\in K$. The proof of these results follows from the locality 
		property of the Moreau envelopes (see Proposition \ref{Prop.3.Mor})
		and by a similar argument to the proof of  Theorem \ref{Theo.LocLwTr}$(i)$ and $(ii)$. 
\end{itemize}	
\end{nota}

For applications to scattered data approximation and image inpainting \cite{ZCO16a,ZCO18},
we assume $K\subset\Omega\subset \R^n$ to be a compact set, $M>0$ and $f:K\to\mathbb{R}$
to be a bounded function, and introduce the following auxiliary functions
 
\begin{subequations}\label{Def.01.AuxFnct.M}
\begin{align}
	&f_{\overline\Omega,K}^M:\,x\in\overline\Omega\mapsto f_{\overline\Omega,K}^M(x)=\left\{\begin{array}{ll}
 				f(x)&x\in K\\[1.5ex]
 				M&x\in \Omega\setminus K\\[1.5ex]
 				\inf_{K}f&x\in \partial\Omega
 				\end{array}\right.\label{Def.01.AuxFnct.M.Loc1}\\[2.5ex]
 	&f^M_{\mathbb{R}^n,K}:\,x\in\R^n\mapsto f^M_{\mathbb{R}^n,K}(x)=\left\{\begin{array}{ll}
 				f(x)&x\in K\\[1.5ex]
 				M&x\in \Omega\setminus K\\[1.5ex]
 				\inf_{K}f&x\in \mathbb{R}^n\setminus\Omega
 				\end{array}\right.\label{Def.01.AuxFnct.M.Glob1}\\[2.5ex]
	&f_{\overline\Omega,K}^{-M}:\,x\in\overline\Omega\mapsto f_{\overline\Omega,K}^{-M}(x)=\left\{\begin{array}{ll}
 				f(x)&x\in K\\[1.5ex]
 				-M&x\in \Omega\setminus K\\[1.5ex]
 				\sup_{K}f&x\in \partial\Omega
 				\end{array}\right.\label{Def.01.AuxFnct.M.Loc2}\\[2.5ex]
 	&f^{-M}_{\mathbb{R}^n,K}:\,x\in\R^n\mapsto f^{-M}_{\mathbb{R}^n,K}(x)=\left\{\begin{array}{ll}
 				f(x)&x\in K\\[1.5ex]
 				-M&x\in \Omega\setminus K\\[1.5ex]
 				\sup_{K}f&x\in \mathbb{R}^n\setminus\Omega
 				\end{array}\right.\label{Def.01.AuxFnct.M.Glob2}
\end{align}
\end{subequations}
where $\inf_{K}f$ and $\sup_{K}f$ denote the infimum and supremum of $f$ over $K$, respectively.

\begin{nota}
With the notation given above, for the inpainting problem the set $D=\Omega\setminus K$ will be the 
non-empty open subset of $\Omega$ representing the damaged region of the image of domain $\Omega$ 
which must be repaired, whereas for scattered data approximations $K$ will be the set of sample points.
\end{nota}

\begin{teo}\label{Thm.3.ImpNew}
Let $\Omega\subset\R^n$ be a bounded open set and $K\subset\Omega$ a compact set. 
Assume $f:K\to \mathbb{R}$ to be bounded, and, for $M>0$, consider  the auxiliary functions defined 
by \eqref{Def.01.AuxFnct.M}.Suppose $\lambda>0$, then, for any $x\in\overline\Omega$,

\begin{equation}\label{Sec3.Eq.13}
 \begin{split}
 	M_{\lambda,\Omega}(f_{\overline\Omega,K}^{M})(x)&=M_{\lambda}(f^{M}_{\mathbb{R}^n,K})(x)\,,\\[1.5ex]
 	M_{\Omega}^{\lambda}(f_{\overline\Omega,K}^{-M})(x)&=M^{\lambda}(f^{-M}_{\mathbb{R}^n,K})(x)\,,
 \end{split}
 \end{equation}
 and

 \begin{equation}\label{Sec3.Eq.14}
 \begin{split}
 	M^{\lambda}_{\Omega}(M_{\lambda,\Omega}(f_{\overline\Omega,K}^{M}))(x)&=M^{\lambda}(M_{\lambda}(f^{M}_{\mathbb{R}^n,K}))(x)\,,\\[1.5ex]
 	M_{\lambda,\Omega}(M^{\lambda}_{\Omega}(f_{\overline\Omega,K}^{-M}))(x)&=M_{\lambda}(M^{\lambda}(f^{-M}_{\mathbb{R}^n,K}))(x)\,.
 \end{split}
 \end{equation}
Furthermore,  let $D=\Omega\setminus K$ and assume that $\Omega^c\cap\overline D=\varnothing$.
Define $O_f=sup_{K}f-\inf_{K}f$ and $\dist(\partial D,\,\partial\Omega)=\inf_{x\in\partial  D}\{\dist(x,\,\partial\Omega)\}$.
If $M>\sup_{\Omega\setminus D}\,f+\lambda\diam^2(D)$ and $\dist^2(\partial D,\,\partial \Omega)>O_f/\lambda$,
the following locality properties hold:
 \begin{itemize}
 	\item[$(i)$] If $x\in\Omega$ is such that 
 		$\dist^2(x,\partial\Omega)>O_f/\lambda$
		then $M_{\lambda,\Omega}(f_{\overline\Omega,K}^{M})(x)$ $[\text{resp. }	M_{\Omega}^{\lambda}(f_{\overline\Omega,K}^{-M})(x)]$
		is determined by values of $f|_{\Omega\setminus D}$. More precisely, if 
 		
		\begin{equation*}
 			\begin{split}
 			&M_{\lambda,\Omega}(f_{\overline\Omega,K}^{M})(x)=f_{\overline\Omega,K}^{M}(z_x)+\lambda|x-z_x|^2\\[1.5ex]
 			&[\text{resp. }	M_{\Omega}^{\lambda}(f_{\overline\Omega,K}^{-M})(x)=f_{\overline\Omega,K}^{-M}(z_x)-\lambda|z_x-x|^2
 			]\,,
 			\end{split}
 		\end{equation*}
 		for some $z_x\in\overline\Omega$, then $z_x\in\Omega\setminus D$.
 	\item[$(ii)$] If $x\in\Omega$ is such that $\dist^2(x,\partial\Omega)>4O_f/\lambda$ then $M^{\lambda}_{\Omega}(M_{\lambda,\Omega}(f_{\overline\Omega,K}^{M}))(x)$ 
 		[resp. $M_{\lambda,\Omega}(M^{\lambda}_{\Omega}(f_{\overline\Omega,K}^{-M}))(x)$]
 		is determined by values of $f|_{\Omega\setminus D}$. More precisely, if
 		
		\begin{equation*}
 			\begin{split}
 			&M^{\lambda}_{\Omega}(M_{\lambda,\Omega}(f_{\overline\Omega,K}^{M}))(x)=M_{\lambda,\Omega}(f_{\overline\Omega,K}^{M})(z_x)-\lambda|x-z_x|^2\\[1.5ex]
 			&[\text{resp. }	
 			M_{\lambda,\Omega}(M^{\lambda}_{\Omega}(f_{\overline\Omega,K}^{-M}))(x)=M^{\lambda}_{\Omega}(f_{\overline\Omega,K}^{-M})(z_x)+\lambda|x-z_x|^2
 			]\,,
 			\end{split}
 		\end{equation*}
 		for some  $z_x\in\overline\Omega$, 
		then 
 		
		\begin{equation*}
 			\begin{split}
				&M_{\lambda,\Omega}(f_{\overline\Omega,K}^{M})(z_x)=f_{\overline\Omega,K}^{M}(y_x)+\lambda|y_x-z_x|^2\\[1.5ex]
 				&[\text{resp. }
				M^{\lambda}_{\Omega}(f_{\overline\Omega,K}^{-M})(z_x)=f_{\overline\Omega,K}^{-M}(y_x)-\lambda|y_x-z_x|^2
				]\,,
 			\end{split}
 		\end{equation*}
		with $z_x,\,y_x\in\Omega\setminus D$.
 \end{itemize}
 \end{teo}
 
 \begin{nota}
 \begin{itemize}
 	\item[$(i)$] Given the definition of $f_{\overline\Omega,K}^{M}$,
 		it is not difficult to show that if $M>\sup_K\,f+\lambda\diam^2(\Omega)$ and if for $x\in\overline\Omega$,
 		$M_{\lambda,\Omega}(f_{\overline\Omega,K}^{M})(x)=f_{\overline\Omega,K}^{M}(z_x)+\lambda|z_x-x|^2$ 
 		with $z_x\in\overline\Omega$, then $z_x$ must belong to $K\cup\partial\Omega$.
 		The locality property of Theorem \ref{Thm.3.ImpNew} is \textit{de facto} making more precise this result 
 		by stating that $z_x\in K$.
 		The proof that, in general, $z_x\in K\cup\partial\Omega$ can be realized by contradiction.
 		Assume $z\in\overline\Omega\setminus(K\cup\partial\Omega)$, then $f_{\overline\Omega,K}^M(z)=M$, thus we have
 		
		\[
 			\begin{split}
 				M &\leq M_{\lambda,\Omega}(f_{\overline\Omega,K}^M)(x)=f_{\overline\Omega,K}^M(z_x)+\lambda|z_x-x|^2   \\[1.5ex]
 				  &= M+\lambda|z_x-x|^2\\[1.5ex]
 				  &\leq f_{\overline\Omega,K}^M(z)+\lambda|z-x|^2\quad\text{for any }z\in\overline\Omega\,.
 			\end{split}
 		\]
 		In particular, the inequality holds also for $y\in K$ such that $|y-x|^2=\dist^2(x,\,K)$, that is,
 		
		\[
 			\begin{split}
 				M&\leq  f_{\overline\Omega,K}^M(y)+\lambda|y-x|^2\\[1.5ex]
 				 &\leq  \sup_{K}\,f+\lambda\dist^2(x,\,K)\\[1.5ex]
 				 &\leq  \sup_{K}\,f+\lambda\diam^2(\Omega)
 			\end{split}
 		\]
 		which is a contradiction. A similar argument can be made for $M^{\lambda}_{\Omega}(f_{\overline\Omega,K}^{-M})(x)$.
 	\item[$(ii)$] The first part of Theorem \ref{Thm.3.ImpNew} on the equality between the local and global Moreau envelopes of
 		our auxiliary functions is, in fact, a consequence of Theorem \ref{Theo.LocLwTr} applied to the functions \eqref{Def.01.AuxFnct.M} 
 		which are bounded in $\Omega$
 		and are of the type \eqref{Eq.01.Ext.Local} and \eqref{Eq.01.Ext.Glob} considered in Theorem \ref{Theo.LocLwTr}. 
 		However, by Theorem \ref{Theo.LocLwTr}, we could only conclude that, for instance, 
 		for $z_x$ as in $(i)$, $(ii)$, we have $z_x\in\Omega$. 
 		Thus, the relevance of the results stated in Theorem \ref{Thm.3.ImpNew} is in the locality properties, that is, 
 		under the conditions on $M$ and $\dist(\partial D,\,\partial\Omega)$ we can conclude that $z_x\in K=\Omega\setminus D$.
	\item[$(iii)$] Theorem \ref{Thm.3.ImpNew}$(ii)$ can be used to get an estimate of the value of $\lambda$ so that 
		the globally defined compensated convex transforms via the mixed critical Moreau envelopes
		coincides with the local compensated convex transforms apart from a boundary layer. For instance, if we want 
		to restric these boundary effects only to a one-pixel wide boundary layer, that is
		for $x\in\Omega$ such that $\dist(x,\partial\Omega)<1$,
		of an image with values in the range $[0,\,1]$
		so that $O_f=1$, we need to take at most $\lambda\geq 4 O_f/\dist^2(x,\partial\Omega)=4$.
 \end{itemize}
 \end{nota}
 

\setcounter{equation}{0}
\section{Numerical Scheme}\label{Sec.NumSchm}
Given the equalities \eqref{Sec3.Eq.03}, \eqref{Sec3.Eq.03a}, the numerical realization of the local compensated convex transforms
boils down to computing the Moreau envelope of a function defined in $\R^n$. Without loss of generality, given the 
relation \eqref{Eq.RelMor} between the upper and  lower Moreau envelopes, we will refer in the following only to the
computation of the lower Moreau envelope. 
To compute the upper Moreau envelope, it is  not difficult to adapt the algorithm proposed,  
or to use the relation \eqref{Eq.RelMor} between the two envelopes.

The computation of the Moreau envelope is an established task in the field of computational convex analysis
\cite{Luc09b} and references therein, that has been tackled by various different approaches aimed at reducing the complexity of a 
direct brute force implementation of the transform. The methods developed in \cite{Luc06,Luc09a,Cor96}, for instance, 
are based on a dimensional reduction. 
The authors exploit the property that the Moreau envelope can be factored by $n$ 1d Moreau envelopes and its
relationship with the Legendre-Fenchel transform \cite[Example 11.26]{RW98}. The factorization of the Moreau envelope is
also used in \cite{FH12} where the construction of the \textit{inf-}convolution is reduced 
to the computation, in constant time, of the envelope of parabolas. 

Moreau envelopes can also be obtained by mathematical morphology operations which can be
particularly useful in the case when $f$ represents an image. Such a class of methods 
can be obtained by an appropriate modification of the ones that compute the 
Euclidean distance transform of binary images.
Here we develop such a method
that generalizes the one used in \cite{HM94,Shi09,SM92} to compute the discrete Euclidean distance transform. 

The fundamental idea is 
the characterization \eqref{Sec2.Eq.EucDistMor} of the Euclidean distance in terms of the 
Moreau envelope of the characteristic function.
The Euclidean distance transform is the erosion of the characteristic function by the 
quadratic structuring element whereas the Moreau envelope is the erosion of the image $f$. 
Thus, one can think of generalizing the Euclidean distance transform of binary images, by replacing the
binary image by an arbitrary function on a grid.
The decomposition of the structuring element which yields the exact Euclidean distance transform \cite{SM92}
into basic ones, yields a simple and fast algorithm  where the discrete lower Moreau envelope can be 
computed by a sequence of local operations, using one-dimensional neighborhoods. 
We will use the same structuring element as in \cite{HM94}
and show that we recover the exact discrete lower Moreau envelope. 

Unless otherwise stated, in the following, $i,\,j,\,k,\, r,\,s,\, p,\, q\in \mathbb{Z}$ denote 
integers whereas $m,\, n\in\mathbb{N}$ are non-negative integers. 
Given $n\geq 1$, we introduce grid of points of the space $\mathbb{R}^n$ with regular spacing $h>0$
denoted by $x_k\in\R^n$, $k\in\mathbb{Z}$. 

\begin{defi}\label{Def.4.DiscrMor}
Suppose $f:\mathbb{R}^n\to\mathbb{R}\cup\{\infty\}$ satisfies
$f(x)\geq -c_0|x|^2-C_1$ with $c_0,\, C_1\geq 0$. Let $h>0$, $n\geq 1$ and denote by $x_k$ a point
of the grid of $\mathbb{R}^n$ of size $h$. Then the discrete Moreau lower envelope at 
$x_k\in\mathbb{R}^n$ is defined by

\begin{equation}\label{Sec4.Def.DiscrMor} 
	M_\lambda^h(f)(x_k)=
	\inf\{ f(x_k+rh)+ \lambda h^2|r|^2,\; r\in \mathbb{Z}^n\}.
\end{equation}
for $\lambda>c_0$.
\end{defi}

The approximation of Moreau lower envelope by the discrete Moreau lower envelope  
is quantified by the following estimation result. We consider the case where $f$ is a 
uniformly continuous function first. 

\begin{teo}\label{Sec4.Theo.CnvrgMor}
Let $f:\mathbb{R}^n\to\mathbb{R}$ be a function with modulus of continuity $\omega_f$ such that 
$\omega_f(t)\leq a t+b$
with $a>0$ and $b\geq 0$ for any $t\geq 0$. Assume $h>0$ and $n\geq 1$. Then,
for any grid point $x_k$ of the grid of $\mathbb{R}^n$ of size $h$, 

\begin{equation}\label{Sec4.Eq.ErrEst}
	\left|M_\lambda^h(f)(x_k)-M_\lambda(f)(x_k)\right|\leq  \omega_f(h \sqrt{n})+2\lambda h^2{n}+2h\sqrt{\lambda}\,d(\lambda)
\end{equation}
where $d(\lambda)=\sqrt{\omega_f\left(a/\lambda+\sqrt{b/\lambda}\right)}$.
\end{teo}

\begin{nota}
Since the modulus of continuity $\omega_f=\omega_f(t)$ tends to zero as $t\to 0+$ 
(see Proposition \ref{Sec2.Eq.Pro.MoC}$(i)$),
it follows that the error bound in \eqref{Sec4.Eq.ErrEst} tends to zero when $h\to 0+$.
\end{nota}

The rate of convergence in \eqref{Sec4.Eq.ErrEst} can be improved for $L-$Lipschitz functions $f$. In this case, 
we have the following result.

\begin{coro}\label{Sec4.Cor.CnvrgMorLips}
Let $f$ be an $L-$Lipschitz function in $\R^n$. Assume $h>0$ and $n\geq 1$. Then,
for any grid point $x_k$ of the grid of $\mathbb{R}^n$ of size $h$, 

\begin{equation}\label{Sec4.Eq.ErrEstLip}
	\left|M_\lambda^h(f)(x_k)-M_\lambda(f)(x_k)\right|\leq (2+\sqrt{n})Lh+2\lambda h^2{n}\,.
\end{equation}
\end{coro}

\begin{nota}
	Under the same regularity conditions, \cite{Luc06} obtains the same converegnce rate 
	for a real function of one real variable. However, in contrast to the scheme \cite{Luc06}
	we obtain directy the values of the Moreau envelope at the grid points
	by a scheme which has the same complexity as the distance transform \cite{Shi09}
	and can be easily implemented and applied to any dimension.
\end{nota}

In Definition \ref{Def.4.DiscrMor}, the infimum is taken over infinitely many grid points thus its
computation is not practical. Therefore we introduce the $m$-th approximation of the discrete 
Moreau lower envelope where the order 
of approximation $m$ is related to the number of nodes that are taken to compute the infimum 
in the definition of the envelope. 

\begin{defi}\label{Def.4.AprxDiscrMor}
Let $n\geq 1$ and denote by  $x_k$ a point of the grid of $\mathbb{R}^n$ of size $h$. 
The $m$-th approximation of the $n-$th dimensional discrete Moreau lower envelope
$M_\lambda^h(f)$ is given by

\begin{equation}\label{Sec4.Eq.02} 
	g_m(x_k)=\inf\{ f(x_k+rh)+ \lambda h^2|r|^2,\;  r\in\mathbb{Z}^n,\,|r|_{\infty} \leq m\}
\end{equation}
for $\lambda>c_0$, where $|r|_\infty$ is the infinity norm of $r\in \mathbb{Z}^n$.
\end{defi}

Given $m\geq 1$, to evaluate $g_m(x_k)$
at any point $x_k$ of the grid of $\R^n$ of size $h$, we can, in fact, consider the 
values $f_m(x_k)$ that are obtained by applying Algorithm \ref{Algo:MoreauEnv}: 

\begin{algorithm}[H]
\begin{algorithmic}[1]
	\STATE{Set $\displaystyle i=1,\, m\in\mathbb{N}$}
	\STATE{$\displaystyle \forall x_k,\, f_0(x_k)=f(x_k)$}
	\WHILE{$\displaystyle i<m$}
	\STATE{$\displaystyle \tau_i=2i-1$
		}
	\STATE{$\displaystyle f_{i}(x_k)=\min\{ f_{i-1}(x_k+rh)+\lambda h^2|r|^2\tau_{i}: \,
		r\in\mathbb{Z}^n,\,|r|_{\infty}\leq 1\}$}
	\STATE{$i\leftarrow i+1$}
	\ENDWHILE
\end{algorithmic}
\caption{\label{Algo:MoreauEnv} Computation of $f_m(x_k)$ at the points $x_k$
of the grid of $\R^n$ of size $h$ for given $m\geq 1$.}
\end{algorithm}

The relation between $f_m(x_k)$ and  $g_m(x_k)$ is described in the following theorem.

\begin{teo}\label{Sec4.Theo.DisSch}
Let $f$ be bounded in $\R^n$. Assume $\lambda>0$. Then for all the grid points $x_k$ 
of the grid of $\R^n$ of size $h$, and $m=0,1,2,\dots$,

\[	
	f_m(x_k)=g_m(x_k)
\]
with $f_m(x_k)$ computed by applying Algorithm \ref{Algo:MoreauEnv}. 
\end{teo}

By definition, for every grid point $x_k$, $g_m(x_k)$ is decreasing in $m$. Since $g_m(x_k)$ is bounded from below
by $M^h_\lambda(f)(x_k)$, $g_m$ will then converge as $m$ goes to $\infty$.
The following result actually shows that it will take only finitely many iterations for $g_m(x_k)$
to reach $M^h_\lambda(f)(x_k)$.

\begin{prop}\label{Sec4.Prop.Est}
Let $f$ be bounded in $\R^n$. 
Assume $h>0$ and denote by $x_k$ a point
of the grid of $\mathbb{R}^n$ of size $h$. Then 

\[
	M^h_\lambda(f)(x_k)=g_m(x_k)	\,,
\]
for $m\geq\floor{\tfrac{1}{h}\sqrt{\tfrac{\osc(f)}{\lambda}}}+1$, where $\floor{x}$ denotes the integer part of $x$.
\end{prop}

\begin{nota}
For an $8-$bit image with $h=1$ the pixel size and $\osc(f)=255$, if we take $m\geq\floor{16/\sqrt{\lambda}}+1$,
we will have $g_m(x_k)=M_{\lambda}^hf(x_k)$ at any grid point $x_k$. 
\end{nota}

For completeness, we conclude this section by giving the scheme that we use for the implementation of 
the convex based definition of the compensated transforms. The scheme is a generalization of 
the one introduced in \cite{Obe08} that is briefly summarized in Algorithm \ref{Algo:CnvxEnv} 
and described below. 
Given a uniform grid of points $x_k\in\mathbb{R}^n$, equally spaced with grid size $h$,  denote by
$S_{x_k}$ the $d-$point stencil of $\mathbb{R}^n$ with center at $x_k$ defined as 
$S_{x_k}=\{x_k+hr, |r|_{\infty}\leq 1, r\in\mathbb{Z}^n\}$ with $|\cdot|_{\infty}$ the $\ell^{\infty}$-norm of 
$r\in\mathbb{Z}^n$ and $d=\#(S)$, cardinality of the finite set $S$. At each grid point $x_k$ we compute an 
approximation of the 
convex envelope of $f$ at $x_k$ by an iterative scheme where each iteration step $m$ is given by

\[
	(\co f)_m(x_k)=\min\Big\{ 
		f(x_k),\,\sum\lambda_i (\co f)_{m-1}(x_i):\,\,\,
		\sum \lambda_i=1,\,\lambda_i\geq 0,\,x_i\in S_{x_k}\Big\}
\]
with the minimum taken between $f(x_k)$ and only some convex combinations of $(\co f)_{m-1}$ 
at the stencil grid points $x_i$ of $S_{x_k}$.
Likewise \cite{Obe08}, the scheme can be shown to converge  
but there is no estimate of the rate of convergence.

\begin{algorithm}
\begin{algorithmic}[1]
	\STATE{Set $\displaystyle m=1,\,(\co f)_{0}=f,\,\,tol$}
	\STATE{$\displaystyle \epsilon=\|f\|_{L^2}$}
	\WHILE{$\displaystyle \epsilon>tol$}
	\STATE{$\displaystyle \forall x_k, \quad (\co f)_{m}(x_k)=\min\Big\{ 
		f(x_k),\,\sum\lambda_i (\co f)_{m-1}(x_i):\,\,\,
		\sum \lambda_i=1,\,\lambda_i\geq 0,\,x_i\in S_{x_k}\Big\}$
		}
	\STATE{$\displaystyle \epsilon=\|(\co f)_m-(\co f)_{m-1}\|_{L^2}$}
	\STATE{$m\leftarrow m+1$}
	\ENDWHILE
\end{algorithmic}
\caption{\label{Algo:CnvxEnv} Computation of the convex envelope of $f$ according to \cite{Obe08}}
\end{algorithm}

 
\setcounter{equation}{0}
\section{Numerical Experiments}\label{Sec.NumExprm}

In this section we first present a one-dimensional and a two-dimensional prototype example with analytical expression for 
$C_{\lambda,\Omega}^l(f_{\overline\Omega}^-)$  which we use: 
\begin{itemize}
	\item[$(i)$] to verify the numerical scheme introduced in Section \ref{Sec.NumSchm}; 
	\item[$(ii)$] to compare the compensated convex transforms using the local Moreau envelope based definition 
			and the convex based definition.
	\item[$(iii)$] to analyze different behaviors of the transformations  
		$C_{\lambda,\Omega}^l(f_{\overline\Omega}^-)$ and $C_{\lambda}^l(f^{\infty})$ at the boundary of the domain. 
\end{itemize}
For computing the compensated convex transforms by using the Moreau envelope based definition, we use the iterative scheme described by 
Algorithm \ref{Algo:MoreauEnv}, whereas for the realization of the 
convex based definition of the compensated
convex transforms, we apply Algorithm \ref{Algo:CnvxEnv}, which has already been
employed to carry out numerical examples of \cite{ZCO15c,ZOC15b}.
In the application of Algorithm \ref{Algo:MoreauEnv}, rather than fixing the number $m$ of iterations,
we introduce the convergence check on the $\ell^{\infty}$ norm of the error between two succesive iterates
such as we do when we apply Algorithm \ref{Algo:CnvxEnv}.
We then describe numerical experiments on applications of the local 
compensated convex transforms on a bounded closed convex domain for the extraction of
the multiscale medial axis map, the extraction of Hausdorff stable intersections of smooth manifolds
and finally, for the interpolation and approximation of sampled functions.


\begin{ex}\label{Sex5.NumEx.1d}\textbf{A one-dimensional prototype example.}
We  consider a one-dimensional model problem given by the piecewise affine double well model

\begin{equation}
	f(x)=\dist(x,\,\{-1,\,1\})=\min\{|x-1|,\,|x+1|\}
\end{equation}
for $x\in\overline\Omega:=\{x\in\mathbb{R},\; |x|\leq 2\}$. Let the corresponding 
$f_{\overline\Omega}^-$ be defined by \eqref{Eq.01.Ext.Local.LW}. For  $\lambda\geq 1$, it is  not difficult 
to show that 

\begin{equation}\label{Ex.1d.ExLwTr}
C^l_{\lambda,\Omega}(f^{-}_{\overline\Omega})(x)=\left\{
		\begin{array}{ll}
	\displaystyle A(|x|-1),	   & \displaystyle \left||x|-\frac{2+x_1}{2}\right|\leq\frac{2-x_1}{2}\,,\\[2.5ex]
	\displaystyle f(x),	   & \displaystyle \left||x|-\frac{x_1+x_2}{2}\right|\leq\frac{x_1-x_2}{2}\,,\\[1.5ex]
	\displaystyle 1-\frac{1}{4\lambda}-\lambda x^2, & \displaystyle |x|\leq x_2,\\[1.5ex]
	\displaystyle 0, &\displaystyle \text{otherwise}\,,
		\end{array}\right.
\end{equation}
where $x_1=2-\sqrt{\lambda}/\lambda$, $x_2=1/(2\lambda)$ and 
$A(x)=-\lambda x^2+(2\lambda-2\sqrt{\lambda}+1)x-\lambda+2\sqrt{\lambda}-1$.
Given the definition of $f_{\overline\Omega}^-$, at the  boundary nodes $\partial\Omega=\{-2,\,2\}$, we have that 
$M_{\lambda,\Omega}(f_{\overline\Omega}^-)(x)=M_{\Omega}^{\lambda}(M_{\lambda,\Omega}(f_{\overline\Omega}^-))(x)=\inf f$ 
for $x\in\partial\Omega$. After choosing a uniform grid of $\Omega$ with grid size $h$ which we denote next as $\Omega_h$, 
we can therefore run Algorithm \ref{Algo:MoreauEnv} only at the interior points $x_k\in\Omega_h$ and assume 
$f_i=\inf f$ when the scheme is 
applied at the first grid point $x_k$ of $\Omega_h$ next to $-2$ and $2$.
For a given grid size $h$ and for any given $m$, number of iterations used in the application of Algorithm \ref{Algo:MoreauEnv},
we compute the $\ell^{\infty}$-norm of the error defined as

\[
	\|e\|_{\ell^{\infty}}=\max\{|C_{\lambda,\Omega}^l(f^{-}_{\overline\Omega})(x_k)-
	M^{\lambda,h}_{\Omega}(M_{\lambda,\Omega}^h(f^{-}_{\overline\Omega}))(x_k)|,\,x_k\in\Omega_h\}\,.
\]

\begin{figure}[htbp]
\begin{center}
	$\begin{array}{cc}
		\includegraphics[width=0.48\textwidth]{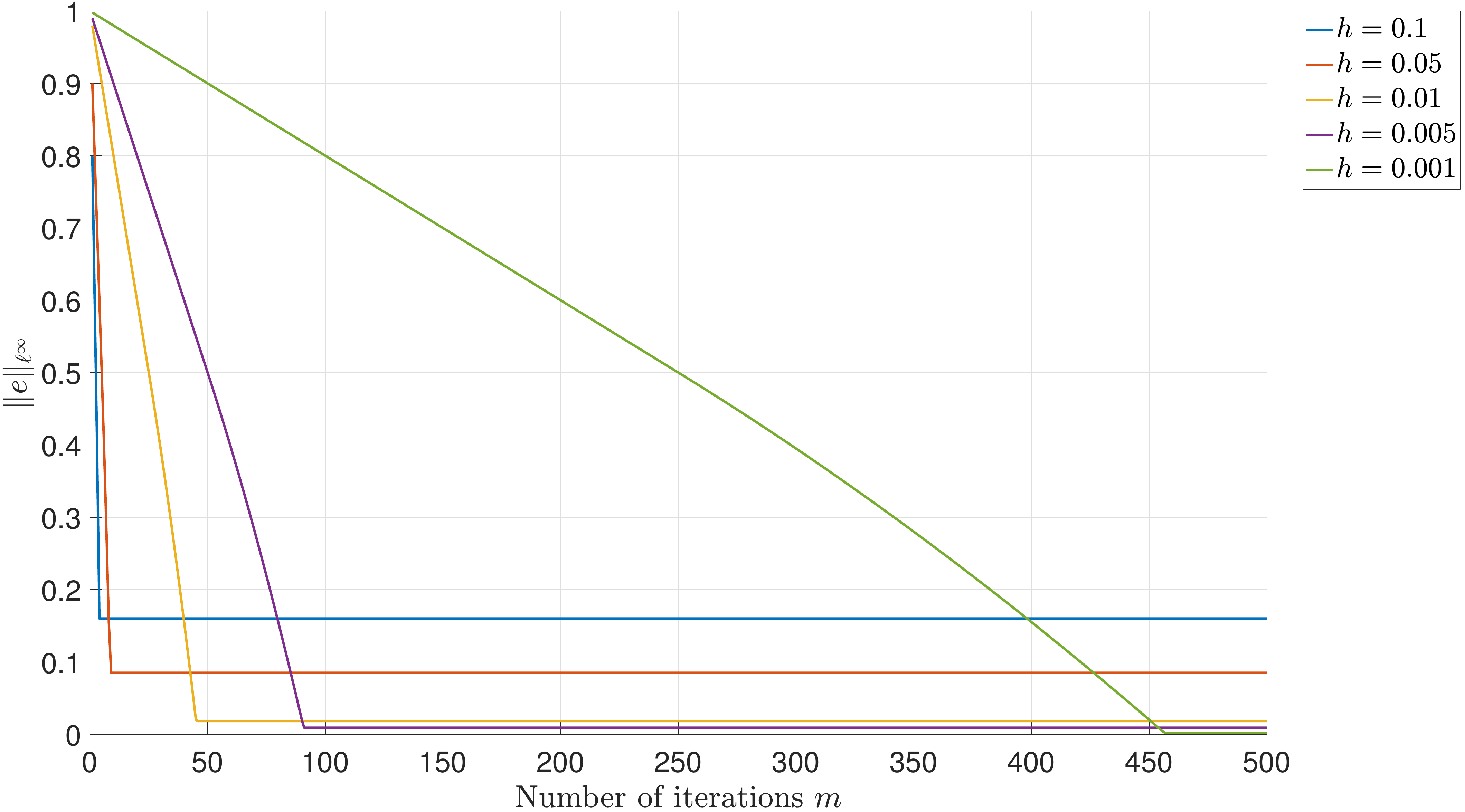}&
		\includegraphics[width=0.48\textwidth]{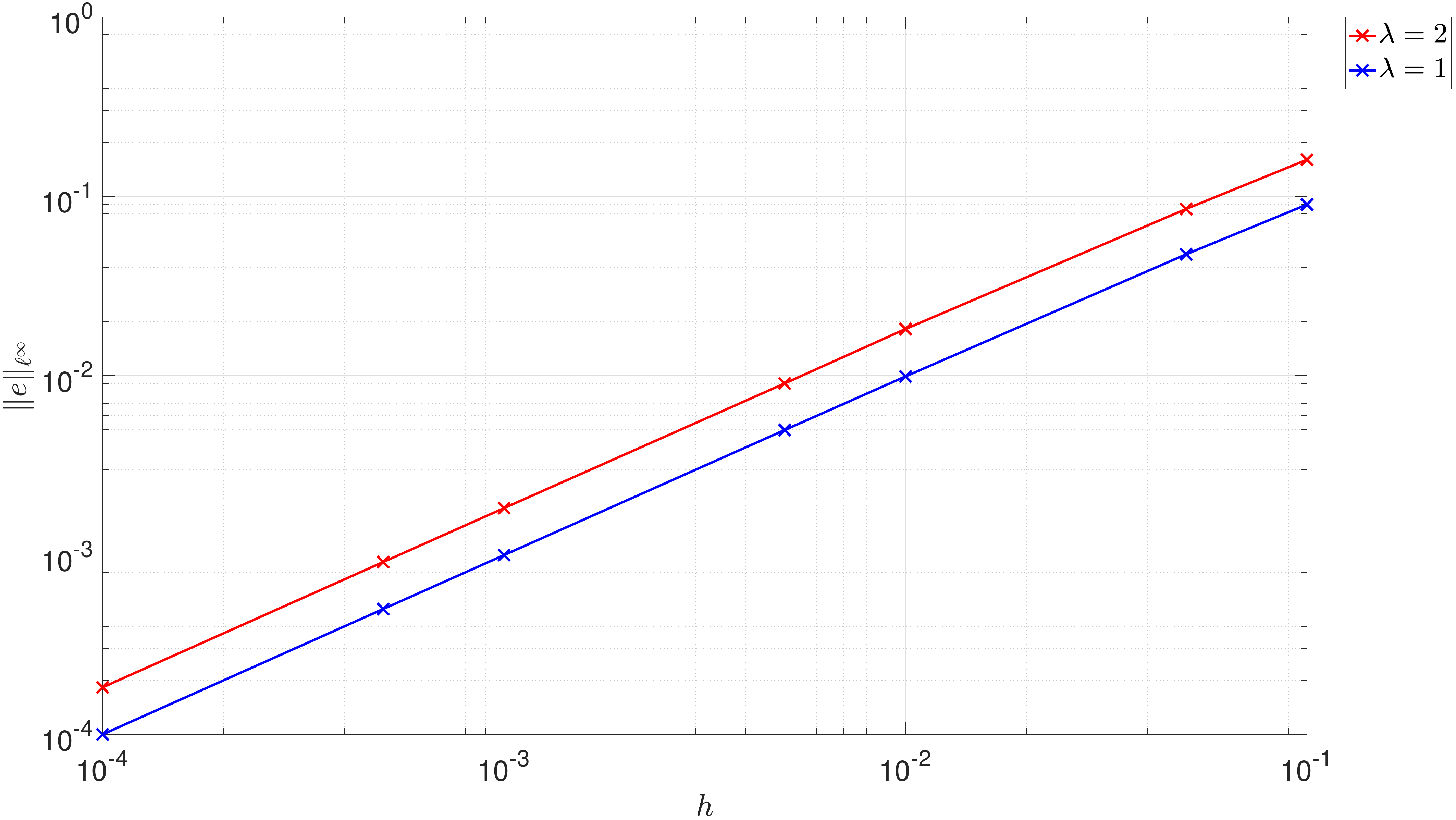}\\
		$(a)$ & $(b)$
	\end{array}$
\end{center}
   \caption{\label{Ex01.2Wel:LclCmp}
$(a)$ Variation of $\|e\|_{\ell^{\infty}}$ with the number $m$ of iterations and for different values of the 
grid size $h$ and $\lambda=2$; $(b)$ Convergence plot of the error with the grid size $h$ and 
$\lambda=1$ and $\lambda=2$.}
\end{figure}

Figure \ref{Ex01.2Wel:LclCmp}$(a)$ displays  the convergence plot with respect to the number of iterations for 
different values of the grid size $h$ and for given values of the parameter $\lambda$, $\lambda=1$ and $\lambda=2$, 
whereas Figure \ref{Ex01.2Wel:LclCmp}$(b)$ 
shows the convergence plot
with respect to the grid size $h$, using for each $h$ the value of $m$ such that the $\ell^{\infty}$
norm of the error between two iterates is not greater than $10^{-7}$ and for different values of $\lambda$.
We observe that the number $m$ of iterations to obtain convergence increases as 
$h$ is reduced, consistently with the theoretical finding of Proposition \ref{Sec4.Prop.Est}, whereas 
Figure \ref{Ex01.2Wel:LclCmp}$(b)$ exhibits the linear convergence rate of the scheme as predicted by 
Corollary \ref{Sec4.Cor.CnvrgMorLips}.
The graph of \eqref{Ex.1d.ExLwTr} and $M^{\lambda,h}_{\Omega}(M_{\lambda,\Omega}^h(f^{-}_{\overline\Omega}))$, 
with the latter corresponding to the grid size $h=0.01$ and $\lambda=2$ are shown in Figure 
\ref{Ex01.2Wel:GraphLclCmp}.

\begin{figure}[htbp]
\centerline{$\includegraphics[width=0.50\textwidth]{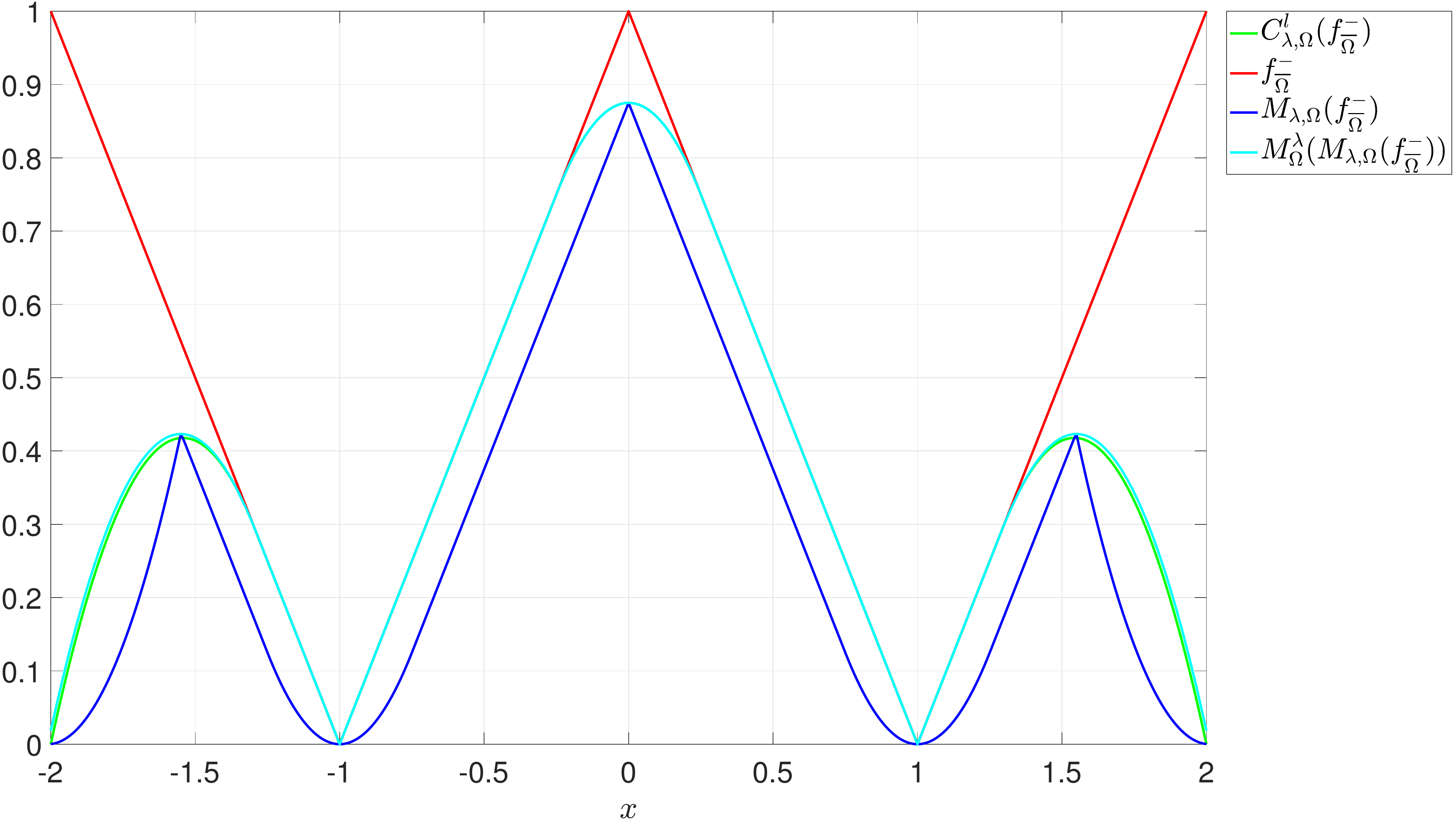}$}
   \caption{\label{Ex01.2Wel:GraphLclCmp}
Graph of \eqref{Ex.1d.ExLwTr}, $M_{\lambda,\Omega}^h(f^{-}_{\overline\Omega})$ and 
$M^{\lambda,h}_{\Omega}(M_{\lambda,\Omega}^h(f^{-}_{\overline\Omega}))$, 
with $h=0.01$ and $\lambda=2$.}
\end{figure}

To compare the Moreau computation of $C_{\lambda,\Omega}^l(f^{-}_{\overline\Omega})$ given by 
$M^{\lambda,h}_{\Omega}(M_{\lambda,\Omega}^h(f^{-}_{\overline\Omega}))$, 
to the convex based definition
which relies on the equality $C_{\lambda,\Omega}^l(f^{-}_{\overline\Omega})=C_{\lambda}^l(f^-_{\R^n})$, 
Table \ref{Ex01.Tab:1d} reports, for different grid size $h$, the error $\|e\|_{\ell^{\infty}}$ and 
the number of iterations $m$ that compute
$M^{\lambda,h}_{\Omega}(M_{\lambda,\Omega}^h(f^{-}_{\overline\Omega}))$ and $C_{\lambda}^l(f^-_{\R^n})$ such that  
the $\ell^{\infty}$ norm of the error between two iterates is not greater than $10^{-7}$.
Inspection of Table \ref{Ex01.Tab:1d} shows that for a given grid size $h$ the Moreau based computation
uses a much  lower number of iterations, especially for small $h$,
and the discrete Moreau based lower transform is more accurate than
the discrete convex based lower transform. Furthermore, for small values of $h$,
we also note that the convex based scheme (Algorithm \ref{Algo:CnvxEnv}) appears not to be numerically stable. 
The number of iterations $m$ relative to the computation of
$M^{\lambda,h}_{\Omega}(M_{\lambda,\Omega}^h(f^{-}_{\overline\Omega}))$ is the total number of iterations which sums up 
the iterations for computing the lower Moreau envelope and the upper Moreau envelope.


\begin{table}[tbhp]
\centerline{
\begin{tabular}{c|c|c|c|c|c|c|c|c|}
\cline{2-9}
	& \multicolumn{4}{|c|}{ $\lambda=1$} & \multicolumn{4}{|c|}{ $\lambda=2$}\\
\cline{2-9}
	& \multicolumn{2}{c|}{ $M^{\lambda,h}_{\Omega}(M^{h}_{\lambda,\Omega}(f_{\overline\Omega}^-))$}    & \multicolumn{2}{c|}{ $C_{\lambda}^l(f_{\R^n}^{-})$ } & \multicolumn{2}{c|}{ $M^{\lambda,h}_{\Omega}(M^{h}_{\lambda,\Omega}(f_{\overline\Omega}^-))$}    & \multicolumn{2}{c|}{ $C_{\lambda}^l(f_{\R^n}^{-})$ } \\ \hline
\multicolumn{1}{ |c| }{$h$}	         & $m$     & $\|e\|_{\ell^{\infty}}$   & $m$     &  $\|e\|_{\ell^{\infty}}$ & $m$     & $\|e\|_{\ell^{\infty}}$   & $m$     &  $\|e\|_{\ell^{\infty}}$ \\ \hline
\multicolumn{1}{ |l| }{$10^{-1}$}	 & $14$  & $0.09$     & $254$   &  $0.10$      & $12$  & $0.16$   & $130$   &  $0.1771431$	\\ \hline
\multicolumn{1}{ |l| }{$5\cdot 10^{-2}$} & $24$  & $0.00475$  & $905$   &  $0.05$      & $21$  & $0.085$  & $472$   &  $0.09$	\\ \hline
\multicolumn{1}{ |l| }{$10^{-2}$}        & $104$ & $0.0099$   & $15871$ &  $0.0100032$ & $95$  & $0.0182$ & $8669$  &  $0.0182276$\\ \hline
\multicolumn{1}{ |l| }{$5\cdot 10^{-3}$} & $204$ & $0.004975$ & $52053$  & $0.0050064$ & $185$ & $0.00905$ & $28870$    &  $0.009132$	\\ \hline
\multicolumn{1}{ |l| }{$10^{-3}$}        & $1004$& $0.000999$ & $651472$ & $0.0104497$ & $917$ & $0.0001826$ & $394000$ &  $0.0059748$\\ \hline
\multicolumn{1}{ |l| }{$5\cdot 10^{-3}$} & $2004$& $0.0004997$& $1572046$& $0.0377351$ & $1831$& $0.0009135$& $1037066$ &  $0.0199438$	\\ \hline
\multicolumn{1}{ |l| }{$10^{-4}$}        & $9999$& $0.0000999$& $-$      &  $-$	       & $9143$& $0.000128$ & $-$       &  $-$	\\ \hline
\end {tabular}
}
\caption{\label{Ex01.Tab:1d} Number of iterations $m$ and values of the error 
$\|e\|_{\ell^{\infty}}$ of  
$M^{\lambda,h}_{\Omega}(M^{h}_{\lambda,\Omega}(f_{\overline\Omega}^-))$  and $C_{\lambda}^l(f_{\R^n}^{-})$ computed by applying 
Algorithm \ref{Algo:MoreauEnv} and Algorithm \ref{Algo:CnvxEnv}, respectively, for different values of 
the grid size $h$ and $\lambda$.
The number of iterations $m$ given in the table corresponds to the termination criteria with 
the $\ell^{\infty}$ norm of the error between two successive iterates not greater than $10^{-7}$.
The reported values of $m$ for computing $M^{\lambda,h}(M^{h}_{\lambda}(f_{\overline\Omega}^-))$ are 
the total number of iterations.}
\end{table}


We conclude this example by looking at the behavior of  $C^l_\lambda(f^\infty)$ with $f^{\infty}$ defined by 
\eqref{Eq:Intro:01}. For $\lambda\geq 1/2$, we have then the following explicit formula for $C^l_\lambda(f^\infty)$,

 \begin{equation}\label{Ex1.LwTrInftyEx}
 C^l_\lambda(f^\infty)(x)=\left\{
 		\begin{array}{ll}
 	\displaystyle 1-\frac{1}{4\lambda}-\lambda x^2,  & \displaystyle |x|\leq \frac{1}{2\lambda},\\[1.5ex]
 	\displaystyle f^\infty(x),          & \displaystyle |x|\geq \frac{1}{2\lambda}.
 		\end{array}\right.
 \end{equation}
 The graph of $C^l_\lambda(f^\infty)(x)$ is displayed in Figure \ref{Ex01.2Wel:InftCmp}$(a)$ along with that  of 
 $C^l_{\lambda,\Omega}(f^{-}_{\overline\Omega})$ given by \eqref{Ex.1d.ExLwTr}. 
 Note  that the two transforms differ at the boundary of $\Omega$, as a result of 
 the type of singularity therein introduced by the definition of $f^{\infty}$ and 
 $f^{-}_{\overline\Omega}$, respectively.
 
 \begin{figure}[htbp]
 \begin{center}
 	$\begin{array}{cc}
 		\includegraphics[width=0.45\textwidth]{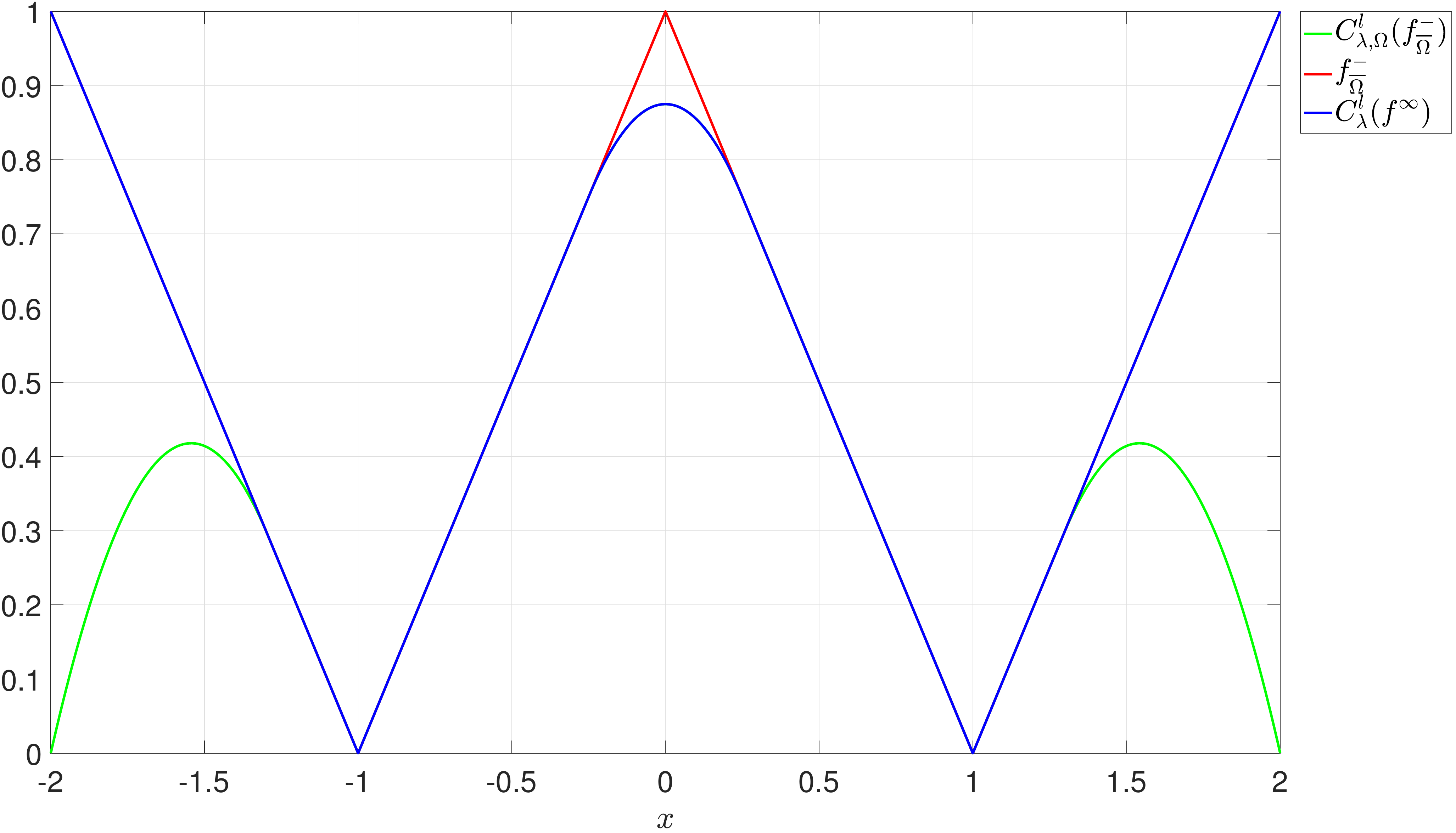}&
 		\includegraphics[width=0.45\textwidth]{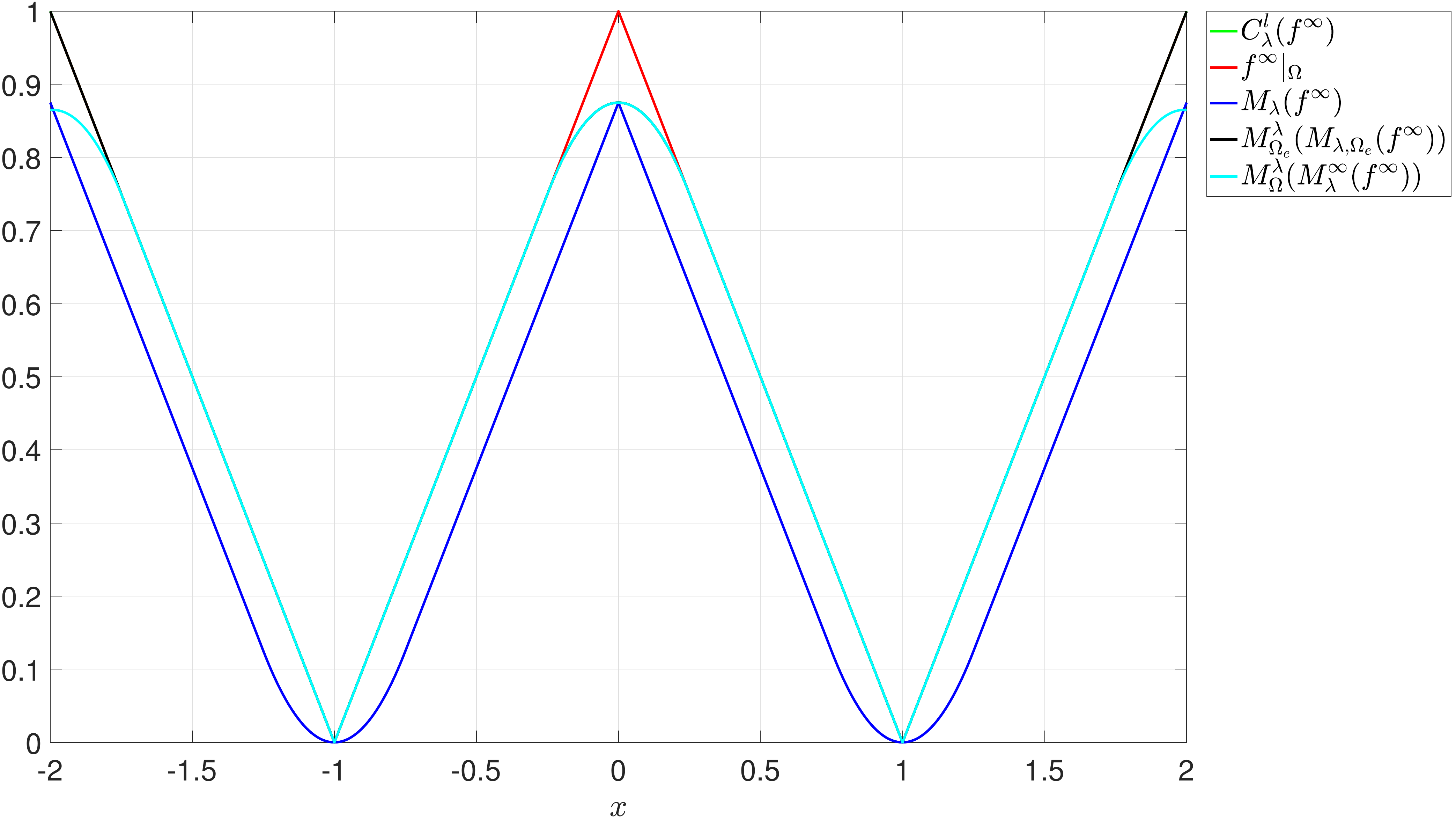}\\
 		$(a)$ & $(b)$ 
 	\end{array}$
 \end{center}
    \caption{\label{Ex01.2Wel:InftCmp}
 $(a)$ Graph of $C^l_\lambda(f^\infty)(x)$ and $C^l_{\lambda,\Omega}(f^{-}_{\overline\Omega})$
 for $\lambda=2$ as given by  \eqref{Ex1.LwTrInftyEx} and \eqref{Ex.1d.ExLwTr}, respectively. 
 $(b)$ Graph of $C^l_\lambda(f^\infty)(x)$, $M^h_{\lambda,\Omega_e}(f^{\infty})$ and 
 of $M_{\Omega_e}^{\lambda,h}(M^h_{\lambda,\Omega_e}(f^{\infty}))$ for $\Omega_e=\Omega$ and $\Omega_e$ 
 given by \eqref{Ex.1d.ExtDom}, and by taking $h=0.01$ and $\lambda=2$. For the simulations, in place
 of $f^{\infty}$ we have considered $f^M(x)=\chi_{\Omega}(x)f(x)+(1-\chi_{\Omega}(x))M$ with $M=10^3$. 
 }
 \end{figure}
 
 At  first sight, the computation of $C^l_\lambda(f^\infty)$ as $M^{\lambda}(M_{\lambda}(f^\infty))$ would be faced 
 with the problem of having to evaluate $M_{\lambda}(f^\infty)$ in $\R^n$. In this case, one could therefore 
 think of resorting to the convex based definition by the assumption that 
 $C^l_\lambda(f^\infty)=f^\infty$ on $\partial \Omega$ and the application of a scheme that computes the convex 
 envelope of a function. Although there are several such schemes in the literature,
 such as the quickhull algorithm \cite{BDH96} or the one introduced in \cite{Obe08}, 
 their application can be quite cumbersome  when applied to compute the convex envelope of a function defined in $\R^n$ with $n>1$, or 
 it can exhibit a slow and unknown rate of convergence as in \cite{Obe08} which would not allow any 
 prediction for the rate of error reduction. As a result, even in the case of $C^l_\lambda(f^\infty)$
 one might wonder whether it is possible to use Algorithm \ref{Algo:MoreauEnv} to obtain an approximation
 of $C^l_\lambda(f^\infty)$. Given the localization effect of the \text{inf-} and \text{sup-}convolution
 with quadratic perturbations, we can also  use the Moreau based definition of $C^l_\lambda(f^\infty)$
 and apply Algorithm \ref{Algo:MoreauEnv} provided that one computes $M_{\lambda}(f^{\infty})$
 and $M^{\lambda}(M_{\lambda}(f^{\infty}))$ over an extended domain $\Omega_e$ that contains $\Omega$.
 Figure \ref{Ex01.2Wel:InftCmp}(b) displays the graph of the transforms in the case of $\Omega_e=\Omega$
 and of $\Omega\subset\Omega_e$ with 
 
 \begin{equation}\label{Ex.1d.ExtDom}
 	\Omega_e=]-2-1/(2\lambda),\,2+1/(2\lambda)[\,.
 \end{equation}
 In the first case, we are actually computing the transform 
 $M_{\Omega}^{\lambda}(M_{\lambda}^{\infty}(f^{\infty}))$ where we have set
 
 \[
 	M_{\lambda}^{\infty}(f^{\infty})(x)=\left\{\begin{array}{ll}
 						M_{\lambda}(f^{\infty})(x)\,, & x\in\Omega\,,\\[1.5ex]
 						\infty\,, &\text{otherwise,}
 						\end{array}\right.
 \]
 which produces a boundary error, whereas in the second case we obtain an excellent approximation of 
 $C_{\lambda}^l(f^{\infty})$ with $\|e\|_{\infty}$=0.0002. Note that how big the domain $\Omega_e$ must be to ensure that 
 $C^l_\lambda(f^\infty)(x)=M_{\Omega_e}^{\lambda}(M_{\lambda,\Omega_e}(f^{\infty}))(x)$ for $x\in\Omega$,
 is an open question.
 
\end{ex}


\begin{ex}\label{Sex5.NumEx.2d}\textbf{A two-dimensional prototype example.}
Let $\Omega=B(O;\,2)$ be the open ball with center at the origin $O\in\R^2$ and radius $r=2$.
Consider the squared distance of $x\in\Omega$ to the boundary $\partial B(O;\,1)$ given by 

\begin{equation}
	f(x,y)=\dist^2(x,\,\partial B(O;\,1)),\quad x\in B(O;\,2)\,,
\end{equation}
and the functions $f_{\overline\Omega}^{-}$ and $f^{\infty}$ defined by \eqref{Eq.01.Ext.Local.LW}
and \eqref{Eq:Intro:01}, respectively.
Given the radial symmetry of $f_{\overline\Omega}^{-}$ and $f^{\infty}$, it is not difficult to verify
that, for $\lambda\geq 1$, 

\begin{equation}\label{Eq.2d.LwTrEx}
	C_{\lambda,\Omega}^l(f_{\overline\Omega}^{-})(x)=\left\{\begin{array}{ll}
		\displaystyle 0,&\displaystyle r>2\,,\\[1.5ex]
		\displaystyle -|m(r-2)|+4\lambda-\lambda r^2,& \displaystyle \left|r-\frac{2+x_p}{2}\right|\leq \frac{2-x_p}{2}\,,\\[2.5ex]
		\displaystyle f_{\overline\Omega}^{-}(r), & \displaystyle \left|r-\frac{x_s+x_p}{2}\right|\leq \frac{x_p-x_s}{2}\,,\\[2.5ex]
		\displaystyle \frac{\lambda}{1+\lambda}-\lambda r^2, &\displaystyle r<x_s\,,
		\end{array}\right.
\end{equation}
and

\begin{equation}\label{Eq.2d.LwTrInftyEx}
	C_{\lambda}^l(f^{\infty})(x)=\left\{\begin{array}{ll}
		\displaystyle f^{\infty}(r),&\displaystyle r>x_s\,,\\[1.5ex]
		\displaystyle \frac{\lambda}{1+\lambda}-\lambda r^2,& \displaystyle r\leq x_s\,,
		\end{array}\right.
\end{equation}
where $r=|x|$, $x_p=2-1/\sqrt{1+\lambda}$, $x_s=1/(1+\lambda)$ and $m=2(1+2\lambda)-2\sqrt{1+\lambda}$.
Though for the computation of \eqref{Eq.2d.LwTrEx} and \eqref{Eq.2d.LwTrInftyEx} we could exploit the symmetry of 
$f_{\overline\Omega}^{-}$ and $f^{\infty}$ and reduce their evaluation to $1d$ problems, in order to verify our scheme for 
$2d$ applications, we will not take the symmetry into account and will refer to $f_{\overline\Omega}^{-}$ and $f^{\infty}$
as generic functions of $x\in\R^2$. Let $D\subset\R^2$ be a box that contains $\Omega$, for instance, let us take

\[
	D=]-2.5,\,2.5[\,\times\, ]-2.5,\,2.5[
\] 
and consider $f_{\overline\Omega}^{-}$ and $f^{\infty}$ to be extended over $D$ by setting
$f_{\overline\Omega}^{-}(x)=0$ for $x\in D\setminus\Omega$, respectively. 
Since Algorithm \ref{Algo:MoreauEnv} is formulated for a square lattice, we will refer to the above extensions over $D$ 
for the application of the algorithm and still denote them by $f_{\overline\Omega}^{-}$ and $f^{\infty}$.

\begin{figure}[htbp]
\begin{center}
	$\begin{array}{cc}
		\includegraphics[width=0.45\textwidth]{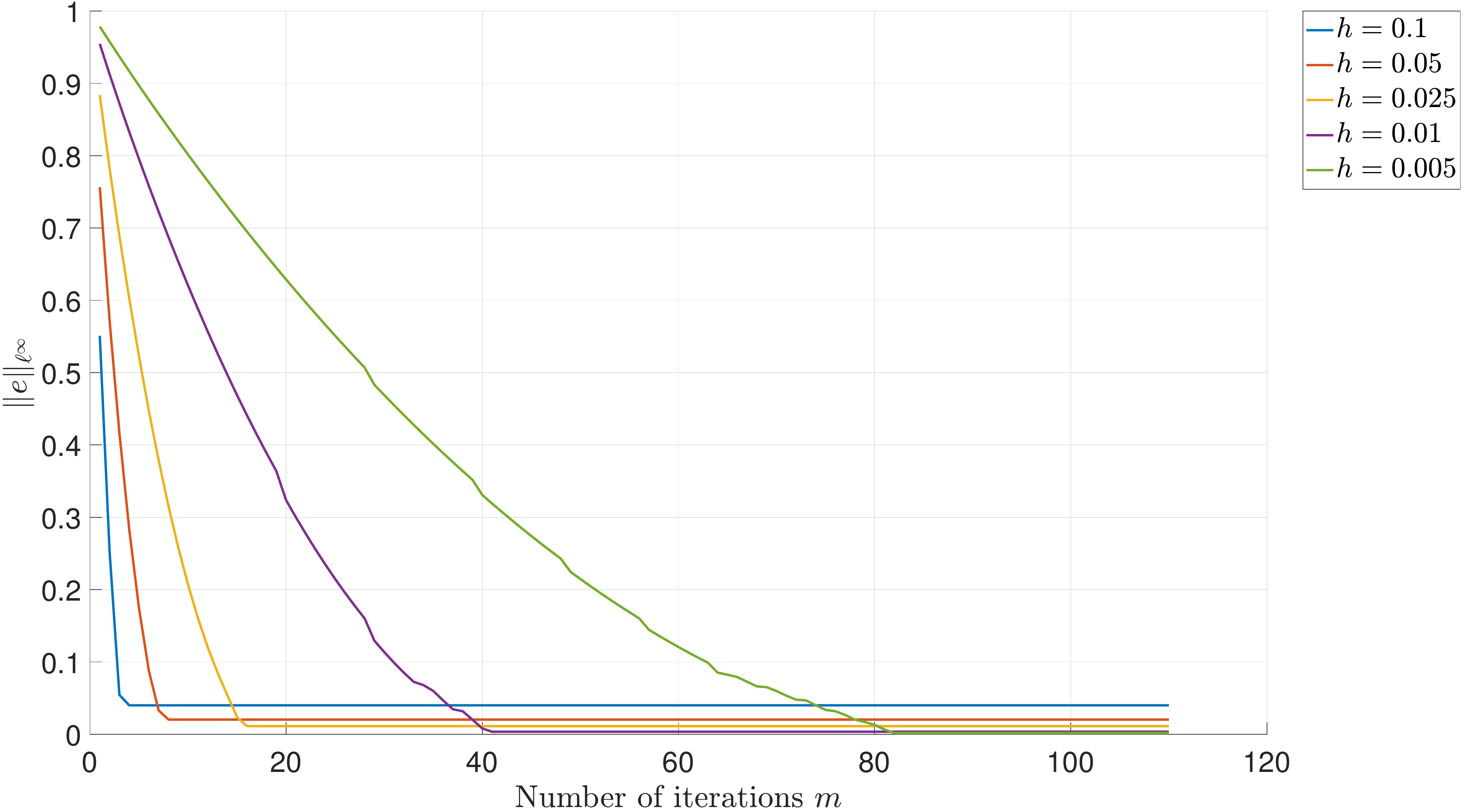}&
		\includegraphics[width=0.45\textwidth]{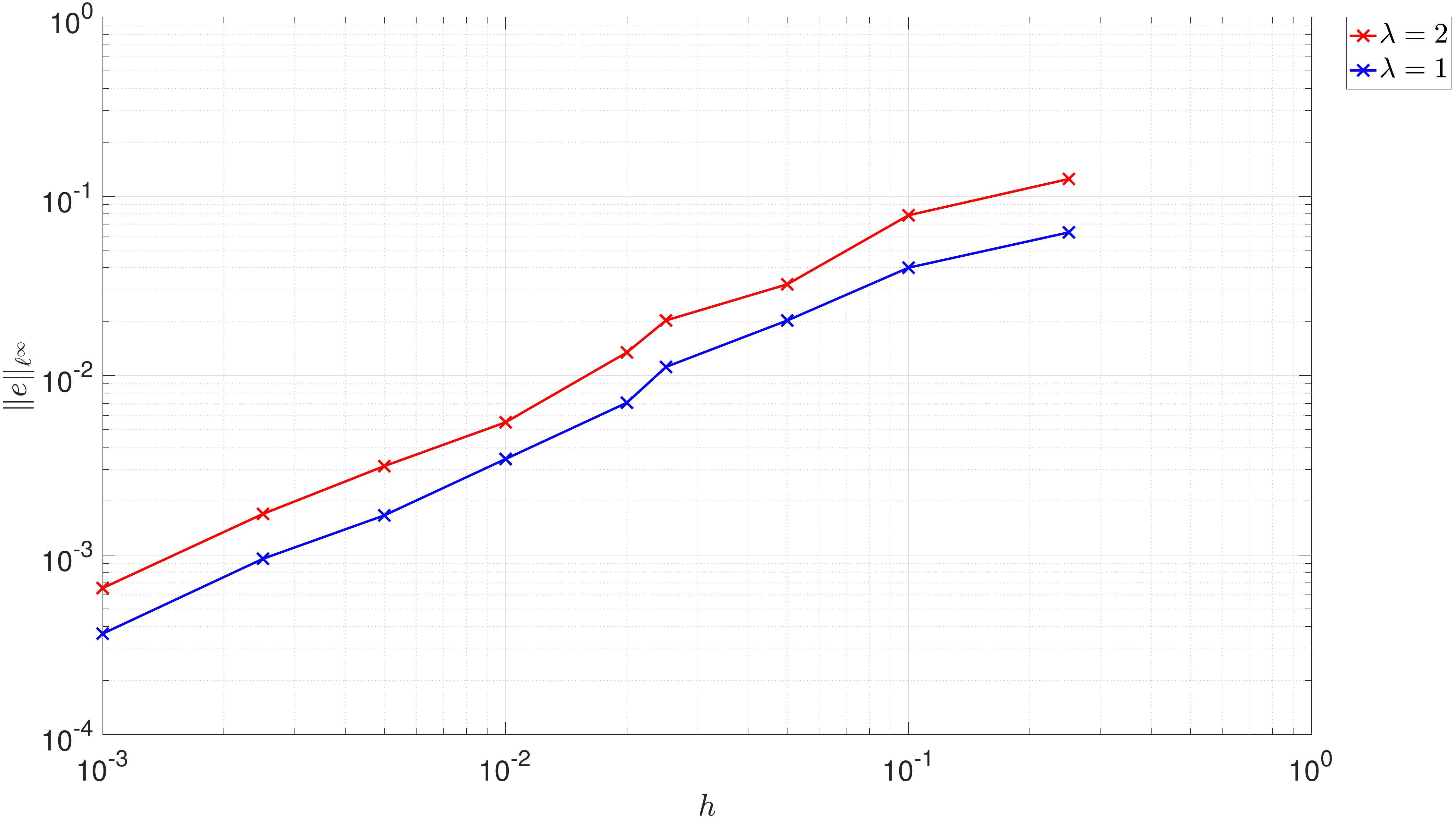}\\
		$(a)$ & $(b)$
	\end{array}$
\end{center}
   \caption{\label{Ex02.2d.CnvrgPlots}
$(a)$ Variation of the $\ell^{\infty}$ norm of the error versus the number of iterations $m$ for different 
grid sizes $h$ and $\lambda=1$.
 $(b)$ Convergence plot of the error with respect to the grid size and for different values of $\lambda$.
}
\end{figure}

As in the previous example, Figure \ref{Ex02.2d.CnvrgPlots}$(a)$ displays the convergence plot of the 
$\ell^{\infty}$ norm of the error
versus the number of iterations for different grid sizes, with a number of iterations 
that increases by reducing $h$, whereas Figure \ref{Ex02.2d.CnvrgPlots}$(b)$ shows the linear convergence of the 
error with the grid size.  Table \ref{Ex01.Tab:2d} contains for different grid sizes $h$,
the total number of iterations $m$ and the error $\|e\|_{\ell^{\infty}}$
of the approximations of $M^{\lambda,h}_{\Omega}(M^{h}_{\lambda,\Omega}(f_{\overline\Omega}^-))$ and 
$C_{\lambda}^l(f_{\R^n}^{-})$, for the comparison of the Moreau and 
convex based computation of the local lower transform, respectively. Also  observe that 
the number of iterations needed to compute the local lower transform using the Moreau based definition is much  
less than the one needed for the computation based on the convex envelope. For small grid size,
while the computation of $M^{\lambda,h}_{\Omega}(M^{h}_{\lambda,\Omega}(f_{\overline\Omega}^-))$ takes just a few seconds, 
the computing time of $C_{\lambda}^l(f_{\R^n}^{-})$ using a tolerance for the convergence of the scheme equal to $10^{-7}$ 
for the $\ell^{\infty}$ norm of the difference between two succesive iterates, is much longer. Furthermore,
for a given grid size, the discrete Moreau based lower transform is much more accurate than the convex based
lower transform and, as in the previous example, for small values of $h$ the convex based scheme
appears  to be unstable, showing oscillations in the error.

\begin{table}[tbhp]
\centerline{
\begin{tabular}{c|c|c|c|c|c|c|c|c|}
\cline{2-9}
	& \multicolumn{4}{c|}{$\lambda=1$} & \multicolumn{4}{c|}{$\lambda=2$} \\ \cline{2-9}
& \multicolumn{2}{c|}{ $M^{\lambda,h}_{\Omega}(M^{h}_{\lambda,\Omega}(f_{\overline\Omega}^-))$}    & \multicolumn{2}{c|}{ $C_{\lambda}^l(f_{\R^n}^{-})$ } & \multicolumn{2}{c|}{ $M^{\lambda,h}_{\Omega}(M^{h}_{\lambda,\Omega}(f_{\overline\Omega}^-))$}    & \multicolumn{2}{c|}{ $C_{\lambda}^l(f_{\R^n}^{-})$ } \\ \hline
\multicolumn{1}{ |l| }{$h$}	& $m$      &  $\|e\|_{\ell^{\infty}}$   & $m$     &  $\|e\|_{\ell^{\infty}}$  & $m$      &  $\|e\|_{\ell^{\infty}}$   & $m$     &  $\|e\|_{\ell^{\infty}}$ \\ \hline
\multicolumn{1}{ |l| }{$0.25$}	& $7$	   &  $0.0629657$		& $22$	  &  $0.1147375$ & $6$  &  $0.125$	& $18$   &  $0.2127817$	\\ \hline
\multicolumn{1}{ |l| }{$0.1$}   & $12$	   &  $0.04$			& $127$	  &  $0.0450623$ & $10$ &  $0.0784089$	& $94$   &  $0.0805668$	\\ \hline
\multicolumn{1}{ |l| }{$0.05$}	& $22$	   &  $0.0203122$		& $466$	  &  $0.0214302$ & $17$ &  $0.0323148$	& $340$  &  $0.0324987$	\\ \hline
\multicolumn{1}{ |l| }{$0.025$}	& $42$	   &  $0.0111877$		& $1634$  &  $0.0166887$ & $31$ &  $0.0203275$	& $1152$ &  $0.0228483$	\\ \hline
\multicolumn{1}{ |l| }{$0.02$}	& $52$	   &  $0.0070677$		& $2453$  &  $0.0170046$ & $39$ &  $0.0134798$	& $1742$ &  $0.0237278$	\\ \hline
\multicolumn{1}{ |l| }{$0.01$}	& $102$	   &  $0.0034324$		& $8073$  &  $0.0161028$ & $75$  &  $0.0054984$	& $5949$ &  $0.0219233$	\\ \hline
\multicolumn{1}{ |l| }{$0.005$}	& $202$	   &  $0.001664$		& $-$     &  $-$	 & $149$ &  $0.0031323$	& $-$    &  $-$	\\ \hline
\multicolumn{1}{ |l| }{$0.0025$}& $402$	   &  $0.00095$			& $-$	  &  $-$	 & $295$ &  $0.0016961$	& $-$	 &  $-$	\\ \hline
\multicolumn{1}{ |l| }{$0.001$} & $1002$   &  $0.0003644$		& $-$	  &  $-$	 & $735$ &  $0.0006545$	& $-$	 &  $-$	\\ \hline
\end {tabular}
}
\caption{\label{Ex01.Tab:2d} Number of iterations $m$ and values of the error $\|e\|_{\ell^{\infty}}$ 
of approximations of 
$M^{\lambda,h}_{\Omega}(M^{h}_{\lambda,\Omega}(f_{\overline\Omega}^-))$ and $C_{\lambda}^l(f_{\R^n}^{-})$ 
by applying Algorithm \ref{Algo:MoreauEnv} and 
Algorithm \ref{Algo:CnvxEnv}, respectively, for different values of the grid size $h$.
The number of iterations $m$ given in the table, has been obtained by taking as termination criteria of the 
two algorithms the check on the $\ell^{\infty}$ norm of the error between two successive iterates, which was
set not greater than $10^{-7}$.  
The values of $m$ shown in the table denote the total number of iterations. There are no values 
reported for the convex based scheme in the case of $h\leq 0.005$ due to the long running time.}
\end{table}


We conclude the discussion of this example with some observations on the computation of $C_{\lambda}^l(f^{\infty})$.
By means of numerical simulations, we show that there exists a domain $D_e$ that contains $\Omega$
such that $C_{\lambda}^l(f^{\infty})(x)=M^{\lambda}_{D_e}(M_{\lambda,D_e}(f^M))(x)$ for $x\in\Omega$ with $f^M$ defined below. 
Assume $a>0$ and consider the extended domain 

\[
	D_e=[-2-a,\,2+a]\,\times\,[-2-a,\,2+a]\,.
\]
Given $M>0$, define the following auxiliary function

\begin{equation}
	f^M(x)=\left\{\begin{array}{ll}
			\displaystyle f(x),&\displaystyle x\in B(O;\,2),\\[1.5ex]
			\displaystyle M,   &\displaystyle x\in D_e\setminus B(O;\,2)\,.
		\end{array}\right.
\end{equation}

Then we apply Algorithm \ref{Algo:MoreauEnv} to the interior grid points of $D_e$ to compute the lower and 
upper Moreau envelope. Table \ref{Ex01.Tab:2d.Infty} reports the value of the $\ell^{\infty}$ norm of the error between
$M^{\lambda}_{D_e}(M_{\lambda,D_e}(f^M))(x)$ and $C_{\lambda}^l(f^{\infty})(x)$ respectively in $\Omega$ for different values 
of the grid size $h$, the  domain extension parameter $'a'$ and the parameter $\lambda$. We observe that by choosing $M$ large enough
and by a suitable choice of $a$, we get an excellent agreement in $\Omega$. Also in this case, 
we conjecture that this extension depends on $\lambda$, but obtaining a formula for it is an open issue at present.
\begin{table}[tbhp]
\centerline{
\begin{tabular}{|c|c|c|c|c|c|c|c|}\hline
\multicolumn{4}{|c|}{ $\lambda=1$}&\multicolumn{4}{|c|}{ $\lambda=2$}\\ \hline
\multicolumn{2}{|c|}{ $h=0.01$}	& \multicolumn{2}{c|}{ { $h=0.005$} }&\multicolumn{2}{|c|}{ $h=0.01$}			      & \multicolumn{2}{c|}{ { $h=0.005$} } \\ \hline
 $a$  & $\|e\|_{\ell^{\infty}}$ & $a$ &  $\|e\|_{\ell^{\infty}}$ & $a$ &  $\|e\|_{\ell^{\infty}}$  & $a$       &  $\|e\|_{\ell^{\infty}}$ \\ \hline
$0$	 & $0.4960309$	& $0$	 & $0.4982693$ & $0$	 & $0.3302767$ & $0$	 & $0.3319488$	\\ \hline
$0.2$	 & $0.3165244$	& $0.2$	 & $0.3184523$ & $0.2$	 & $0.1180116$ & $0.2$	 & $0.1190902$	\\ \hline
$0.5$	 & $0.1225861$	& $0.5$  & $0.1239016$ & $0.34$	 & $0.0329538$ & $0.34$	 & $0.0336187$	\\ \hline
$1$	 & $0.0002370$  & $1$    & $0.0000963$ & $0.5$	 & $0.0004378$ & $0.5$	 & $0.0001693$	\\ \hline
$1.2$	 & $0.0002370$  & $1.2$	 & $0.0000963$ & $1$	 & $0.0004378$ & $1$	 & $0.0001693$	\\ \hline
\end {tabular}
}
\caption{\label{Ex01.Tab:2d.Infty} Values of the error $\|e\|_{\ell^{\infty}}$ 
in $\Omega$ between $C_{\lambda}^l(f^{\infty})$ and $M^{\lambda}_{D_e}(M_{\lambda,D_e}(f^M))$ for 
different values of $h$, the extension $a$ and $\lambda$. The results refer to $M\geq 10^3$.
}
\end{table}

\end{ex}


\begin{ex}\label{Se5.NumEx.MMAM}\textbf{The Multiscale Medial Axis Map.}
We present an application of Corollary \ref{Cor.SetCmp} to find the multiscale medial axis map of 
the closed set $K$ represented in Figure \ref{Sec5.NumEx.MMAM}$(a)$. 
The open set $\Omega$ is taken, in this case, as the domain of the whole image and $A=\Omega\setminus K$. 
By Corollary \ref{Cor.SetCmp} the quadratic multiscale medial axis map with scale $\lambda$
of the closed set $K$ can then be computed for $x\in\overline\Omega$ as 

\begin{equation*}
	\mathcal{M}(\lambda;\,K)(x)=\Big(f_{\overline\Omega}^{-}(x)-C_{\lambda,\Omega}^l(f_{\overline\Omega}^{-})(x)\Big)
\end{equation*}
where $f_{\overline\Omega}^{-}(x)=\dist^2(x;\,K^c)$ for $x\in\overline\Omega$. Note that in this case 
$f_{\overline\Omega}^{-}(x)=0$ for $x\in\partial\Omega$. 

\begin{figure}[htbp]
\centerline{
	$\begin{array}{cc}
		\fbox{\includegraphics[width=0.45\textwidth]{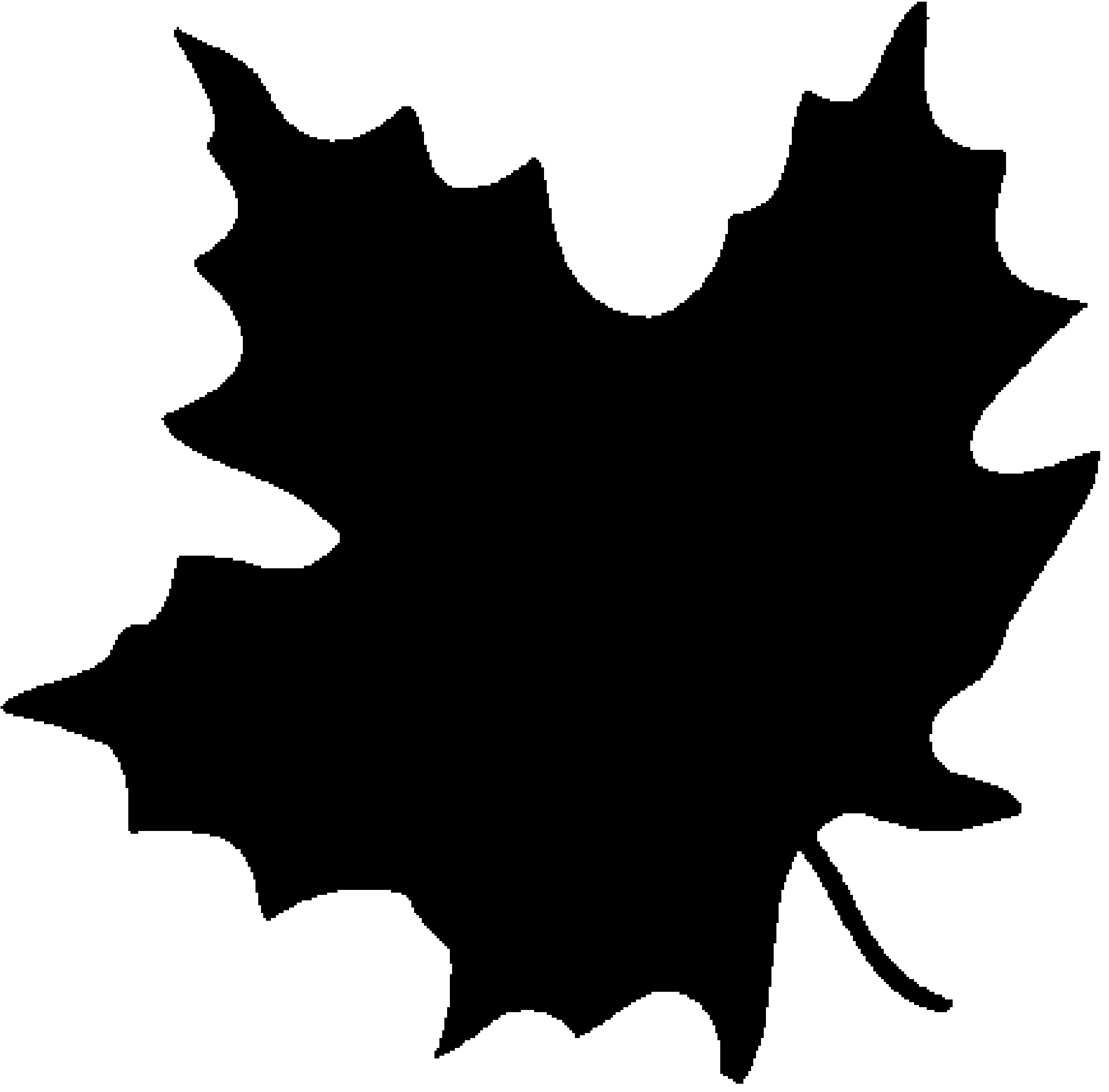}}&
		\includegraphics[width=0.45\textwidth]{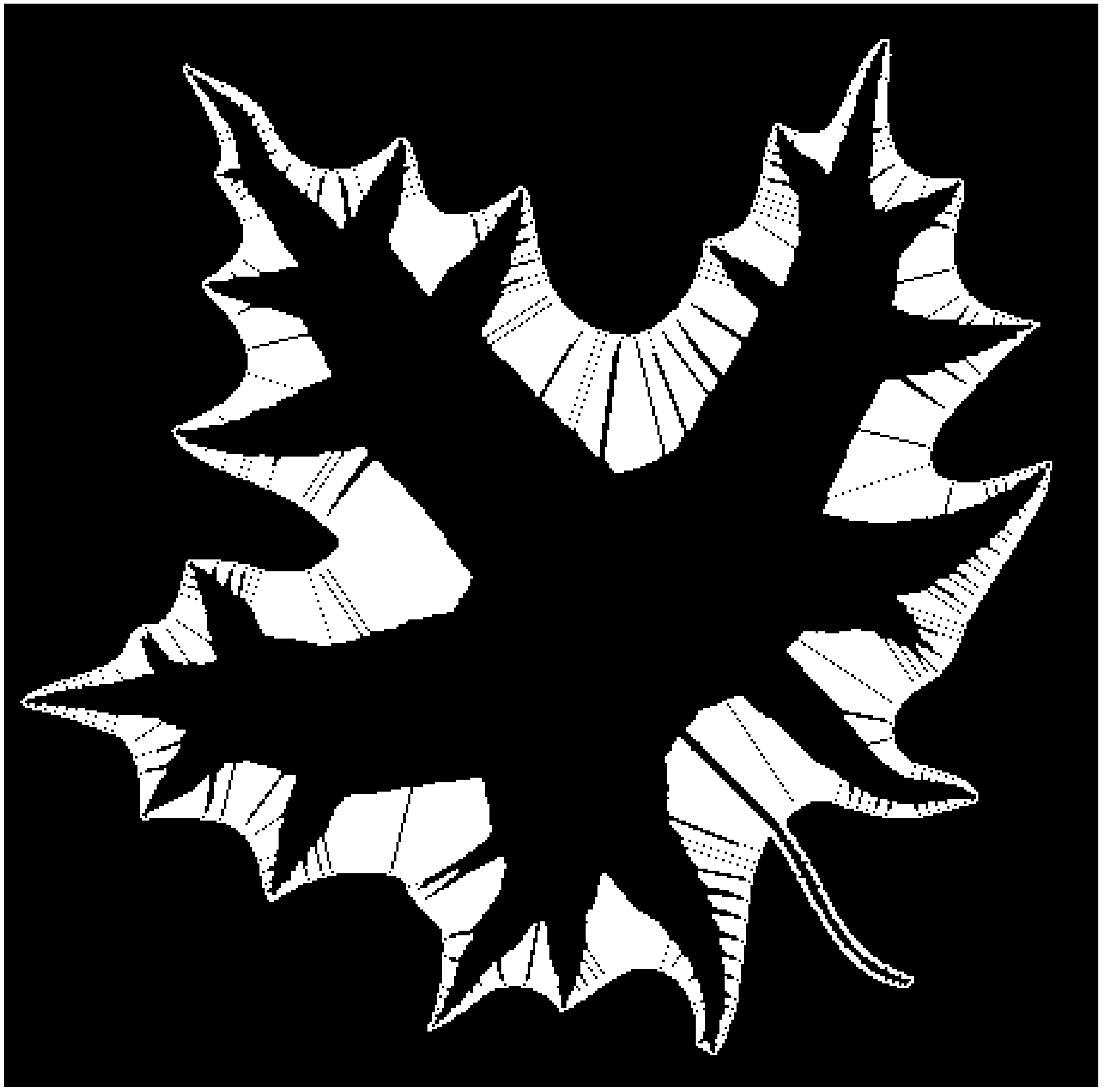}\\
			(a)&(b)	\\
		\includegraphics[width=0.45\textwidth]{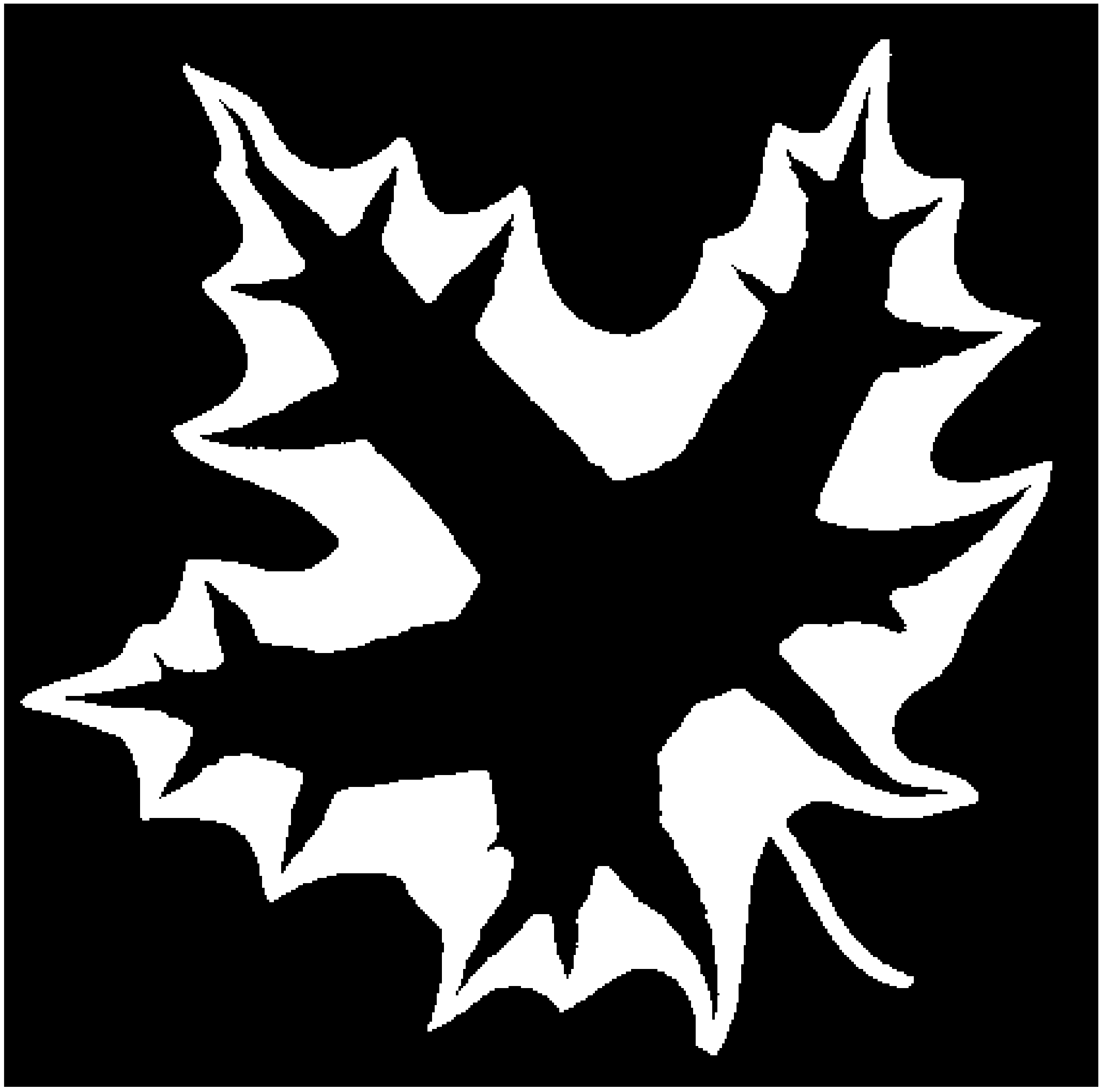}&
		\includegraphics[width=0.45\textwidth]{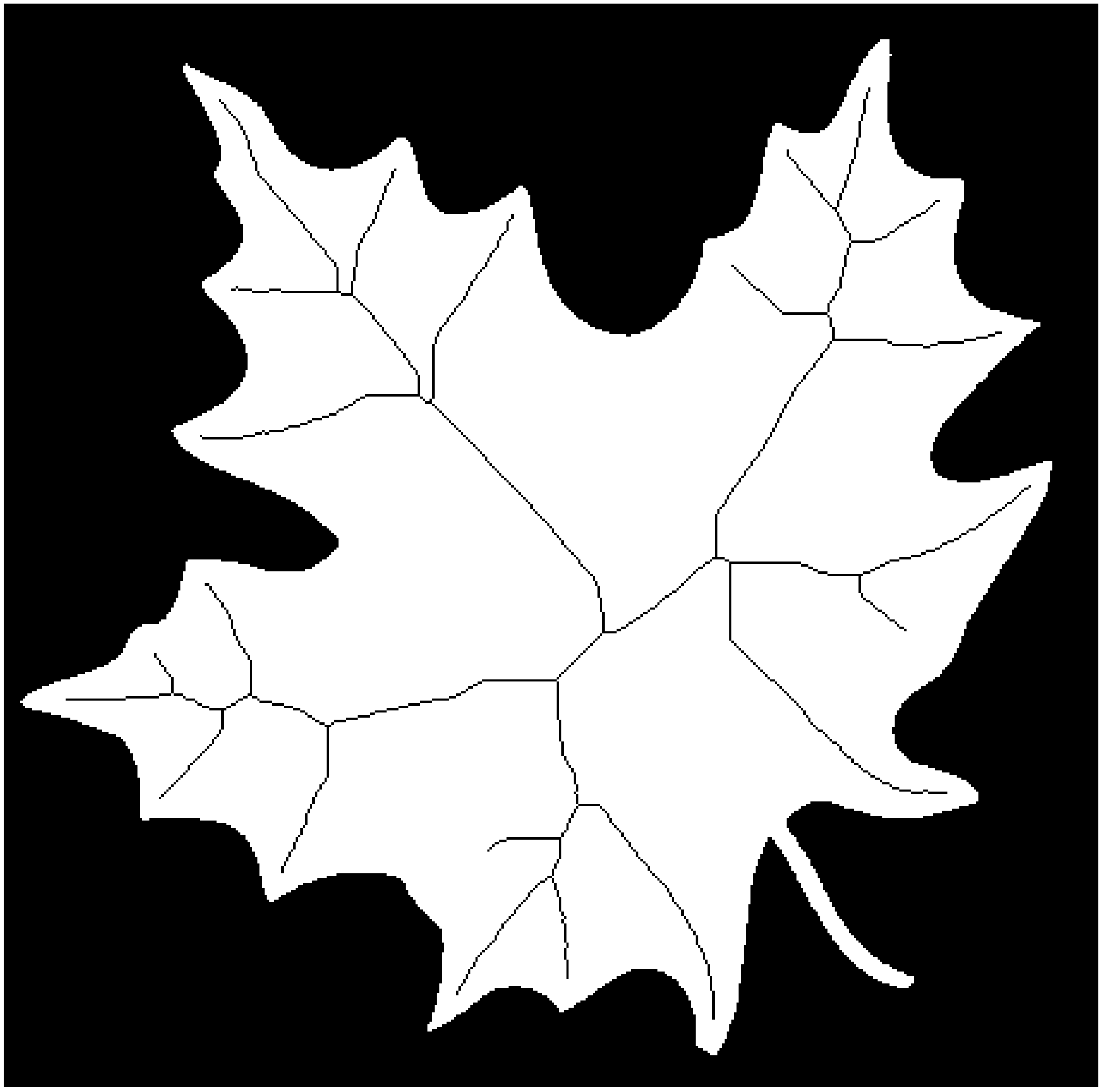}\\
			(c)&(d)
	\end{array}$
}
\caption{\label{Sec5.NumEx.MMAM} Example \ref{Se5.NumEx.MMAM}.  
Computation of the multiscale medial axis map. 
$(a)$ Characteristic function of the set $\Omega\setminus K$;
$(b)$ Support of the multiscale medial axis map for $\lambda=1$ with the display of all the fine brances generated 
by the steps on the boundary. $(c)$ Suplevel set $\{x\in\Omega: \mathcal{M}(\lambda;\,K)(x)> t\}$ for $t=1$, displaying 
the stable branches corresponding to the level $t$. $(d)$ Morphological thinning of $(c)$ by the structuring element 
described in \cite[page 879, bottom of first column through top of second column]{LLS92} implemented in \textsc{MATLAB}.
}
\end{figure}

The support of $\mathcal{M}(\lambda;\,K)(x)$ with all its fine branches is shown
in Figure \ref{Sec5.NumEx.MMAM}$(b)$. Figure \ref{Sec5.NumEx.MMAM}$(c)$ depicts the stable branches 
that correspond to the suplevel set of $\mathcal{M}(\lambda;\,K)(x)$, given by 
$\{x\in\Omega: \mathcal{M}(\lambda;\,K)(x)> t\}$ with $t>0$ measure of the branch height, 
whereas Figure \ref{Sec5.NumEx.MMAM}$(d)$ displays the results of the morphological thinning of the 
suplevel set shown in Figure \ref{Sec5.NumEx.MMAM}$(c)$ using the algorithm described in 
\cite[page 879, bottom of first column through top of second column]{LLS92} implemented in \textsc{MATLAB}.
To produce Figure \ref{Sec5.NumEx.MMAM}$(c)$ we have taken $t=1$ after normalizing $\mathcal{M}(\lambda;\,K)$
to the range $[0,\,255]$. While the medial representation of the leaf stem is present in the support of $\mathcal{M}(\lambda;\,K)$ 
(see Figure \ref{Sec5.NumEx.MMAM}$(b)$), this disappears in the suplevel set relative to $t=1$. 
The very small values of $\mathcal{M}(\lambda;\,K)$ at such points is the result of the 
small value of the separation angle, which is, in turn, related to 
the values of $\mathcal{M}(\lambda;\,K)$ (see the bound (3.11) in \cite{ZCO15c}). 
The application of Algorithm \ref{Algo:MoreauEnv} and Algorithm \ref{Algo:CnvxEnv}, when $f$ is an image,
simplifies by taking the digitized image as the grid, and $h=1$ equal to the pixel size.
To compare the performance of the two algorithms, Figure \ref{Sec5.NumEx.MMAM-CnvrgPlot} displays the 
variation of the $L^2-$norm of $C_{\lambda,\Omega}^l(f_{\overline\Omega}^{-})$
with the number of iterations $m$ used to compute 
$M^{\lambda}_{\Omega}(M_{\lambda,\Omega}(f_{\overline\Omega}^{-}))$ by Algorithm \ref{Algo:MoreauEnv} 
and to compute $C_{\lambda}^l(f_{\R^n}^{-})$ by Algorithm \ref{Algo:CnvxEnv}.
By referring to the Moreau based definition,  
convergence of the scheme is achieved after a finite number of iterations ($m=118$) which is much lower than those needed to compute the 
convex envelope based definition of the lower transform ($m=22783$). 

\begin{figure}[htbp]
	\centerline{$\includegraphics[width=0.45\textwidth]{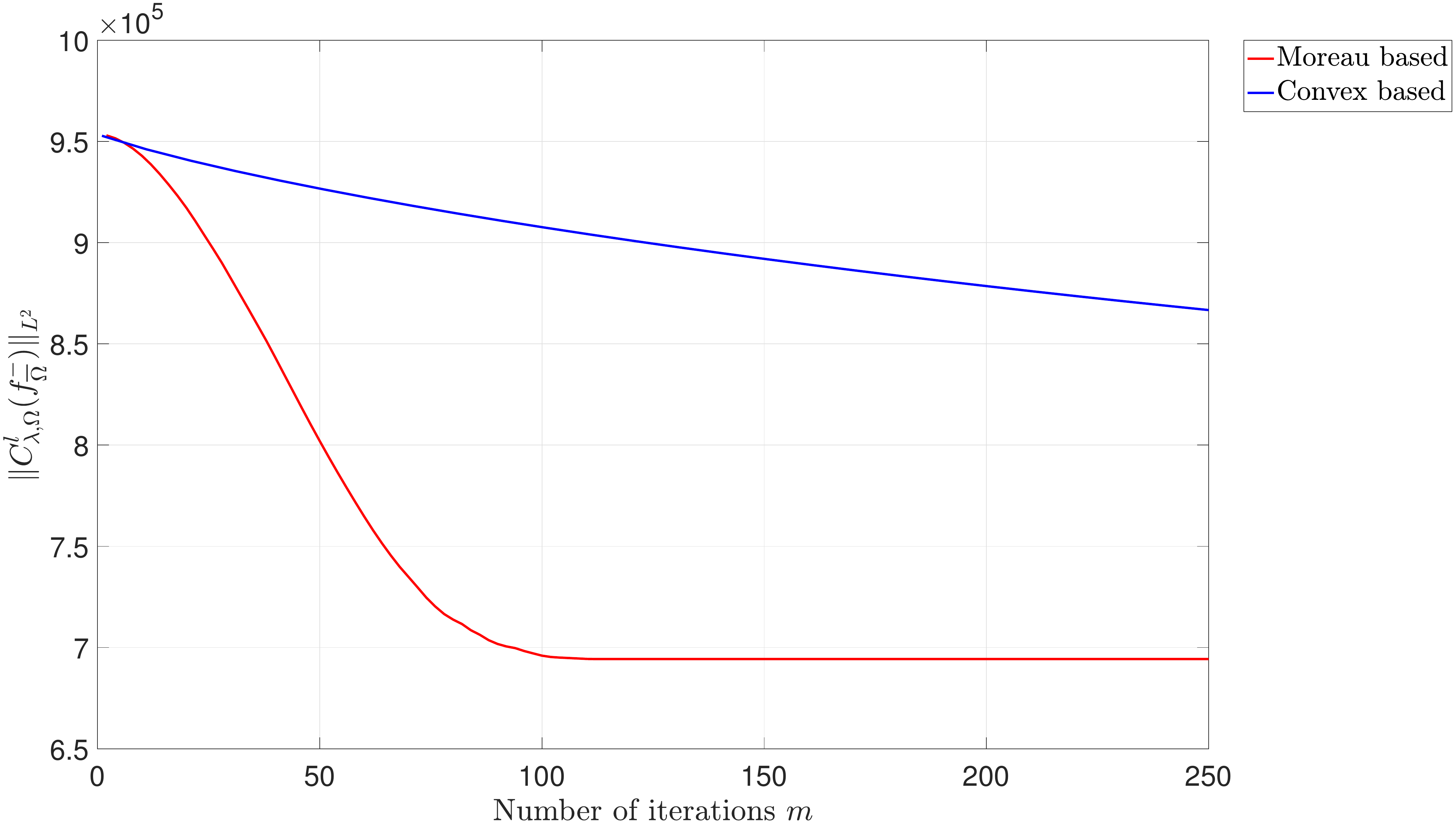}$}
\caption{\label{Sec5.NumEx.MMAM-CnvrgPlot}  
	Example \ref{Se5.NumEx.MMAM}.
	Variation of the $L^2$-norm of $C_{\lambda,\Omega}^l(f_{\overline\Omega}^{-})$ with the number $m$ of iterations 
	using the Moreau based definition and the convex based definition. $\lambda=1$.
}
\end{figure}

\end{ex}

\begin{ex}\label{Se5.NumEx.CurvInt}\textbf{Shape interrogation.}
As an application of Theorem \ref{Thm.3.ChrctFnct}, we consider  the computation of the 
intersection extraction filter $I_{\lambda}(\cdot;\,K)$ introduced in \cite{ZOC15a} with the digitized set  
$K$ as the input image. Given a non-empty compact set $K\subset \R^n$, 
and taking $\Omega$ as a reference bounding box such that $K\subset\Omega$ and with $K$ distant 
enough from the boundary of $\Omega$, by Theorem \ref{Thm.3.ChrctFnct}, the filter $I_{\lambda}(\cdot;\,K)$ 
can be expressed in terms of the local transforms as

\begin{equation}\label{Eq.Sec5.NumEx:filt}
		I_{\lambda,\Omega}(x;\,K)=\Big| C_{4\lambda,\Omega}^{u}(\chi_K^{\Omega})(x)-
		2\Big(C_{4\lambda,\Omega}^{u}(\chi_K^{\Omega})(x)-
		C_{\lambda,\Omega}^{l}(C_{\lambda,\Omega}^{u}(\chi_K^{\Omega}))(x)
		\Big)\Big|\,.
\end{equation}
For the digitized set $K$ given by the collection of curves shown in Figure \ref{Sec5.NumEx.ShapeInte}$(a)$, 
the local maxima of  $I_{\lambda,\Omega}(\cdot;\,K)$ coincide with all the crossing and turning points of 
the set $K$.
The filter defined by \eqref{Eq.Sec5.NumEx:filt} can also be applied to $3d$ geometries. Due to the Hausdorff stability 
of \eqref{Eq.Sec5.NumEx:filt}, we can also consider $K$ represented by point clouds. 
Figure \ref{Sec5.NumEx.ShapeInte}$(b)$ displays the intersection between manifolds of different 
dimensions with each manifold sampled by point clouds.

\begin{figure}[htbp]
\centerline{$\begin{array}{cc}
	\begin{overpic}[scale=0.45]{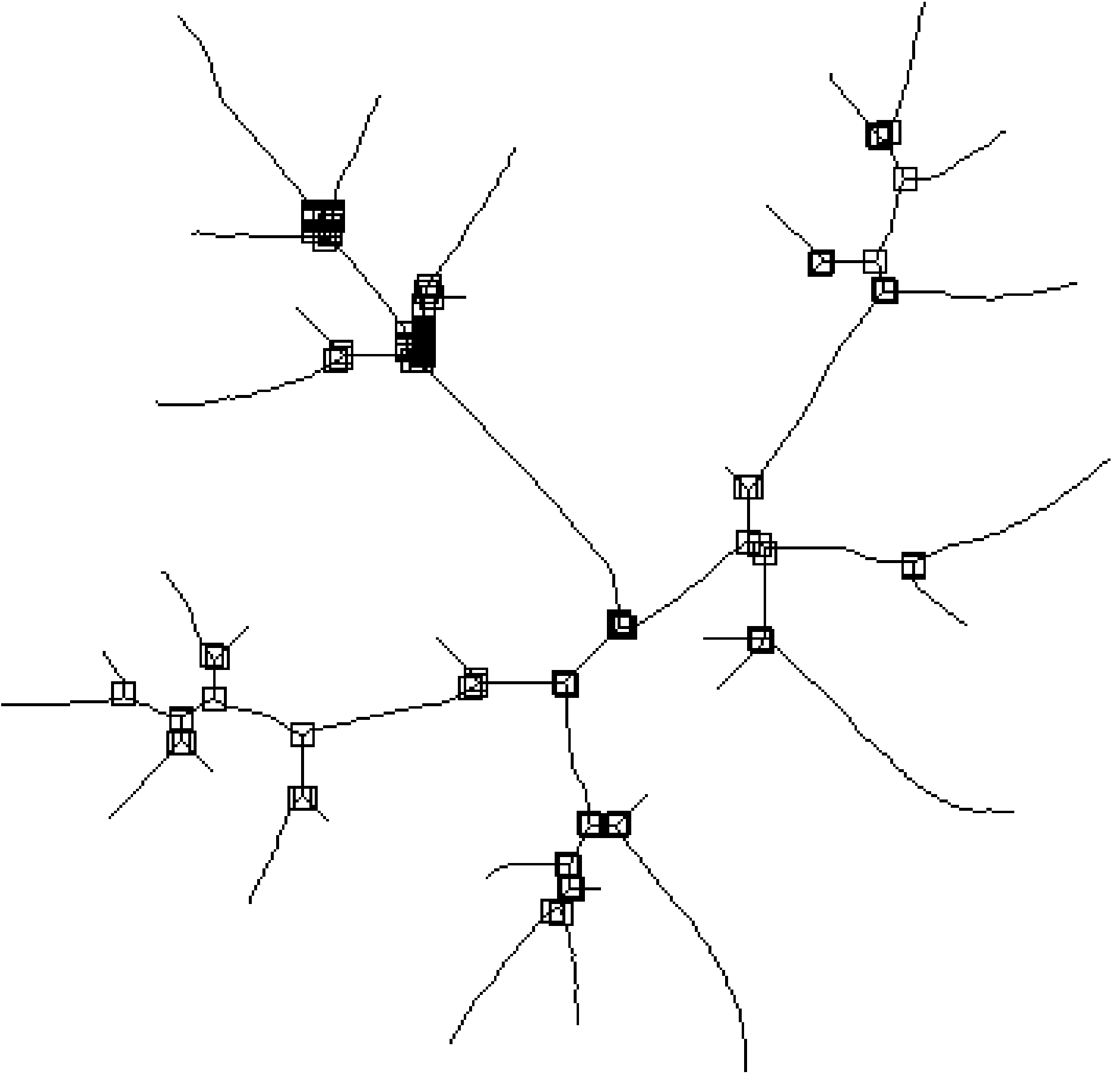}
		\put(80,1){\includegraphics[scale=0.15]{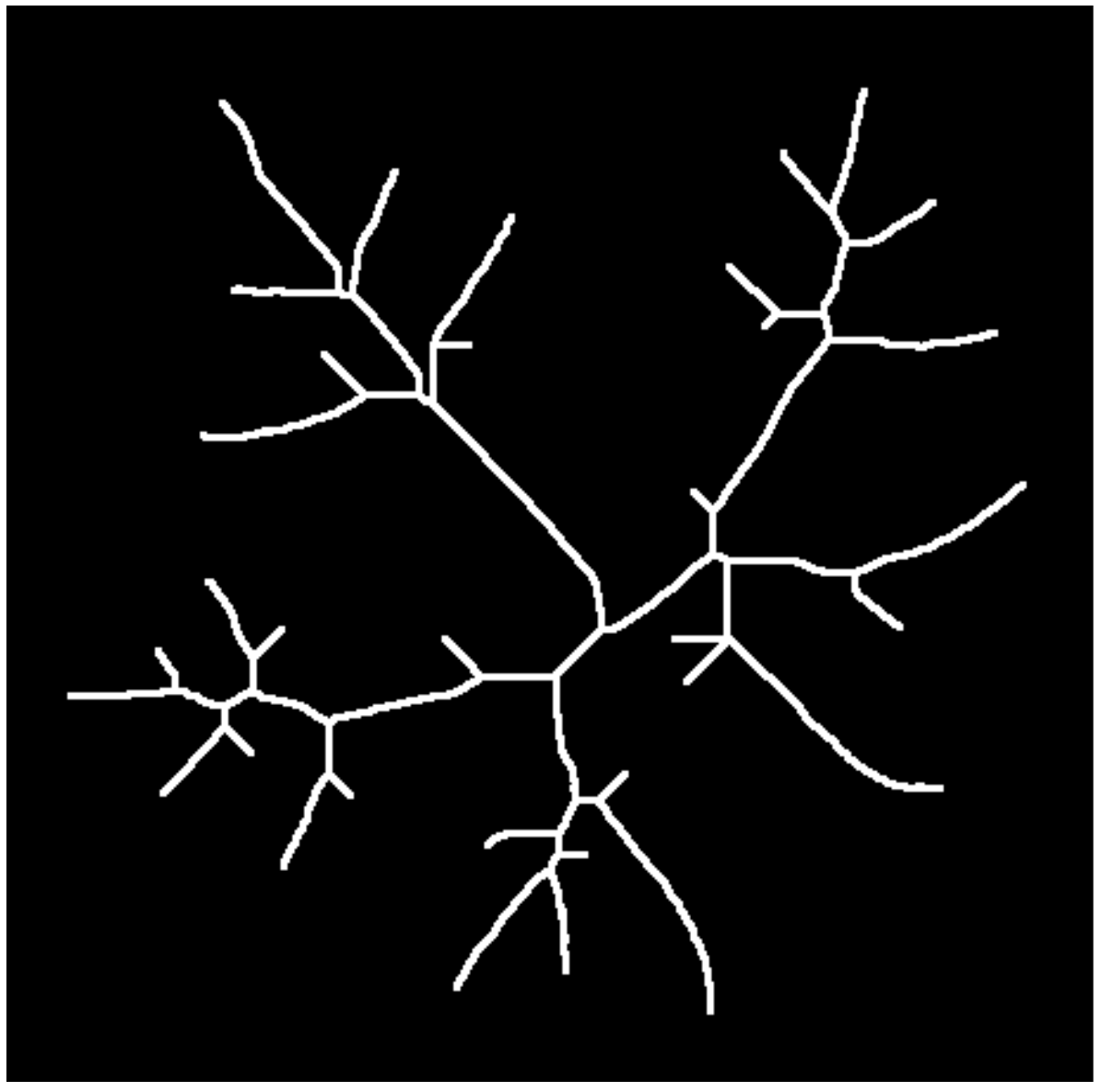}}
	\end{overpic}&
		\includegraphics[width=0.45\textwidth]{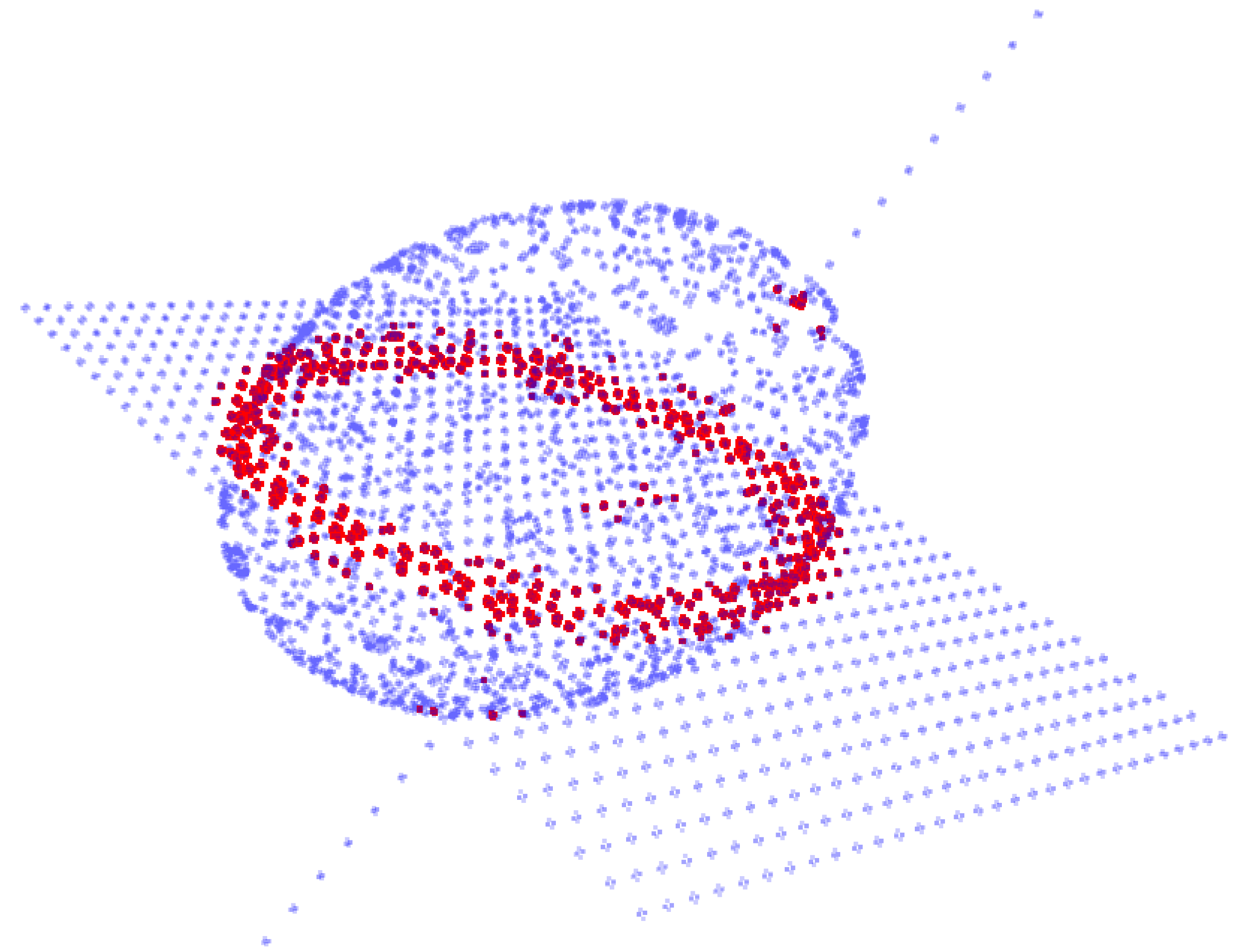}\\
	(a)&(b)
\end{array}$
}
\caption{\label{Sec5.NumEx.ShapeInte} Example \ref{Se5.NumEx.CurvInt}.
Applications of the filter defined by \eqref{Eq.Sec5.NumEx:filt} for shape interrogation to detect: $(a)$ Intersection points of curves with square markers at the
local maxima of the filter \eqref{Eq.Sec5.NumEx:filt}; $\lambda=5$. The network is displayed as inlaid picture.
$(b)$ Intersections between an ellipsoid, a plane and a line; $\lambda=0.01$.
}
\end{figure}

\end{ex}


 \begin{ex}\label{Sec5.ImInp}\textbf{Image Inpainting.}
 Let $\Omega$ be the domain of the whole image, $D\subset \Omega$ the set of missing/damaged pixels  	 
 and $K=\Omega\setminus D$ the set of the true pixels. The image inpainting problem consists in
 reconstructing the image over $D$ from knowing $f_K$, if we denote by $f$ the original image.
 In \cite{ZCO18} it is shown that the restored image can be obtained by the average compensated convex transform,
 which under the assumptions of Theorem \ref{Thm.3.ImpNew}, can be computed for $x\in\overline\Omega$ as

 \begin{equation}\label{Sec5.NumEx.Eq:Imp}
 		A_{\lambda,\,\Omega}^M(x)=\frac{1}{2}\Big(C_{\lambda,\Omega}^l(f_{\overline\Omega,K}^M)(x)+
 			C_{\lambda,\Omega}^u(f_{\overline\Omega,K}^{-M})(x)\Big)
 		\,,
 \end{equation}
 where $f_{\overline\Omega,K}^{M}$ and $f_{\overline\Omega,K}^{-M}$ have been defined in \eqref{Def.01.AuxFnct.M.Loc1} and 
 \eqref{Def.01.AuxFnct.M.Loc2}, respectively.
 We give here an application of \eqref{Sec5.NumEx.Eq:Imp} for the problem of removing scratches overprinted over an image such as
 the one displayed in Figure \ref{Sec5.NumEx.Fig:Imp}$(a)$.
 By using the Moreau based definition of the upper and lower transforms that enter \eqref{Sec5.NumEx.Eq:Imp} 
 and by applying Algorithm \ref{Algo:MoreauEnv},
 we compute with $\lambda=5$ and $M=10^4$ the averaged image $A_{\lambda,\Omega}^{M}$ defined by \eqref{Sec5.NumEx.Eq:Imp} and consider 
 the restored image given by
 
 \begin{equation}\label{Sec5.NumEx.Eq:Imp.Rest}
 	I(x)=\left\{\begin{array}{ll}
 		A_{\lambda,\Omega}^M(x)& x\in D\\[1.5ex]
 		f_K(x) & x\in K.
 		\end{array}\right.
 \end{equation}
The restored image $I$ is shown in Figure \ref{Sec5.NumEx.Fig:Imp}$(b)$. The number of iterations needed to obtain 
$A_{\lambda,\Omega}^M$ is $m=34$ for a tolerance equal to $10^{-7}$ on the $\ell^{\infty}$ norm of the difference between two
succesive iterates. Taking the $\mathrm{PSNR}$ (peak-to-signal ratio) as a measure of the quality of the restoration,
which is expressed in the units of $\mathrm{dB}$ and, for an $8-$bit image,
is defined by

\begin{equation*}
	\mathrm{PSNR}=10\log_{10}\displaystyle \frac{255^2}{\frac{1}{mn}\sum_{i,j}|f_{i,j}-r_{i,j}|^2}
\end{equation*}
where $f_{i,j}$ and $r_{i,j}$ denote the pixels values of the original and restored image, respectively,
and $m,\,n$ denote the size of the image $f$,
we find a value of $\mathrm{PSNR}$ equal to $33.165\,\mathrm{dB}$. Figure \ref{Sec5.NumEx.Fig:Imp}$c$ 
displays the restored image based on the convex based transforms, where we use the convex based definition of the transforms and
apply Algorithm \ref{Algo:CnvxEnv} to compute the convex envelope with a $tol=10^{-7}$ for the $\ell^{\infty}$ norm of the difference
between two successive iterates. 
We assume $C_{\lambda}^l(f_K^M)= C_{\lambda}^u(f_K^{-M})=f_K$ on the boundary of the image array.
In this case, we needed $m=1063$ iterations and got $\mathrm{PSNR}=33.239\,\mathrm{dB}$, 
which is slightly higher than the one that uses the Moreau based 
transforms but at the expenses of an higher number of iterations. Figure \ref{Sec5.NumEx.Fig:Imp}$d$ 
displays the TV based restoration obtained by applying the TV inpainting method described in \cite{CHN05}
and solved by the split Bregman method described in \cite{CCM07,Get12}. In this case $m=3531$ and $\mathrm{PSNR}=33.027\,\mathrm{dB}$.
Figure \ref{Sec5.NumEx.Fig:Imp.Det}  compares the  details of the original image and 
of the restored images near the right eye, respectively, showing that $A_{\lambda,\Omega}^{M}$, either by the Moreau based transforms 
or by the Convex based transforms, is able to preserve image details and does not introduce unintended effects.   
 
 \begin{figure}[htbp]
 \centerline{$\begin{array}{cc}
 	\includegraphics[width=0.35\textwidth]{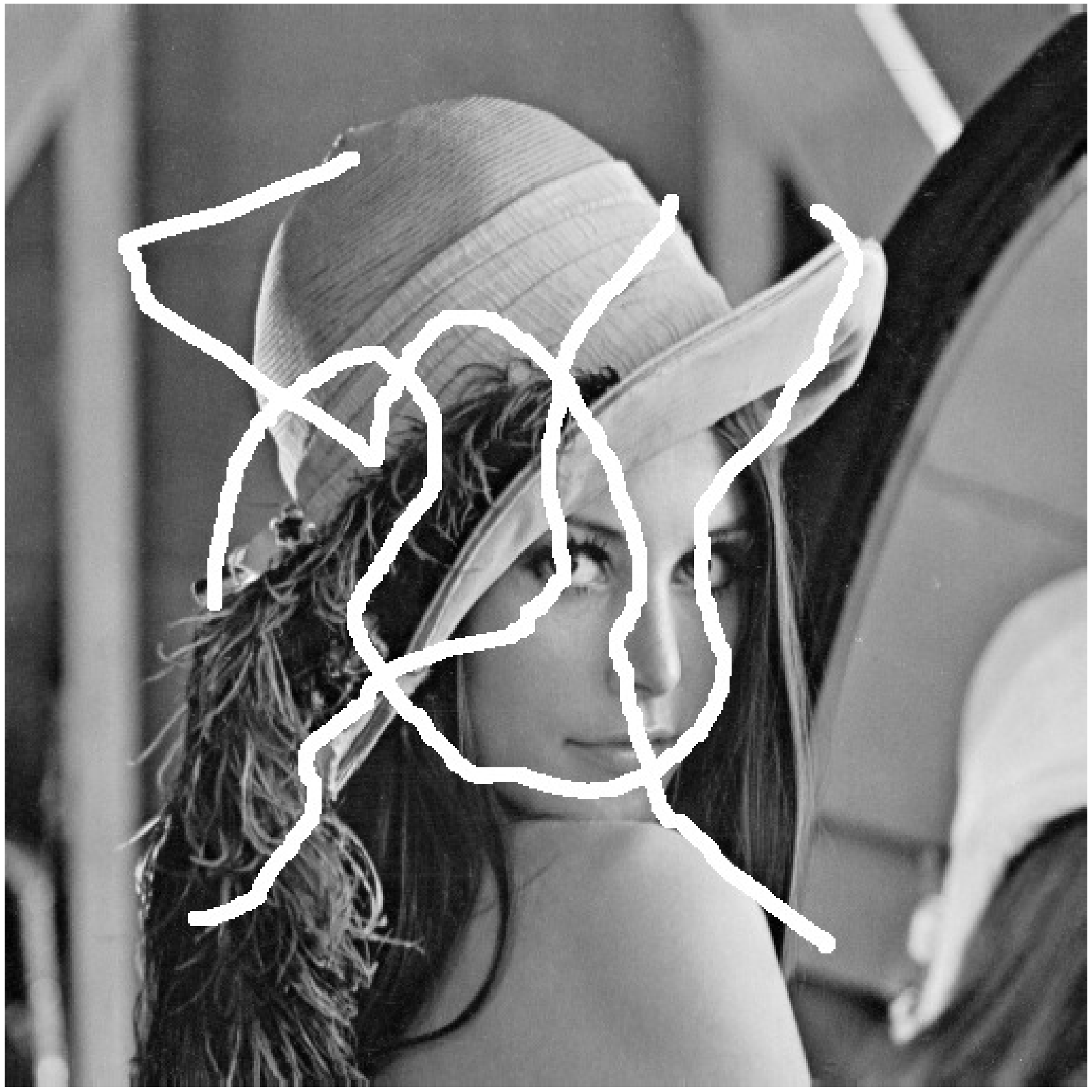}&
 	\includegraphics[width=0.35\textwidth]{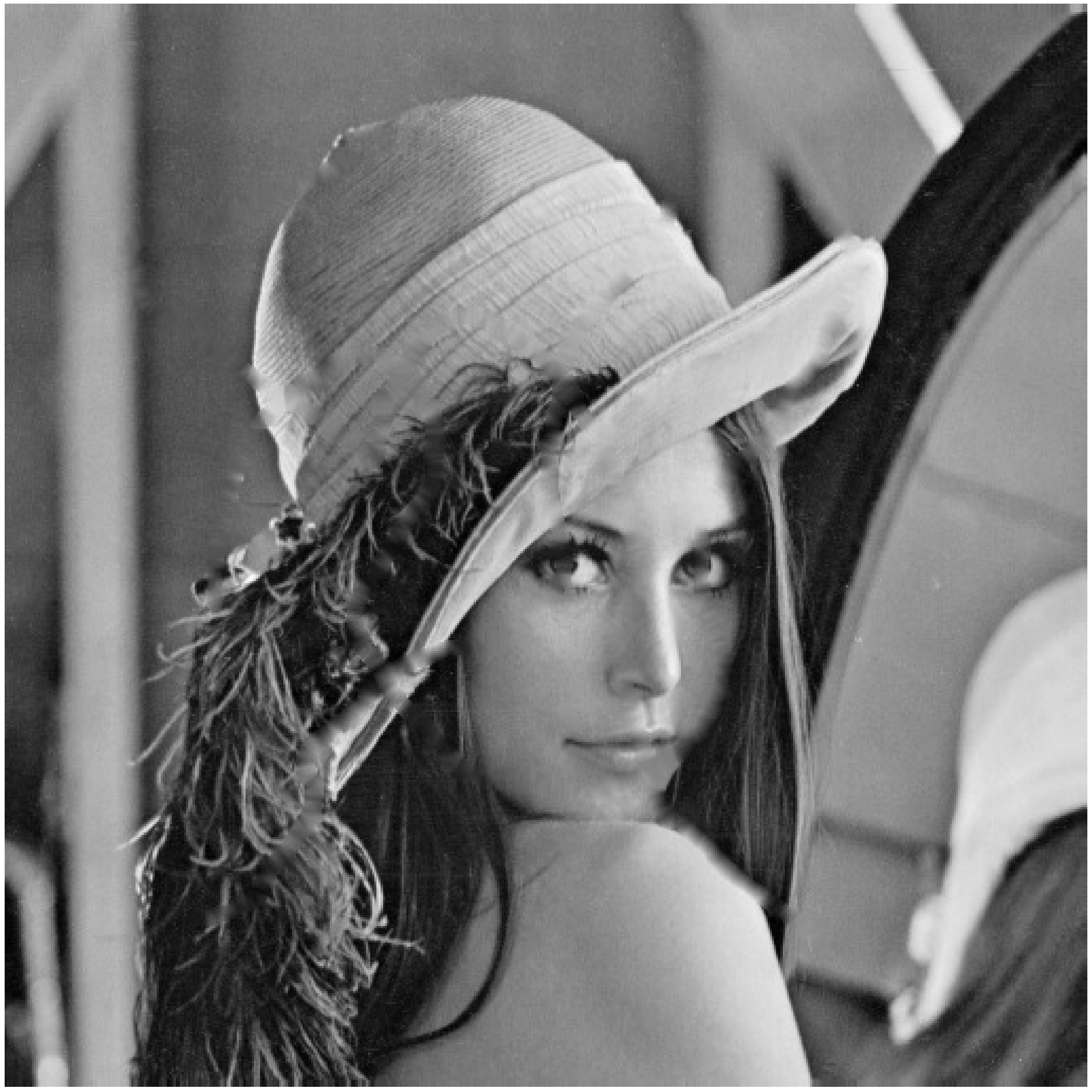}\\
 	(a) & (b) \\
 	\includegraphics[width=0.35\textwidth]{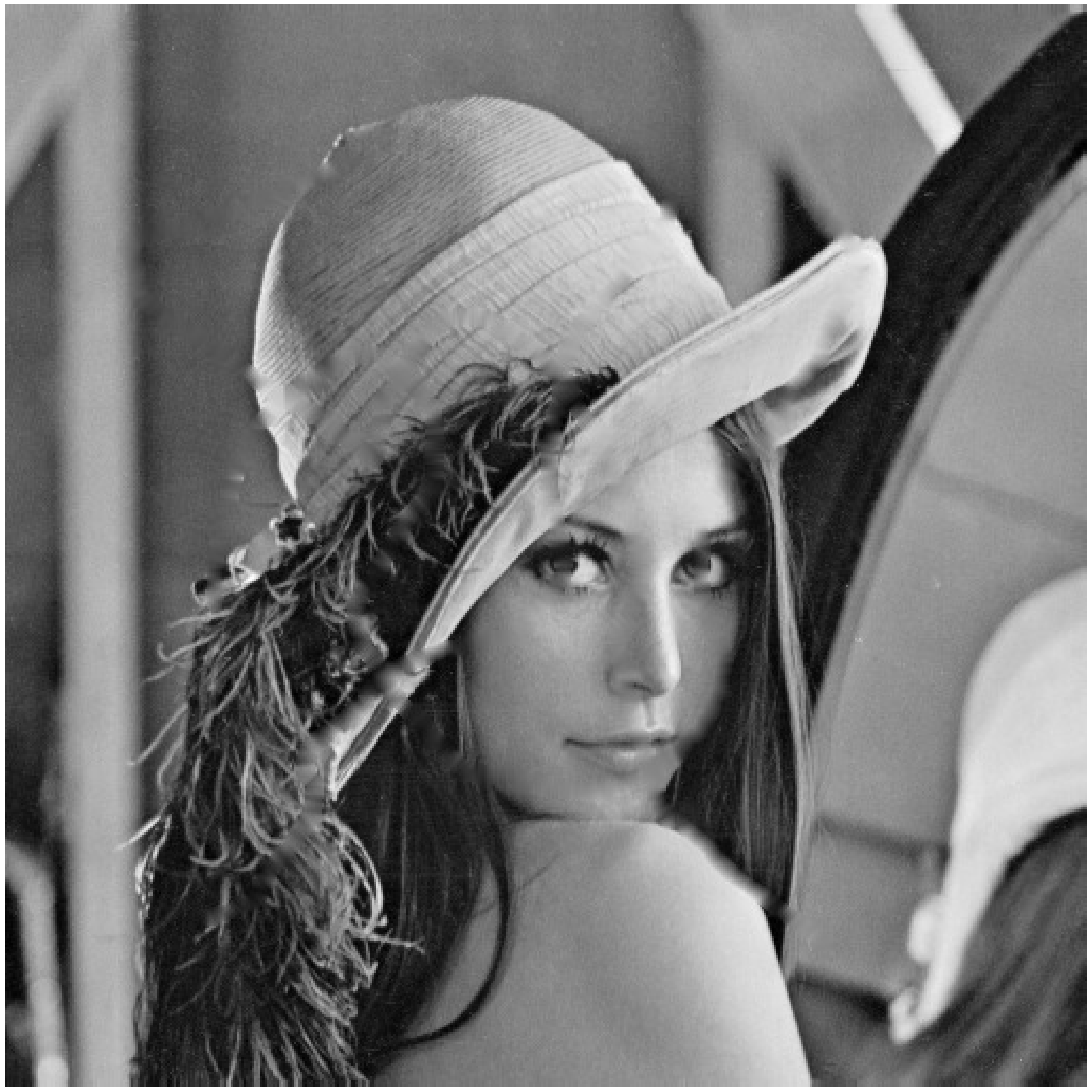}&
 	\includegraphics[width=0.35\textwidth]{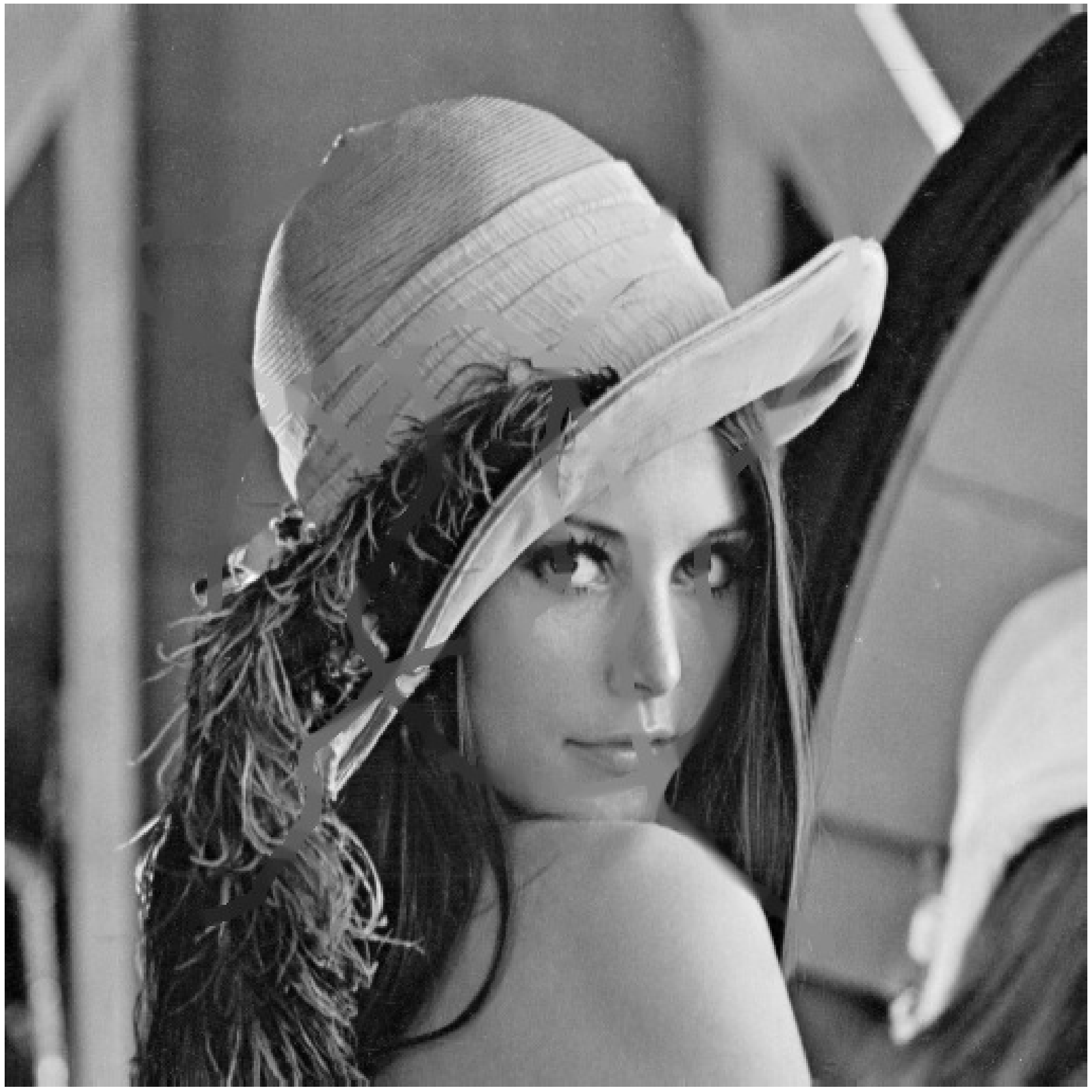}\\
 	(c) & (d)
 \end{array}$
 }\caption{\label{Sec5.NumEx.Fig:Imp} Example \ref{Sec5.ImInp}. 
 	Impainting of scratches over an image: $(a)$ \textsf{Lena} image with scratches; $(b)$ Restored image $I$ as defined by \eqref{Sec5.NumEx.Eq:Imp.Rest}
 	with $\lambda=5$ and $M=10^4$, with the Moreau based trasforms computed using Algorithm \ref{Algo:MoreauEnv}. 
 	Computed value for $\mathrm{PSNR}=33.165\,\mathrm{dB}$. Number of iterations $m=34$ for a tolerance on the error 
	between two succesive iterates equal to $10^{-7}$; $(c)$; Restored image $I$ as defined by \eqref{Sec5.NumEx.Eq:Imp.Rest}
 	with $\lambda=5$ and $M=10^4$, with the convex based trasforms computed using Algorithm \ref{Algo:CnvxEnv}. 
 	Computed value for $\mathrm{PSNR}=33.239\,\mathrm{dB}$. Number of iterations $m=1063$ 
	for a tolerance on the error between two succesive iterates equal to $10^{-7}$; 
 	$(d)$  Restored image by the Split Bregman inpainting method described in \cite{CCM07,Get12}. 
 		 Computed value for $\mathrm{PSNR}=33.027\,\mathrm{dB}$. Number of iterations $m=3531$.}
 \end{figure}
 
 \begin{figure}[ht]
 \centerline{$\begin{array}{ccc}
 	\includegraphics[width=0.33\textwidth]{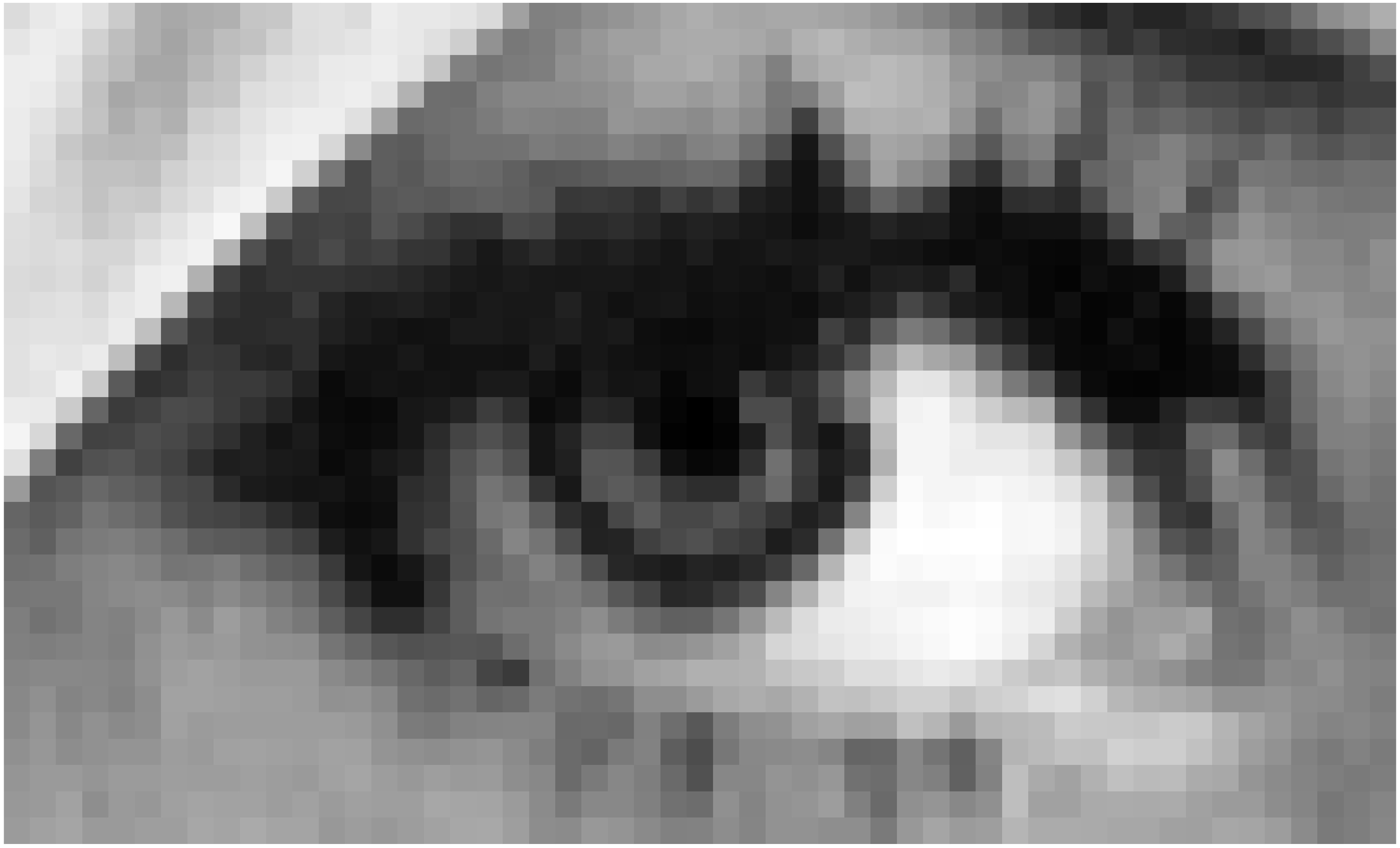}&
 	\includegraphics[width=0.33\textwidth]{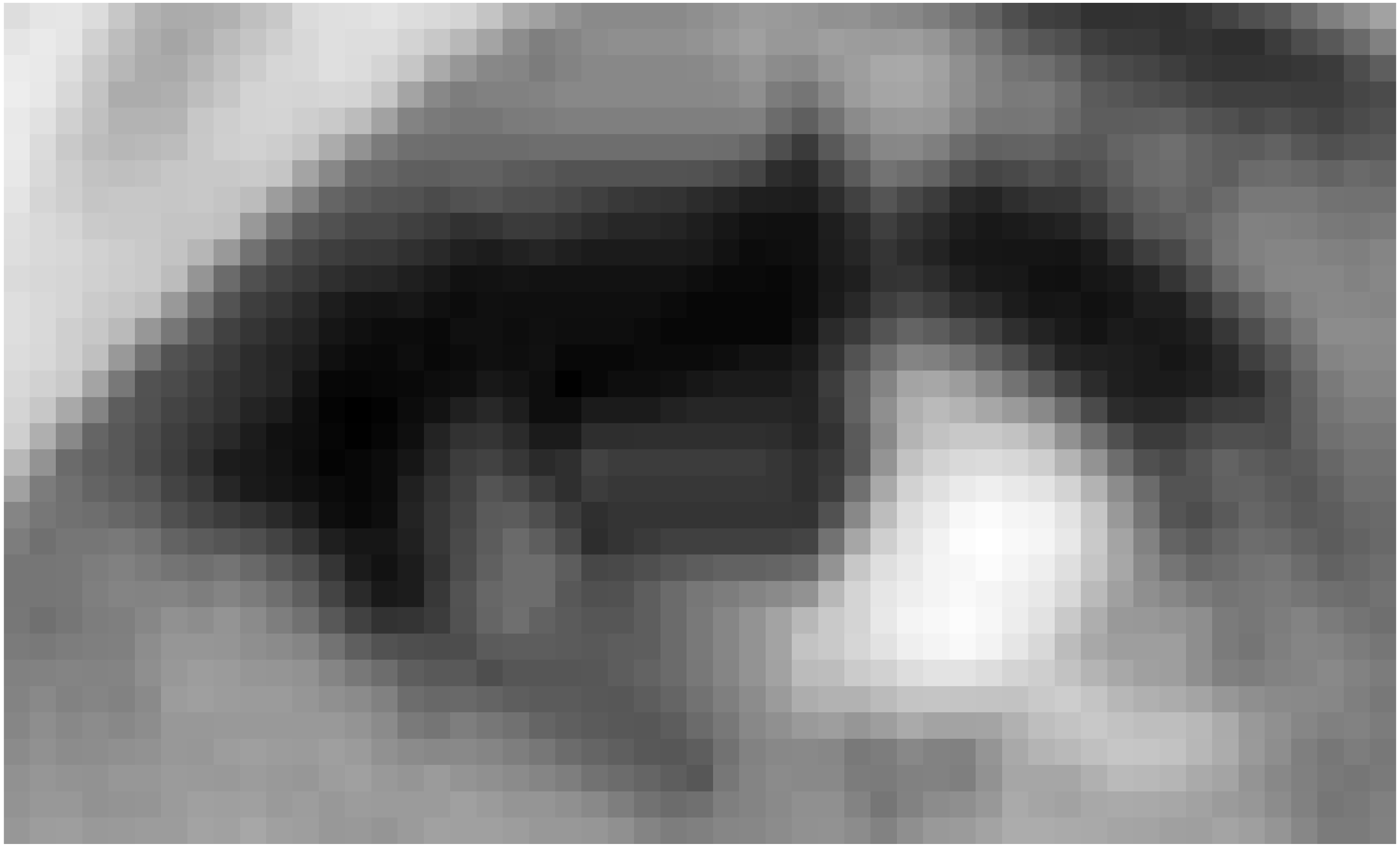}&
 	\includegraphics[width=0.33\textwidth]{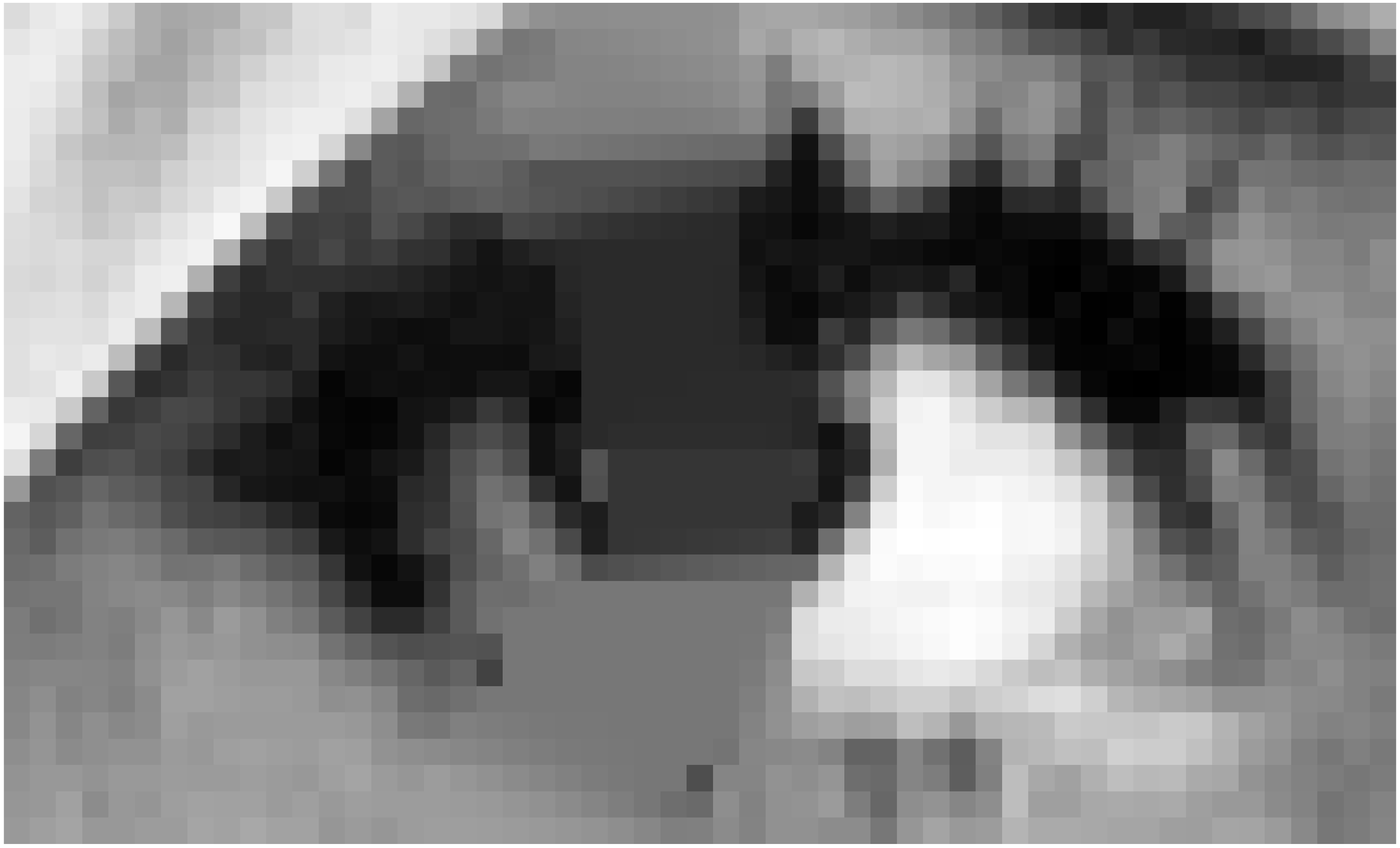}\\
 	(a) & (b) & (c)
 \end{array}$
 }\caption{\label{Sec5.NumEx.Fig:Imp.Det} Example \ref{Sec5.ImInp}.
 	Comparison of a detail of the original image with the corresponding detail of the 
 	restored images according to the local Moreau based compensated convex transforms 
 	and the TV-based method:  $(a)$ Right eye detail of the original image without overprinted scratches.
 	$(b)$ Right eye detail of the restored image $A_{\lambda}^M(f_K)$.
 	$(c)$ Right eye detail of the TV-based restored image.}
 \end{figure}

 \end{ex}


\begin{ex}\label{Sec5.SP}\textbf{Salt \& pepper noise removal.}
As a further example of scattered data approximation and application of Theorem \ref{Thm.3.ImpNew},
we consider the restoration of an image corrupted by
salt \& pepper noise by computing $A_{\lambda,\overline\Omega}^M(f_K)$ defined by \eqref{Sec5.NumEx.Eq:Imp.Rest}. 
Given the picture displayed in Figure \ref{Fig.Sec5.SP.0_7}$(a)$ with size $512\times 512$ pixels and 
damaged by $70\%$ salt \& pepper noise as shown in Figure \ref{Fig.Sec5.SP.0_7}$(b)$, let us denote by $K$ the set of the true pixels,
by $\Omega$ the domain of the whole image and $f_K$ the values of free noise pixels.
Since for $x\in\Omega$ the values of $C_{\lambda,\Omega}^l(f_{\overline\Omega,K}^M)(x)$ and 
$C_{\lambda,\Omega}^u(f_{\overline\Omega,K}^{-M})(x)$ depend on $f_K$ if $\dist(x,\partial\Omega)>4O_f/\lambda$, to 
reduce the boundary error due to the redefinition of $f_K$ on $\partial\Omega$, we consider an enlarged image and 
then restrict the restored image to the original domain. The enlarged image is obtained by padding one pixel before the first 
image element and after the last image element along each dimension, defining $f_K$ thereing equal to $\min_K f$
or $\max_K f$ according to whether we are computing $C_{\lambda,\Omega}^l(f_{\overline\Omega,K}^M)$ or 
$C_{\lambda,\Omega}^u(f_{\overline\Omega,K}^{-M})$, respectively. We apply then Algorithm \ref{Algo:MoreauEnv}
and Algorithm \ref{Algo:CnvxEnv} using a value of $tol=10^{-7}$ to compute the Moreau and convex based upper and 
lower compensated convex transforms that enter \eqref{Sec5.NumEx.Eq:Imp}, respectively. 
For the application of Algorithm \ref{Algo:MoreauEnv} the convergence check was on the $\ell^{\infty}$ norm of the 
error between two successive iterates, whereas in the application of Algorithm \ref{Algo:CnvxEnv}, to limit the computing 
time of the convex envelope scheme,  
we opted to make the convergence check by comparing the value of the $PSNR$ between two successive iterates.  
For $\lambda=15$ and $M=10^{13}$, the restored images with the Moreau and convex based scheme are displayed 
in Figure \ref{Fig.Sec5.SP.0_7}$(c)$ and Figure \ref{Fig.Sec5.SP.0_7}$(d)$, respectively. For the Moreau based
scheme we got a value of $\mathrm{PSNR}$ equal to $28.917\,\mathrm{dB}$ after $m=30$ iterations 
for a tolerance on the $\ell{\infty}$ norm of the difference between two succesive iterates equal to $10^{-7}$,
whereas for the convex based scheme, we got  $\mathrm{PSNR}=29.308\,\mathrm{dB}$ after $m=185$ iterations, 
giving an $\ell^{\infty}$ norm of the error between two succesive iterations equal to
$2.18\cdot 10^{-5}$. We note that the restored image by the convex based scheme presents a
slightly higher value of the $\mathrm{PSNR}$ than the one that uses the Moreau based 
transforms but at the expenses of an higher number of iterations. For both the restored images, there is no visible
effect on the boundary of the image. 
Finally, we consider the case of noise density equal to 99\%. Figure \ref{Fig.Sec5.SP.0_99} displays the
restored images by the Moreau and convex based definitions of the compensated transforms for $\lambda=5$ and $M=10^{13}$
obtained after a number of iterations equal to $m=78$ and $m=3081$, respectively.
The two restored images present comparable values of $\mathrm{PSNR}$, with the Moreau based scheme 
obtained after a much lower number of iterations than the one that uses the convex based definition of the 
lower and upper compensated convex transform. Even in this case there is no visible boundary effects for 
the two restored images, and as a result of the theoretical findings of \cite[Corollary 4.7]{ZCO18}, 
we note $A_{\lambda,\overline\Omega}^M(f_K)$ to be realized by an almost continuous piecewise affine interpolation.

\begin{figure}[htbp]
	\centerline{$\begin{array}{cc}
		\includegraphics[height=5cm]{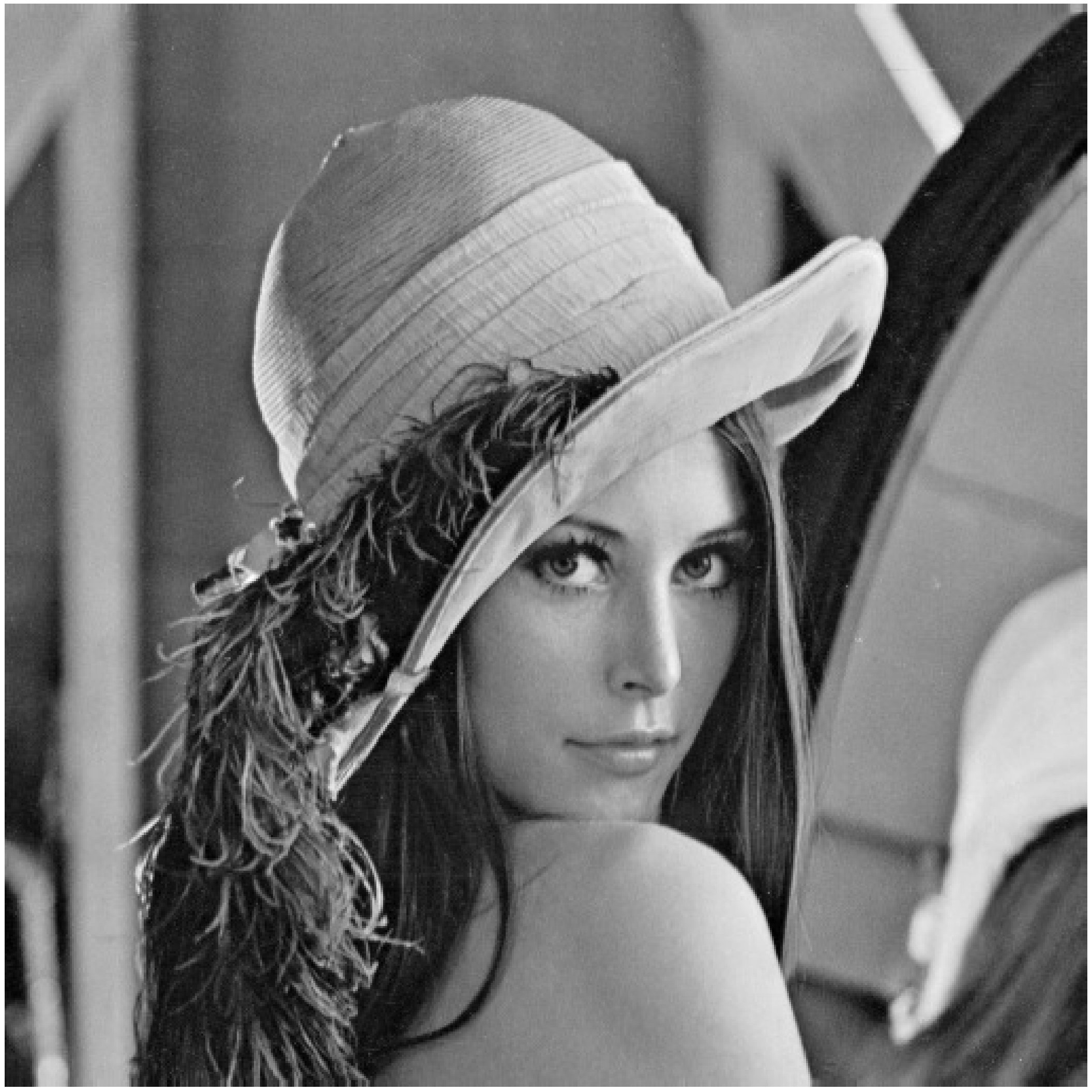}&   
		\includegraphics[height=5cm]{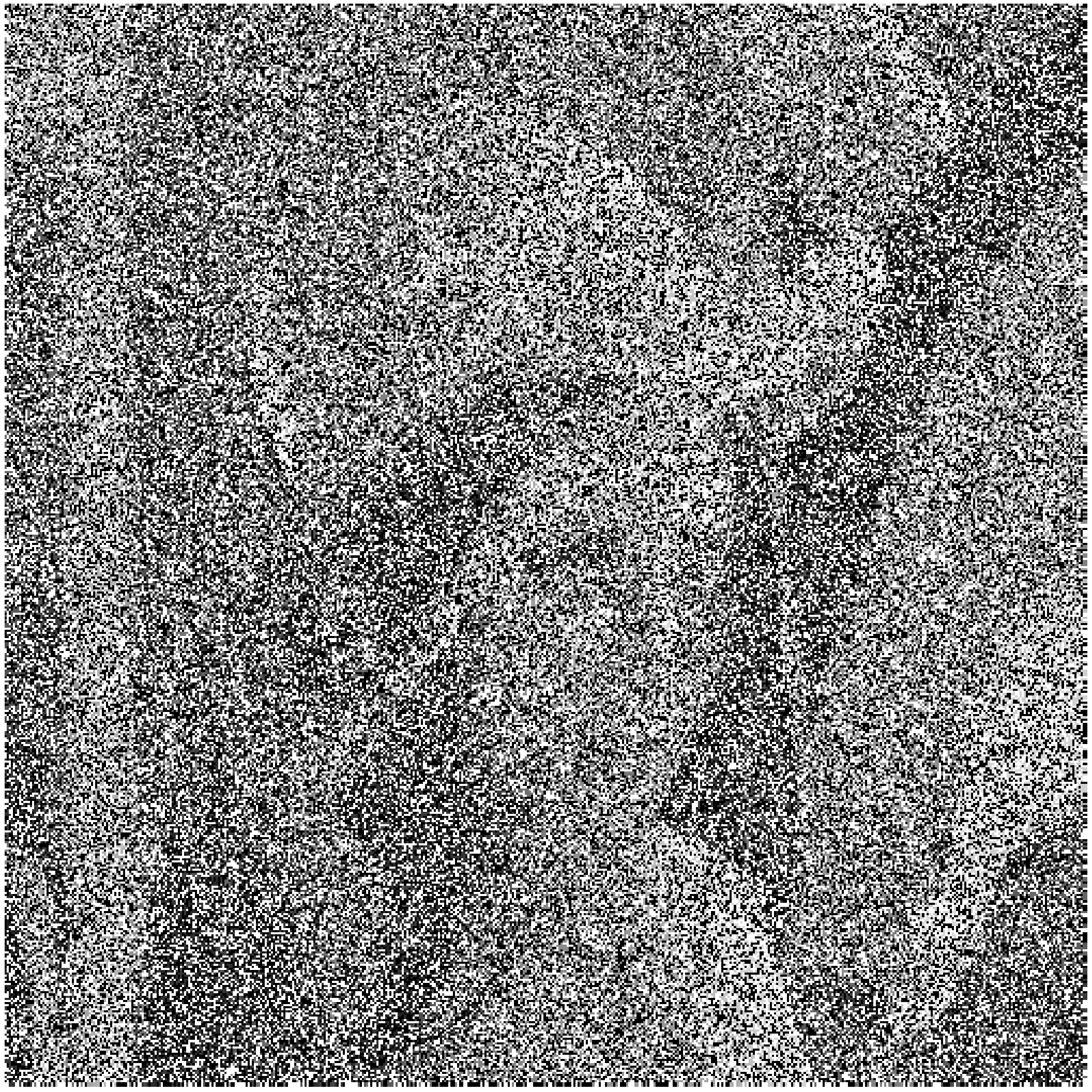}\\   
		(a)&(b)\\
		\includegraphics[height=5cm]{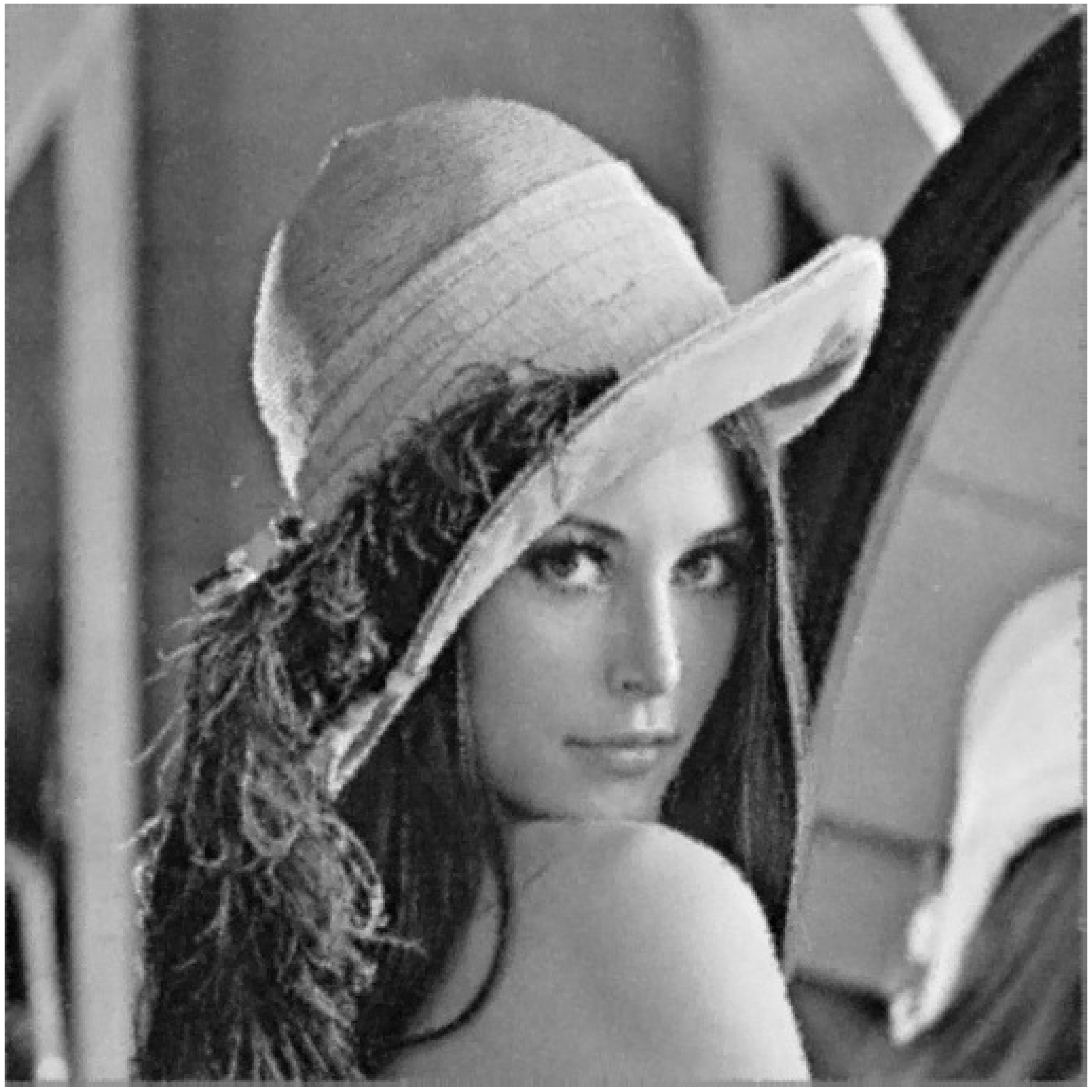}&	
		\includegraphics[height=5cm]{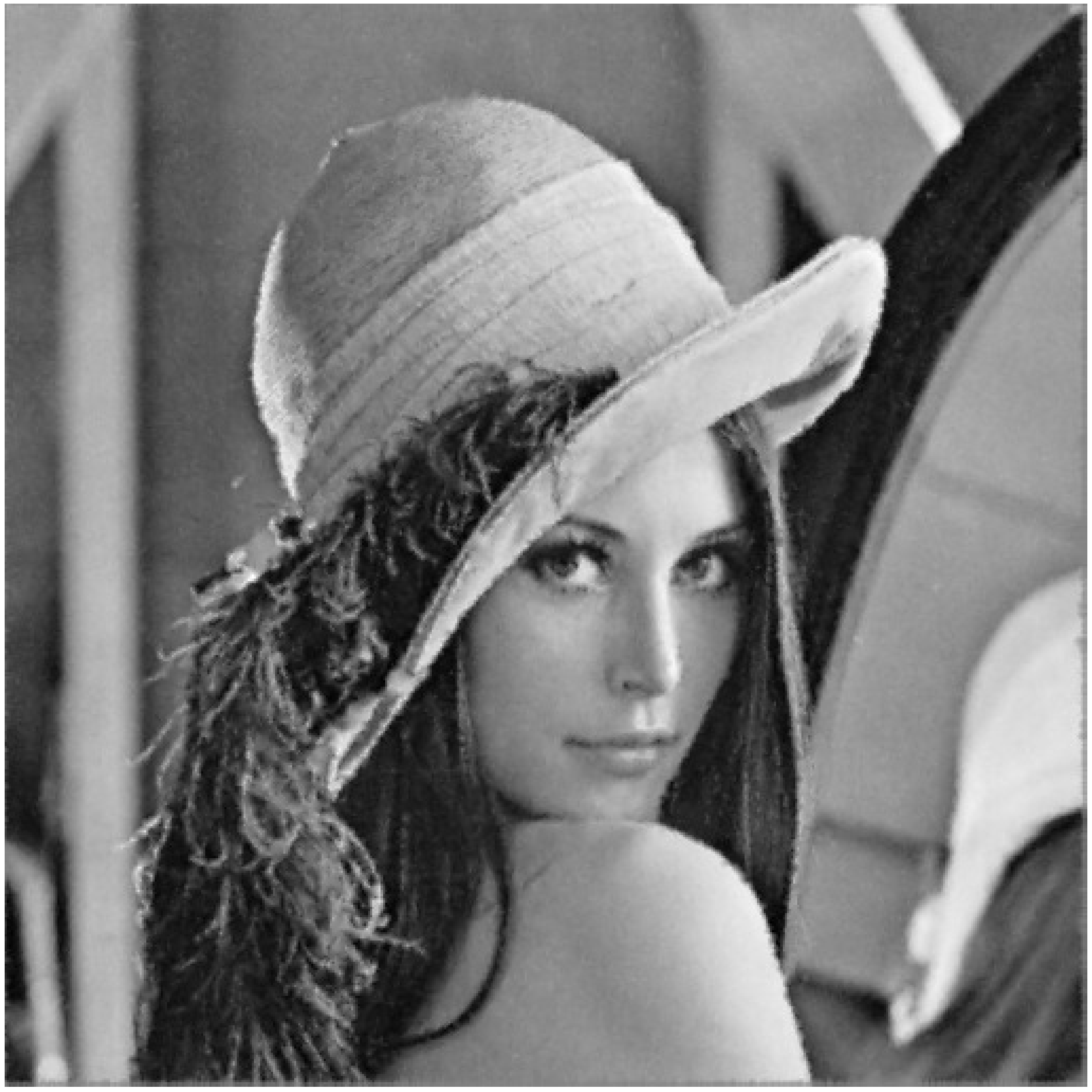}\\	
		 (c)&(d)
	\end{array}$}
	\caption{\label{Fig.Sec5.SP.0_7} Example \ref{Sec5.SP}.
		Restoration of $70\%$ corrupted image.
		$(a)$ Original image with size $512\times 512$;
		$(b)$ Original image covered by a salt \& pepper noise density of $70\%$. $\mathrm{PSNR}=6.998\,\mathrm{dB}$; 
		$(c)$ Restored image $A_{\lambda,\Omega}^M(f_K)$ using the Moreau based transforms computed by applying Algorithm \ref{Algo:MoreauEnv},
			with $\lambda=15$, $M=10^{13}$ and $tol=10^{-7}$. $\mathrm{PSNR}=28.917\,\mathrm{dB}$. Number of iterations $m=30$.
		$(d)$ Restored image $A_{\lambda,\Omega}^M(f_K)$ using the convex based transforms computed by applying Algorithm \ref{Algo:CnvxEnv}
			with $\lambda=15$, $M=10^{13}$ and $tol=10^{-7}$ on the error between the $\mathrm{PSNR}$ of two successive iterates.
			$\mathrm{PSNR}=29.304\,\mathrm{dB}$. 
			Number of iterations $m=185$ for an error between two succesive iterates of the convex based scheme 
			equal to $2.18\cdot 10^{-5}$.
	}
\end{figure}

\begin{figure}[htbp]
	\centerline{$\begin{array}{cc}
		\includegraphics[height=5cm]{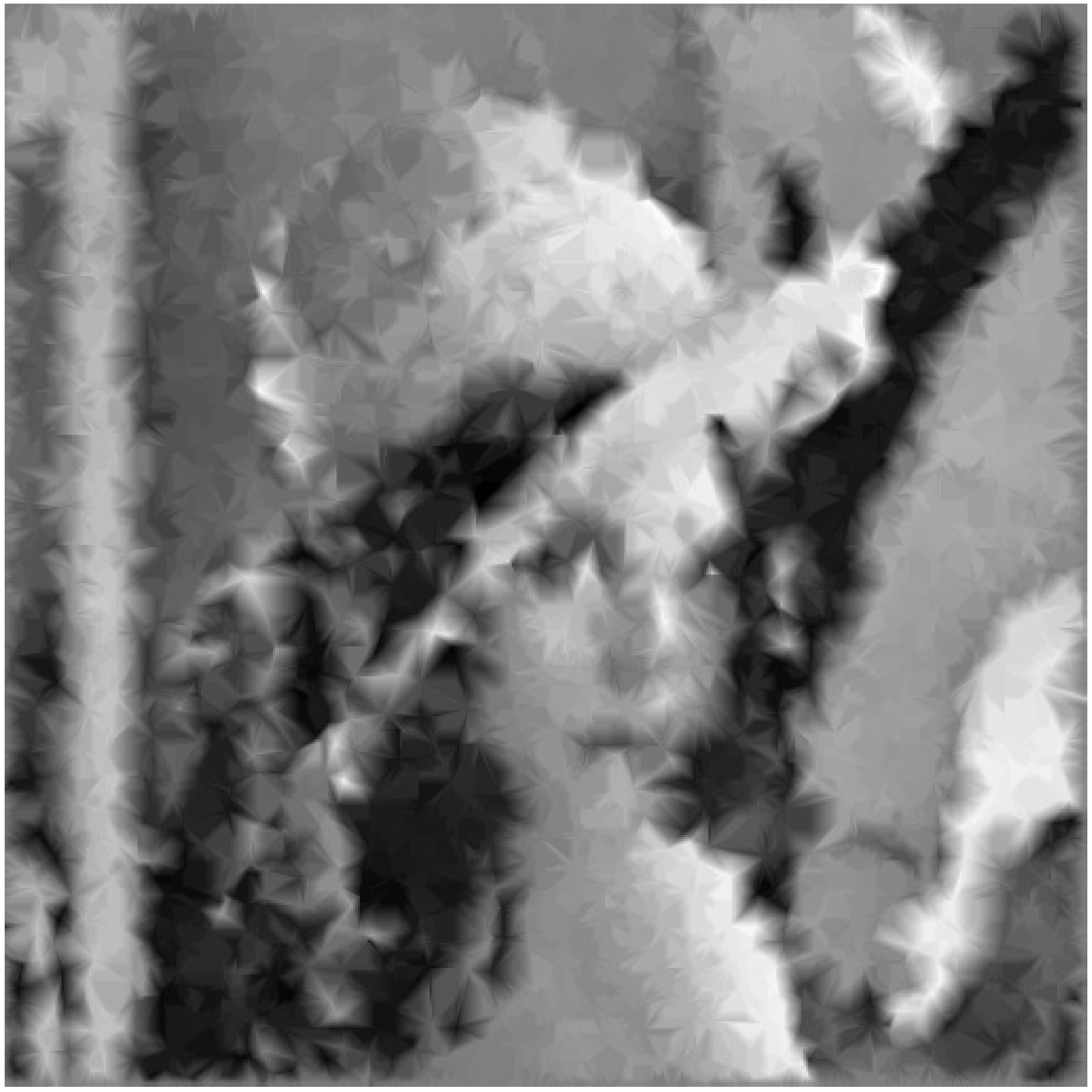}&	
		\includegraphics[height=5cm]{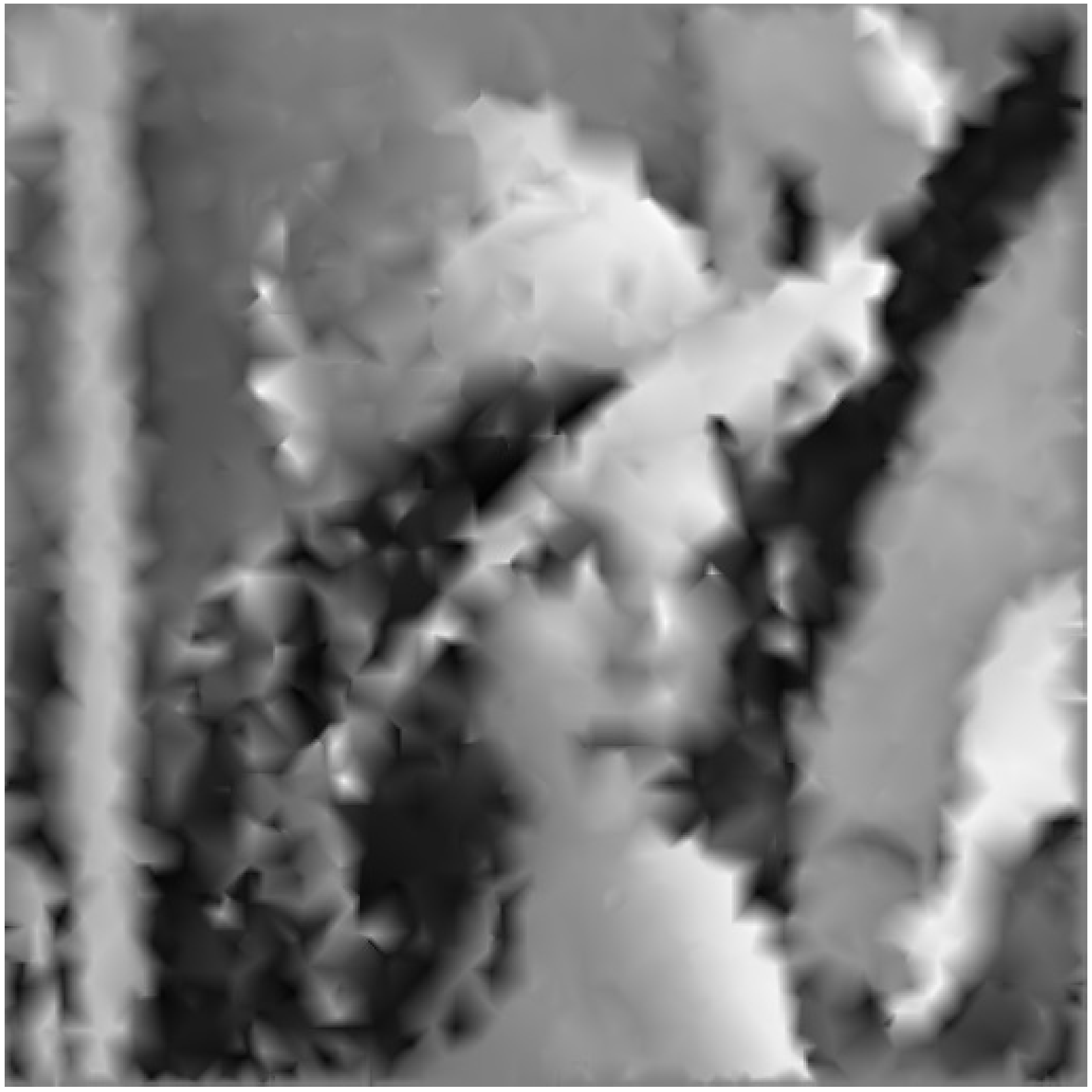}\\	
		(a)&(b)
	\end{array}$}
	\caption{\label{Fig.Sec5.SP.0_99} Example \ref{Sec5.SP}.
		  Restoration of $99\%$ corrupted image ($\mathrm{PSNR}=5.901\,\mathrm{dB}$):	
			$(a)$ Restored image $A_{\lambda}^M(f_K)$ by applying Algorithm  \ref{Algo:MoreauEnv} with 
				$\lambda=5$, $M=10^{13}$, $tol=10^{-7}$. 
				$\mathrm{PSNR}=28.918\,\mathrm{dB}$.  Number of iterations $m=78$.
			$(b)$ Restored Image $A_{\lambda}^M(f_K)$ by applying Algorithm  \ref{Algo:CnvxEnv} with 
				$\lambda=5$, $M=10^{13}$ and $tol=10^{-7}$.  
				$\mathrm{PSNR}=29.304\,\mathrm{dB}$.
		 		Number of iterations $m=3081$ for an error between two succesive iterates of the convex based scheme
				equal to $1.42\cdot 10^{-2}$.	
	}
\end{figure}

\end{ex}

 
\setcounter{equation}{0}
\section{Conclusions}\label{Sec.Conc}


Compensated convex transforms, or also known as proximity hull in the case of the lower transform, 
or grey scale erosion and dilation with quadratic structuring elements in mathematical morphology,
provide a geometric tight--approximation method for general functions that yields novel ways to 
smooth functions, to identify
singularities in functions, and to interpolate and approximate data. Introduced for real valued functions 
defined in $\mathbb{R}^n$, in view of applications to image processing and computer aided geometric design, 
in this paper we have proposed a definition for functions defined on bounded convex open sets such that the 
following equality holds

\begin{equation*}
	C^l_{\lambda,\,\Omega}(f_{\overline\Omega}^{-})(x)=C^l_{\lambda}(f_{\R^n}^-)(x)\quad\text{and}\quad
	C^u_{\lambda,\,\Omega}(f_{\overline\Omega}^{+})(x)=C^l_{\lambda}(f_{\R^n}^+)(x)\quad\text{for any }
	x\in\overline\Omega\,,
\end{equation*}
where $f_{\R^n}^-$ and $f_{\R^n}^+$ are suitable extensions of $f$ to $\mathbb{R}^n$,
$C^l_{\lambda}(f_{\R^n}^-)$ and $C^u_{\lambda}(f_{\R^n}^+)$ are the `global' compensated 
convex transforms \eqref{Eq.01.Def.LwTr} and \eqref{Eq.01.Def.UpTr}, 
whereas $C^l_{\lambda,\,\Omega}(f_{\overline\Omega}^{-})$ and $C^l_{\lambda,\,\Omega}(f_{\overline\Omega}^{+})$ 
are the local versions defined in terms of the Moreau envelopes as \eqref{Def.1.LLCT} and \eqref{Def.1.LUCT}, respectively.
We have thus proposed a new algorithm to compute the Moreau envelope which 
has linear complexity with respect to the number of grid points and we have given applications to
tasks of image processing such as image inpainting and restoration of image with 
high density of salt \& pepper noise, and of shape interrogation such as detection of intersection of 
sampled geometries and multiscale medial axis map. The performance of the methods and the accuracy 
of the results show that, when coupled with an efficient numerical scheme, the theory of 
compensated convex transforms provides a valid and feasible alternative to state-of-art methods especially
for processing data without any a prior information or represented by point cloud.


\section*{Acknowledgements}
The authors are grateful to the anonymous reviewers for their comments and suggestions.
KZ wishes to thank The University of Nottingham for its support, 
AO acknowledges the partial financial support of the Argentinian Research Council (CONICET)
through the project PIP 100,
the National University of Tucum\'{a}n through the project PIUNT CX-E625 and the FonCyT 
through the project PICT 2016 201-0105 Prestamo BID, whereas
EC is grateful for the financial support of the College of Science, Swansea University.


\appendix


\setcounter{equation}{0}
\section{Appendix: Proofs of main results}\label{Sec.Proofs}


\noindent {\bf Proof of \eqref{Eq.01.MixMor}}.
Next we prove only the characterization of the lower compensated convex transform.
The one for the upper transform is obtained by using the relationship \eqref{Eq.RelMor}.
From the definition \eqref{Eq.01.Def.LwTr} of $C_{\lambda}^(f)$ and the characterization \eqref{Eq.Prp.GlCnvx.2}
of the convex envelope, we have
\begin{align}
	C_{\lambda}^l(f)(x)&=\co[f+\lambda|\cdot|^2](x)-\lambda|x|^2\notag	\\[1.5ex]
		   &=\underset{\ell\in Aff(\mathbb{R}^n)}{\sup}\{\ell(x):\,\ell(y)\leq f(y)+\lambda|y|^2\}-\lambda|x|^2	\notag\\[1.5ex]
		   &=\underset{\ell\in Aff(\mathbb{R}^n)}{\sup}\{\ell(x)-\lambda|x|^2:\,\ell(y)-\lambda|y|^2\leq f(y)\}\,,
				   \label{Eq.Apndx.Charact.01}
\end{align}
which shows the lower compensated transform to be equal to the envelope of quadratic functions with given
curvature $\lambda$. If we represent the elements of $\Aff(\R^n)$ as
\[
	\ell(x)=a\cdot x+b,\quad a\in\R^n,\,\,b\in\R\,,
\]
then 
\[
	\begin{split}
	\ell(x)-\lambda|x|^2	&=a\cdot x+b-\lambda|x|^2\\[1.5ex]
				&=-\lambda\left|x-\frac{a}{2\lambda}\right|^2+\frac{|a|^2}{4\lambda}+b\,.
	\end{split}
\]
thus \eqref{Eq.Apndx.Charact.01} reads also as
\[
	\begin{split}
		C_{\lambda}^l(f)(x)&=
		\underset{\substack{z\in\R^n\\c\in\R}}{\sup}
			\left\{-\lambda|x-z|^2+c:\,-\lambda|y-z|^2+c\leq f(y)\right\}\\[1.5ex]
		&=\underset{z\in\R^n}{\sup}\left\{-\lambda|x-z|^2+
		\underset{y\in\R^n}{\inf}\left\{f(y)+\lambda|y-z|^2\right\}
		\right\}\\[1.5ex]
		&=M^{\lambda}(M_{\lambda}(f))(x)\,,
	\end{split}
\]
where we have used the fact that for any $z\in\R^n$, the largest $c\in \R$ such that
\[
	c\leq f(y)+\lambda|y-z|^2
\]
is given by
\[
	\underset{z\in\mathbb{R}^n}{\inf}\left\{f(y)+\lambda|y-z|^2\right\}\,,
\]
which is $M_{\lambda}(f)(z)$. \hfill\qed\\


\noindent {\bf Proof of Proposition \ref{Prop.3.Mor}.} \textit{Part $(i)$}:  Since  
\begin{equation}
	M_\lambda(f)(x)=\underline{f}(y_x)+\lambda|y_x-x|^2=
	\inf\{\underline{f}(y)+\lambda|y-x|^2,\; y\in\mathbb{R}^n\}\leq \underline{f}(x)
\end{equation}
then
\begin{equation}
	\lambda|y_x-x|^2\leq \underline{f}(x)-\underline{f}(y_x)\leq M-m=O_f\,,
\end{equation}
which concludes the proof for $M_\lambda(f)(x)$. The proof for $M^\lambda(f)(x)$ follows 
similar arguments.
\hfill\qed\\


\textit{Part $(ii)$}: Again we have
\begin{equation}
	M_\lambda(f)(x)={f}(y_x)+\lambda|y_x-x|^2=
	\inf\{{f}(y)+\lambda|y-x|^2,\; y\in\mathbb{R}^n\}\leq {f}(x)
\end{equation}
 so that $\lambda|y_x-x|^2\leq f(x)-f(y_x)\leq L|y_x-x|$. Thus
 $|y_x-x|\leq L/\lambda$. The proof for $M^\lambda(f)(x)$ follows similar arguments.
\hfill\qed\\


\noindent {\bf Proof of Proposition \ref{Sec2:Prop:CmprCmpTr}:} 
We first observe that by definition, for $z\in\overline\Omega$
\[
\begin{split}	
	M_\lambda(f^\infty)(z)	&=\inf\{f^\infty(y)+\lambda|y-z|^2: \,y\in\mathbb{R}^n\}\\[1.5ex]
				&=\inf\{\underline{f}(y)+\lambda|y-z|^2: \,y\in\overline\Omega\}\\[1.5ex]
				&=M_{\lambda,\Omega}(\underline{f})(z)\\[1.5ex]
				&=M_{\lambda,\Omega}(f)(z)\,,
\end{split}	
\]
as $f^\infty(y)+\lambda|y-z|^2=+\infty$ if $y\in\mathbb{R}^n\setminus{\overline\Omega}$, thus for 
$x\in\overline\Omega$, we have
\begin{equation}\label{Eq:Proof:Compar1}
\begin{split}	
	C^l_\lambda(f^\infty)(x)&=M^\lambda(M_\lambda(f^\infty))(x)\\[1.5ex]
				&=\sup\{M_\lambda(f^\infty)(z)-\lambda|x-z|^2:\; z\in\mathbb{R}^n\}\\[1.5ex]
				&\geq \sup\{M_\lambda(f^\infty)(z)-\lambda|x-z|^2:\; z\in\overline\Omega\}\\[1.5ex]
				&=\sup\{M_{\lambda,\Omega}(f)(z)-\lambda|x-z|^2:\; z\in\overline\Omega\}\\[1.5ex]
				&=C^l_{\lambda,\Omega}(f)(x)\,,
\end{split}	
\end{equation}
hence $C^l_\lambda(f^\infty)(x)\geq C^l_{\lambda,\Omega}(f)(x)$.
By the definition of lower transform for $x\in \overline\Omega$, we find that
\begin{equation}\label{Eq:Proof:Compar2}
\begin{split}
	C^l_\lambda(f^\infty)(x)&=\co[f^\infty(\cdot)+\lambda|\cdot-x|^2](x)\\[1.5ex]
				&=\inf\Big\{\sum^{n+1}_{j=1}\lambda_j\bigg(\underline{f}(x_j) +\lambda|x_j-x|^2\bigg):\,
					\lambda_j\geq 0,\, x_j\in\overline\Omega,\; j=1,2,\dots,n+1,\\[1.5ex]
				&\phantom{xxxxxxxxxxxxxxxxxxxxxxxxxxxx}\sum_{j=1}^{n+1}\lambda_j=1,\, 
					\sum_{j=1}^{n+1}\lambda_jx_j=x\Big\}\\[1.5ex]
				&\leq \underline{f}(x)\,,
\end{split}
\end{equation}
given that if there exists $x_i\notin \overline\Omega$ and $\lambda_j>0$, then the sum is $+\infty$.
Furthermore, since $M_{\lambda,\Omega}(f_{\overline\Omega}^{-})=M_{\lambda,\Omega}(\underline{f_{\overline\Omega}^{-}})$
and by definition \eqref{Eq.01.Ext.Local.LW}, \eqref{Sec2.Eq.SemClos} and \eqref{Sec2.Eq.SemClosBnd},
for $x\in\overline\Omega$  it is $\underline{f^{-}_{\overline\Omega}}(x)\leq \underline{f}(x)$, 
then it is $M_{\lambda,\Omega}(f^{-}_{\overline\Omega})(x)\leq M_{\lambda,\Omega}(\underline{f})(x)=M_{\lambda,\Omega}(f)(x)$, hence
\begin{equation}\label{Eq:Proof:Compar3}
	\begin{split}
		C^l_{\lambda,\Omega}(f^-_{\overline\Omega})(x)	&= M^{\lambda}_{\Omega}(M_{\lambda,\Omega}(f^{-}_{\overline\Omega}))(x)\\[1.5ex]
								&\leq M^{\lambda}_{\Omega}(M_{\lambda,\Omega}(f))(x)\\[1.5ex]
								&=C^l_{\lambda,\Omega}(f)(x)\,,
	\end{split}
\end{equation}
thus $C^l_{\lambda,\Omega}(f^-_{\overline\Omega})(x)\leq C^l_{\lambda,\Omega}(f)(x)$. 
By comparing \eqref{Eq:Proof:Compar1}, \eqref{Eq:Proof:Compar2} and \eqref{Eq:Proof:Compar3}
we conclude the proof.
\hfill\qed\\

\noindent {\bf Proof of Theorem \ref{Theo.LocLwTr}}. We next present the proofs only for 
the lower transform, given that the results for the upper transform will follow from \eqref{Eq.RelMor}. 

{\bf Proof of \eqref{Sec3.Eq.01}:}
Without loss of generality, suppose that $\inf_{\Omega}f=0$. 
By definition, for $x\in\overline \Omega$, we have
\begin{equation}\label{Eq:Apndx.Ineq}
\begin{split}
	M_{\lambda,\Omega}(f_{\overline{\Omega}}^{-})(x)	&=\inf\{f_{\Omega}^{-}(z)+\lambda|z-x|^2,\; z\in \overline \Omega\}\\[1.5ex]
				&\geq \inf\{f^-_{\R^n}(y)+\lambda|y-x|^2,\; y\in \mathbb{R}^n\}\\[1.5ex]
				&=M_\lambda(f^-_{\R^n})(x)\,,
\end{split}
\end{equation}
given that the second infimum is taken over a larger set. 
In order to show that also the opposite inequality holds, first
let $\left(y_k\right)_{k\in\mathbb{N}}$ be a minimizing sequence of $g(y):= f^-_{\mathbb{R}^n}(y)+\lambda|y-x|^2$ 
in $\mathbb{R}^n$. Clearly $g$ is coercive in $\mathbb{R}^n$ in the sense that $g(y)/|y|\to+\infty$ 
as $|y|\to+\infty$. Thus  $\left(y_k\right)_{k\in\mathbb{N}}$ is bounded in $\mathbb{R}^n$ and, therefore,
there exists a $y_x\in\mathbb{R}^n$ and a convergent subsequence $\left(y_{k_j}\right)_{k_j\in\mathbb{N}}$ 
such that $y_{k_j}\to y_x\in\mathbb{R}^n$.

Next we show that $y_x\in\overline\Omega$ and then we will use $y_x$ to prove that there also holds 
$M_\lambda(f^-_{\R^n})(x)\geq M_{\lambda,\Omega}(f_{\overline{\Omega}}^{-})(x)$.

If $y_x\notin\overline\Omega$, since $f^-_{\mathbb{R}^n}(y)=0=\inf_{\Omega} f$ for 
$y\in \mathbb{R}^n\setminus \Omega$, then  $f^-_{\mathbb{R}^n}$ is lower semicontinuous in 
$\mathbb{R}^n\setminus \Omega$. As a result, $g(y):= f^-_{\mathbb{R}^n}(y)+\lambda|y-x|^2$ is 
lower semicontinuous at $y_x$ and $y_x$ is a minimizer of $g(y)$ in $\mathbb{R}^n$. 
We would have therefore 
\[
	\begin{split}
		M_\lambda(f^-_{\mathbb{R}^n})(x)&=\inf\{f^-_{\mathbb{R}^n}(y)+\lambda|y-x|^2,\;y\in\mathbb{R}^n\}\\[1.5ex]
			&=f^-_{\mathbb{R}^n}(y_x)+\lambda|y_x-x|^2\\[1.5ex]
			&=\lambda|y_x-x|^2\,.
	\end{split}
\]
Since $y_x\neq x$, along the line segment $[x,\,y_x]$, there would then exist a point $\xi_x\in\partial\Omega$
such that 
\begin{equation}
	\begin{split}
		M_\lambda(f^-_{\mathbb{R}^n})(x)=\lambda|y_x-x|^2&>\lambda|\xi_x-x|^2\\[1.5ex]
			&\geq \inf\{f^-_{\overline\Omega}(z)+\lambda|z-x|^2,\;z\in\overline\Omega\}
			=M_{\lambda,\Omega}(f^-_{\overline\Omega})(x)\,,
	\end{split}
\end{equation}
which is a contadiction, and therefore $y_x\in\overline\Omega$.

Let us now examine separately the case of $y_x\in\partial\Omega$ and $y_x\in\Omega$.
If $y_x\in\partial\Omega$, by an argument similar to the above,
we conclude that $g(y)$ is lower semicontinuous at $y_x\notin\Omega$ and therefore
\begin{equation}\label{Eq:MorEnvFRn}
	M_\lambda(f^-_{\mathbb{R}^n})(x)=\lambda|y_x-x|^2\,.
\end{equation}
We now have to distinguish the two cases of $x\in\partial\Omega$ and $x\in\Omega$.
If $x\in\partial \Omega$, from \eqref{Eq:Apndx.Ineq}, \eqref{Eq:MorEnvFRn} 
and the ordering property of the lower Moreau envelope \eqref{Eq:OrderMoreau}, we have
\begin{equation}\label{Eq:Apndx.IneqEx}	
\begin{split}
	\lambda|y_x-x|^2&=M_\lambda(f^-_{\mathbb{R}^n})(x)\\[1.5ex]
			&\leq M_{\lambda,\Omega}(f^{-}_{\overline\Omega})(x)\\[1.5ex]
			&\leq f_{\overline\Omega}^-(x)\\[1.5ex]
			&=0\,,
\end{split}
\end{equation}
that is, $y_x=x$, hence
\begin{equation}\label{Eq:Apndx.Ineq1a}	
	M_{\lambda,\Omega}(f^{-}_{\overline\Omega})(x)\leq f_{\overline\Omega}^-(x)=0=
	\lambda|y_x-x|^2=M_\lambda(f^-_{\mathbb{R}^n})(x)\,.
\end{equation}
If $x\in \Omega$, on the other hand, let us consider the convergent minimizing sequence 
$(y_{k_j})_{j\in\mathbb{N}}\in\R^n$ of $f_{\R^n}^-(y)+\lambda|y-x|^2$.
We note that for $j\in\mathbb{N}$ such that $y_{k_j}\in\overline\Omega$, we have
\begin{equation}\label{Eq:1stcase}
	f^{-}_{\mathbb{R}^n}(y_{k_j})+\lambda|y_{k_j}-x|^2=
	f_{\overline\Omega}^{-}(y_{k_j})+\lambda|y_{k_j}-x|^2\geq M_{\lambda}(f_{\overline\Omega}^{-})(x)\,,
\end{equation}
and for $j\in\mathbb{N}$ such that $y_{k_j}\in\R^n\setminus\overline\Omega$, there exists 
$z_j\in[x,\,y_{k_j}]$ and $z_j\in\partial\Omega$ 
so that 
\begin{equation}\label{Eq:2ndcase}
\begin{split}
	f^{-}_{\mathbb{R}^n}(y_{k_j})+\lambda|y_{k_j}-x|^2&=\lambda|y_{k_j}-x|^2	\\[1.5ex]	
							&\geq \lambda|z_j-x|^2\\[1.5ex] 
							&=f_{\overline\Omega}^{-}(z_j)+\lambda|z_j-x|^2\\[1.5ex]
							&\geq M_{\lambda}(f_{\overline\Omega}^{-})(x)\,.
\end{split}
\end{equation}
Thus, from \eqref{Eq:1stcase} and \eqref{Eq:2ndcase} we get that for any $j\in\mathbb{N}$, it is
\[
	f^{-}_{\mathbb{R}^n}(y_{k_j})+\lambda|y_{k_j}-x|^2\geq 
	M_{\lambda}(f_{\overline\Omega}^{-})(x)\,.
\]
By letting $j\to\infty$, since $(y_{k_j})_{j\in\mathbb{N}}$ is a minimizing sequence of 
$f^{-}_{\mathbb{R}^n}(y)+\lambda|y-x|^2$, we conclude that 
\begin{equation}\label{Eq:IneqBnd}
	M_{\lambda}(f_{\mathbb{R}^n}^-)(x)\geq 
	M_{\lambda}(f_{\overline\Omega}^{-})(x)\,.
\end{equation}
If $y_x\in \Omega$, for $j>0$ sufficiently large, $y_{k_j}\in\Omega$, thus
\[
	\begin{split}
		f(y_{k_j})+\lambda|y_{k_j}-x|^2&\geq \inf\{f^-_{\overline\Omega}(z)+\lambda|z-x|^2,\;z\in\bar\Omega\}\\[1.5ex]
			&=M_{\lambda,\Omega}(f^-_{\overline\Omega})(x)\,.
	\end{split}
\]
Since $(y_{k_j})_{j\in\mathbb{N}}$ is a minimizing sequence of $g(y)$ in $\mathbb{R}^n$, 
if we let $j\to\infty$, we obtain also in this case 
\begin{equation}\label{Eq:Apndx.Ineq2}
	M_\lambda(f^-_{\mathbb{R}^n})(x)\geq M_{\lambda,\Omega}(f^-_{\overline\Omega})(x)\,.
\end{equation}	
By comparing \eqref{Eq:IneqBnd} and \eqref{Eq:Apndx.Ineq2} with \eqref{Eq:Apndx.Ineq}, we conclude that 
$M_\lambda(f^-_{\mathbb{R}^n})(x)= M_{\lambda,\Omega}(f^-_{\overline\Omega})(x)$ for $x\in\overline\Omega$
and this completes the proof.
\hfill \qed\\

{\bf Proof of \eqref{Sec3.Eq.02}:} We first observe that for
$x\notin \overline \Omega$, as $0$ is the minimum value of $f_{\R^n}^-$, we have 
$M_\lambda(f_{\R^n}^-)(x)=0$. Thus 
\begin{equation}
	M_\lambda(f^-_{\R^n})(x)=\left\{
		\begin{array}{ll} 
	\displaystyle M_{\lambda,\Omega}(f_{\overline{\Omega}}^{-})(x),  & \displaystyle x\in\overline \Omega,\\[1.5ex]
	\displaystyle 0,		      & \displaystyle x\notin \overline \Omega.
		\end{array}\right.
\end{equation}
To complete the proof, we need to show that for $x\in \overline \Omega$,
$M^\lambda(M_\lambda(f_{\R^n}^-))(x)=M^{\lambda}_{\Omega}(M_{\lambda,\Omega}(f_{\overline\Omega}^{-}(x))$. 
By definition, for  $x\in\overline \Omega$, we have that
\begin{equation}
	\begin{split}
		&M^\lambda(M_\lambda(f^-_{\R^n}))(x)=\sup\Bigg\{M_\lambda(f^-_{\R^n})(y)-\lambda|y-x|^2,\; y\in \mathbb{R}^n\Bigg\}\\[1.5ex]
		&=\max\Bigg\{\sup\Big\{M_{\lambda}(f_{\overline\Omega}^{-})(y)-\lambda|y-x|^2,\; 
			y\in\overline\Omega\Big\},\; 
		\sup\Big\{M_\lambda(f^-_{\R^n})(y)-\lambda|y-x|^2,\; y\notin \overline \Omega\Big\}
			\Bigg\}\\[1.5ex]
		&\geq M^{\lambda}_{\Omega}(M_{\lambda,\Omega}(f_{\overline\Omega}^{-}(x)).
	\end{split}
\end{equation}
Since $f_{\mathbb{R}^n}^-$ is bounded in $\R^n$, then $M_{\lambda}(f_{\mathbb{R}^n}^-)$ is also a bounded function
in $\R^n$, thus $g(y):=M_{\lambda}(f_{\mathbb{R}^n}^-)(y)-\lambda|y-x|^2$ is coercive, in the sense that
$\lim_{|y|\to\infty}\,g(y)/|y|=-\infty$. Furthermore, $g=g(y)$ is also continuous because 
so is $M_{\lambda}(f_{\mathbb{R}^n}^-)(y)$ \cite[Proposition 1.1]{AA93}, thus a
maximum point, say, $y_x\in\mathbb{R}^n$ of $g=g(y)$ exists. If $y_x\notin\overline \Omega$, we have
\begin{equation}
	M^\lambda(M_\lambda(f_{\R^n}^-))(x)=-\lambda|y_x-x|^2
\end{equation}
as $M_\lambda(f^-_{\R^n})(y_x)=0$. Thus $y_x\neq x$. Again if we take the point $\eta_x\in \partial \Omega$ 
in the line segment $[x,\, y_x]$, then
\begin{equation}
	-\lambda|y_x-x|^2<-\lambda|\eta_x-x|^2\leq 
	\sup\{M_{\lambda,\Omega}(f_{\overline\Omega}^{-})(z)-\lambda|z-x|^2,\; 
	z\in\overline \Omega\}=M^{\lambda}_{\Omega}(M_{\lambda,\Omega}(f_{\overline\Omega}^{-}(x)).
\end{equation}
This is a contradiction. Thus $
M^{\lambda}(M_\lambda(f^-_{\R^n}))(x)=M^{\lambda}_{\Omega}(M_{\lambda,\Omega}(f_{\overline\Omega}^{-}(x))$
for $x\in\overline \Omega$ so that
\begin{equation}
	C^l_{\lambda,\Omega}(f_{\overline\Omega}^{-})(x)=C^l_{\lambda}(f_{\R^n}^-)(x) 
\end{equation}
for $x\in\overline \Omega$ and this concludes the proof.
\hfill \qed\\

As for the estimate in $(i)$, it follows as argued above that there exists a $z_x\in\mathbb{R}^n$ such that 
\[
	M_{\lambda,\Omega}(f_{\overline{\Omega}}^{-})(x)=f_{\overline{\Omega}}^{-}(z_x)+\lambda|z_x-x|^2
\]
which yields the estimate
\[
	|z_x-x|^2\leq \frac{O_f}{\lambda}\,.
\]
Thus if $\dist(x,\partial\Omega)>\sqrt{O_f/\lambda}$, then by the triangle inequality we conclude that $z_x\in\Omega$.

The estimate in $(ii)$ is obtained by a similar argument applied to $M_{\lambda,\Omega}(f_{\Omega}^{-})$.
\hfill \qed\\ 

%
\noindent {\bf Proof of Corollary \ref{Cor.SetBnd}:}
	Since $f(x)=0$ for $x \not\in \overline{\Omega}$, this follows from the fact that $f=f^-_{\R^n}$, 
	the definition of $f_{\overline\Omega}^{-}$ and \eqref{Sec3.Eq.03}.
\hfill \qed\\
%
%

\noindent {\bf Proof of Corollary \ref{Cor.SetCmp}:}
	Same arguments as in the previous corollary apply.
\hfill \qed\\

%

\noindent {\bf Proof of Theorem \ref{Thm.3.ChrctFnct}:}
First we show that if $x\in\R^n\setminus\overline\Omega$, then $M^{\lambda}(\chi_K)(x)=0$. 
Suppose that for some $x\in\R^n\setminus\overline\Omega$,
$M^{\lambda}(\chi_K)(x)\not=0$. By definition \eqref{Eq.01.Def.UpLwMor} of $M^{\lambda}(\chi_K)$ and since $x\not\in K$, we have that 
\[
	M^{\lambda}(\chi_K)(x)\geq \chi_K(x)-\lambda|x-x|^2=0\,.
\]	
Since $\chi_K$ is upper-semicontinuous and $\chi_K(y)-\lambda|y-x|^2\to-\infty$ as $|y|\to\infty$, 
there exists
a $y_0\in\mathbb{R}^n$ such that 
\[
	M^{\lambda}(\chi_K)(x)=\chi_{K}(y_0)-\lambda|y_0-x|^2>0\,, 
\]
thus $y_0\in K$, otherwise $M^{\lambda}(\chi_K)(x)\leq 0$.
So we have that 
\[
	1-\lambda|y_0-x|^2>0\,,
\]
thus 
\[
	|y_0-x|\leq \tfrac{1}{\sqrt{\lambda}}\,,
\]
which contradicts that $\dist(K,\partial\Omega)>\tfrac{1}{\sqrt{\lambda}}$.

Next we show that for any $x\in\overline\Omega$, $M_{\Omega}^{\lambda}(\chi_{K}^{\Omega})(x)=M^{\lambda}(\chi_K)(x)$. 
If, for some $x\in\overline\Omega$, there exists $y_0\not\in\overline\Omega$ such that 
\[
	M^{\lambda}(\chi_K)(x)=\chi_K(y_0)-\lambda|y_0-x|^2
\] 
then 
\[
	\chi_K(y_0)=0\,, 
\]
so that 
\[
	M^{\lambda}(\chi_K)(x)=-\lambda|y_0-x|^2<0\,. 
\]
which contradicts the fact that $M^{\lambda}(\chi_K) (x)\geq 0$ for any $x\in\mathbb{R}^n$. 
Thus $y_0\in\overline\Omega$ and so $M^{\lambda}(\chi_K)(x)=M_{\Omega}^{\lambda}(\chi_{K}^{\Omega})(x)$ for $x\in\overline\Omega$.

Now we consider the mixed transform. Let $x\in\overline\Omega$. If there exists $y_0\not\in\overline\Omega$ such that 
\[
	M_{\lambda}(M^{\lambda}(\chi_K))(x)=M^{\lambda}(\chi_K)(y_0)+\lambda|y_0-x|^2 = \lambda|y_0-x|^2
\]
then 
\[
	\lambda|y_0-x|^2>0
\] 
given that $y_0\not\in\overline\Omega$.
Then consider the segment 
\[
	\ell(t)=y_0+t(x-y_0)
\] 
for $t\in[0,\,1]$. We have $\ell(0)=y_0\not\in\overline\Omega$
and $\ell(1)=x\in\overline\Omega$. So there exists $\xi\in]0,\,1[$ such that 
\[
	\ell(\xi)\in\partial\Omega\,.
\]
We then have
\begin{equation*}
	\begin{split}
	\lambda|y_0-x|^2&=M_{\lambda}(M^{\lambda}(\chi_K))(x)\\[1.5ex]
			&\leq M^{\lambda}(\chi_K)(\ell(\xi))+\lambda|\ell(\xi)-x|^2\\[1.5ex]
			&=\lambda(1-\xi)^2|y_0-x|^2
	\end{split}
\end{equation*}
which is a contradiction as $0\leq (1-\xi)^2< 1$, thus $y_0\in\overline\Omega$ and 
\begin{equation*}
	\begin{split}
		M_{\lambda}(M^{\lambda}(\chi_K))(x)&=\underset{y\in\overline\Omega}{\inf}
			\Big(M^{\lambda}(\chi_K)(y)+\lambda|y-x|^2\Big)\\[1.5ex]
			&=\underset{y\in\overline\Omega}{\inf}
			\Big(M^{\lambda}_{\Omega}(\chi_K^{\Omega})(y)+\lambda|y-x|^2\Big)=M_{\lambda,\Omega}(M^{\lambda}_{\Omega}(\chi^{\Omega}_K))(x)
	\end{split}
\end{equation*}
which concludes the proof.
\hfill \qed\\

\noindent {\bf Proof of Theorem \ref{Thm.3.ImpNew}:}
The statements \eqref{Sec3.Eq.13} and \eqref{Sec3.Eq.14} are a consequence of Theorem \ref{Theo.LocLwTr}
applied to the auxiliary functions \eqref{Def.01.AuxFnct.M} which are bounded functions in $\Omega$ of the type
\eqref{Eq.01.Ext.Local} and \eqref{Eq.01.Ext.Glob} considered in Theorem \ref{Theo.LocLwTr}. 
Therefore, next, we will focus only on the proof of the local properties.

\textit{Part $(i)$:} Assume $x\in D$ and let $z_x$ be such that
\[
	M_{\lambda,\Omega}(f_{\overline\Omega,K}^M)(x)=f_{\overline\Omega,K}^M(z_x)+\lambda|z_x-x|^2\,.
\]
We will first show that $z_x\in(\Omega\setminus D)\cup\partial\Omega$, otherwise we have
\[
	M_{\lambda,\Omega}(f_{\overline\Omega,K}^M)(x)=M+\lambda|z_x-x|^2\,.
\]
Let $y_x\in \partial D$ such that $|y_x-x|=\dist(x,\Omega\setminus D)$ we have then that
\[
	M_{\lambda,\Omega}(f_{\overline\Omega,K}^M)(x)\leq f_{\overline\Omega,K}^M(y_x)+\lambda|y_x-x|^2\leq \underset{K}{\sup}\,f+\lambda\mathrm{diam}^2(D)
\]
so that 
\[
	M\leq \underset{K}{\sup}\,f+\lambda\mathrm{diam}^2(D)
\]
which is a contradiction. Now if $z_x\in\partial \Omega$, we have
\[
	M_{\lambda,\Omega}(f_{\overline\Omega,K}^M)(x)= \underset{K}{\inf}\,f+\lambda|z_x-x|^2\geq \lambda\dist^2(x,\partial\Omega)+
	\underset{K}{\inf}\,f\,.
\]
On the other hand,
\[
	M_{\lambda,\Omega}(f_{\overline\Omega,K}^M)(x) \leq  \underset{K}{\sup}\,f+\lambda\dist^2(x,\Omega\setminus D)
\]
so that,
\[
	\underset{K}{\inf}\,f + \lambda\dist^2(x,\partial\Omega)\leq \underset{K}{\sup}\,f+\lambda\dist^2(x,\partial D)\,.
\]
Let $u_x\in\partial\Omega$ such that $x-u_x=\dist(x,\partial\Omega)$ and $x^{\ast}\in\partial D$ such that $x^{\ast}$ lies on the line segment $[x,\,u_x]$.
Thus, 
\[
	\dist(x,\partial\Omega)=|x-x^{\ast}|+|x^{\ast}-u_x|\geq \dist(x,\partial\Omega)+\dist(\partial D,\partial\Omega)\,,
\]
hence
\[
\begin{split}
	\underset{K}{\sup}\,f-\underset{K}{\inf}\,f +\lambda\dist^2(x,\partial D)
	&\geq \lambda\left(\dist(x,\partial D)+\dist(\partial D,\partial\Omega)\right)^2 \\
	&\geq 
	\lambda\dist^2(x,\partial D) +	\lambda\dist^2(\partial D,\partial\Omega)\,,
\end{split}
\]
so that
\[
	\dist^2(\partial D,\partial\Omega)\leq \frac{O_f}{\lambda}
\]
which contradicts the assumption.

Assume $x\in\Omega\setminus D$ and that $M_{\lambda}(f_{\overline\Omega,K}^M)(x)=f_{\overline\Omega,K}^M(z_x)+\lambda|x-z_x|^2$ 
with $z_x\in\partial \Omega$. Then, since
\[
	M_{\lambda}(f_{\overline\Omega,K}^M)(x)\leq f_{\overline\Omega,K}^M)(x)=f(x)
\]
we have that 
\[
	\lambda\dist^2(x,\partial\Omega)\leq \lambda|z_x-x|^2\leq f(x)- \underset{K}{\inf}\,f \leq O_f\,,
\]
thus 
\[
	\dist^2(x,\partial\Omega)\leq \frac{O_f}{\lambda}\,,
\]
which is a contradiction. Thus $z_x\in\Omega\setminus D$.

\textit{Part $(ii)$:} Assume first $x\in D$ and let $z_x$ be such that
\[
	M^{\lambda}(M_{\lambda}(f_{\overline\Omega,K}^M))(x)=M_{\lambda}(f_{\overline\Omega,K}^M)(z_x)-\lambda|x-z_x|^2\,,
\]
then we need to distinguish the following three cases:
\[
\begin{array}{lll}	
	\displaystyle (a)\,z_x\in D\,; &
	\displaystyle (b)\,z_x\in \partial \Omega\,; &
	\displaystyle (c)\,z_x\in \Omega\setminus D\,. 
\end{array}
\]

\textit{Case $(a)$:} If $z_x\in D$, then by \textit{Part (i)} we conclude that $M_{\lambda}(f_{\overline\Omega,K}^M))(z_x)$
is determined by $f|_{\Omega\setminus D}$.

\textit{Case $(b)$:} If $z_x\in\partial\Omega$, we have
\[
	M_{\lambda}(f_{\overline\Omega,K}^M)(z_x)-\lambda|z_x-x|^2\leq \underset{K}{\inf}\,f-\lambda\dist^2(x,\partial\Omega)\,,
\]
and 
\[
	M^{\lambda}_{\Omega}(M_{\lambda,\Omega}(f_{\overline\Omega,K}^M))(x)\geq M_{\lambda,\Omega}(f_{\overline\Omega,K}^M)(x)
	=f_{\overline\Omega,K}^M(u_x)+\lambda|u_x-x|^2\geq f(u_x)\,,
\]
for some $u_x\in (\Omega\setminus D)\cup\partial\Omega$, thus
\[
	\dist^2(x,\partial\Omega)\leq\frac{O_f}{\lambda}
\]
which is a contradiction, thus $z_x\not\in\partial\Omega$.

\textit{Case $(c)$:} Let $z_x\in\Omega\setminus D$. We have then that
\[
	M_{\lambda}(f_{\overline\Omega,K}^M)(z_x)-\lambda|z_x-x|^2\leq f(z_x)-\lambda|z_x-x|^2\,,
\]
and
\[
	M^{\lambda}_{\Omega}(M_{\lambda,\Omega}(f_{\overline\Omega,K}^M))(x)\geq M_{\lambda,\Omega}(f_{\overline\Omega,K}^M)(x) 
	= f_{\overline\Omega,K}^M(u_x)+\lambda|z_x-x|^2\geq f(u_x)\,,
\]
for some $u_x\in (\Omega\setminus D)\cup \partial\Omega$. Thus
\[
	\lambda|z_x-u_x|^2\leq O_f\,,
\]
that is,
\[
	|z_x-u_x|\leq \sqrt{\frac{O_f}{\lambda}}\,.
\]
Therefore, 
\[
	\dist(z_x,\partial\Omega)\geq \dist(x,\partial\Omega)-|z_x-x|\geq 2\sqrt{\frac{O_f}{\lambda}}-\sqrt{\frac{O_f}{\lambda}}=\sqrt{\frac{O_f}{\lambda}}\,,
\]
and from \textit{Part $(i)$} we conclude that $M_{\lambda,\Omega}(f_{\overline\Omega,K}^M)(z_x)$ is determined by $f|_{\Omega\setminus D}$.

Now, let us assume that $x\in \Omega\setminus D$, $\dist(x,\partial\Omega)>2\sqrt{\tfrac{O_f}{\lambda}}$ and let $z_x$ be such that
\[
	M^{\lambda}_{\Omega}(M_{\lambda,\Omega}(f_{\overline\Omega,K}^M))(x)=M_{\lambda,\Omega}(f_{\overline\Omega,K}^M)(z_x)-\lambda|x-z_x|^2\,,
\]
then, also in this case, we need to analyze the following three claims:
\[
\begin{array}{lll}	
	\displaystyle (a)\,z_x\in D\,; &
	\displaystyle (b)\,z_x\in \partial \Omega\,; &
	\displaystyle (c)\,z_x\in \Omega\setminus D\,. 
\end{array}
\]
\mbox{}

\textit{Case $(a)$:} By \textit{Part $(i)$} we conclude that $M_{\lambda,\Omega}(f_{\overline\Omega,K}^M)(z_x)$ is determined by $f|_{\Omega\setminus D}$.

\textit{Case $(b)$:} As in the similar case seen previusly, we can mke the same arguments and show that $z_x\not\in \partial \Omega$.

\textit{Case $(c)$:} By imitating the arguments made previously for the same case, let $z_x\in\Omega\setminus D$ be such that
\[
	M_{\lambda}(f_{\overline\Omega,K}^M)(z_x)-\lambda|z_x-x|^2\leq f(z_x)-\lambda|z_x-x|^2\,,
\]
and
\[
	M^{\lambda}_{\Omega}(M_{\lambda,\Omega}(f_{\overline\Omega,K}^M))(x)\geq M_{\lambda,\Omega}(f_{\overline\Omega,K}^M)(x) 
	= f_{\overline\Omega,K}^M(u_x)+\lambda|z_x-x|^2\geq f(u_x)\,,
\]
for some $u_x\in (\Omega\setminus D)\cup \partial\Omega$. Thus,
\[
	|z_x-u_x|\leq \sqrt{\frac{O_f}{\lambda}}\,.
\]
Therefore, 
\[
	\dist(z_x,\partial\Omega)\geq \dist(x,\partial\Omega)-|z_x-x|\geq \sqrt{\frac{O_f}{\lambda}}\,,
\]
and from \textit{Part $(i)$} we conclude that $M_{\lambda,\Omega}(f_{\overline\Omega,K}^M)(z_x)$ is determined by 
$f|_{\Omega\setminus D}$. \hfill \qed\\


\noindent {\bf Proof of Theorem \ref{Sec4.Theo.CnvrgMor}:}
We have that 
\begin{equation}\label{Sec6.Theo.CnvrgMor.Eq.01}
	M_{\lambda}f(x_k)=f(z_{x_k})+\lambda|z_{x_k}-x_k|^2\leq f(x_k)\,,
\end{equation}
with some $z_{x_k}\in\mathbb{R}^n$.
By Proposition \ref{Sec2.Pro.MoC}, 
\begin{equation}\label{Sec6.Theo.CnvrgMor.Eq.02}
	\begin{split}
		\lambda|z_{x_k}-x_k|^2	&\leq f(x_k)-f(z_{x_k})\\[1.5ex]
					&\leq \omega_f(|z_{x_k}-x_k|)\\[1.5ex]
					&\leq a|z_{x_k}-x_k|+b
	\end{split}
\end{equation}
which yields the following estimate
\begin{equation}\label{Sec6.Theo.CnvrgMor.Eq.03}
		|z_{x_k}-x_k|\leq \frac{a+\sqrt{a^2+4b\lambda}}{2\lambda} 
		= \frac{a}{\lambda}+\sqrt{\frac{b}{\lambda}}\,.
\end{equation}
Since $\omega_f$ is non-decreasing, from \eqref{Sec6.Theo.CnvrgMor.Eq.02} and  \eqref{Sec6.Theo.CnvrgMor.Eq.03} we obtain
\begin{equation}\label{Sec6.Theo.CnvrgMor.Eq.04}
		|z_{x_k}-x_k|	\leq \sqrt{\frac{\omega_f(a/\lambda+\sqrt{b/\lambda})}{\lambda}}
				= \frac{d(\lambda)}{\sqrt{\lambda}}\,,
\end{equation}
where we have set $d(\lambda)=\sqrt{\omega_f(a/\lambda+\sqrt{b/\lambda})}$. 
Now, it is not difficult to realize that we can give the following representation for $z_{x_k}$
which enters \eqref{Sec6.Theo.CnvrgMor.Eq.01} in terms of the grid points as follows
\footnote{
Let $x_{z_{x_k}}$ be for instance the closest grid point to $z_{x_k}$ and set 
$z_{x_k}=x_k+x_{z_{x_k}}-x_k+z_{x_k}-x_{z_{x_k}}$. Noting that 
$x_{z_{x_k}}-x_k=hr_{x_k}$, $h$ the grid size and $r_{x_k}\in\mathbb{Z}^n$, and 
$z_{x_k}-x_{z_{x_k}}=h\delta$ with $\delta\in\mathbb{R}^n$, $0\leq\delta_i\leq 1$, 
$i=1,\ldots,n$, we obtain \eqref{Sec6.Theo.CnvrgMor.Eq.05}.
}
\begin{equation}\label{Sec6.Theo.CnvrgMor.Eq.05}
		z_{x_k}=x_k+h(r_{x_k}+\delta)
\end{equation}
$r_{x_k} \in \mathbb{Z}^n$ and $\delta\in\mathbb{R}^n$ such that $|\delta|< \sqrt{n}$. 
Thus, we have
\begin{equation}\label{Sec6.Theo.CnvrgMor.Eq.06}
	\begin{split}
		M_{\lambda}^h f(x_k)& \geq M_{\lambda}f (x_k) = \lambda|z_{x_k}-x_k|^2 + f(z_{x_k})\\[1.5ex]
				    & =	\lambda h^2 |r_{x_k}+\delta|^2 +	f(x_k+h(r_{x_k}+\delta))	\\[1.5ex]
				    &   \geq \lambda h^2 |r_{x_k}|^2 + \lambda h^2 |\delta|^2 + 2\lambda h^2 r_{x_k}\cdot \delta
				       + f(x_k+hr_{x_k})-\omega_f(|h\delta|)	\\[1.5ex]
				    & \geq M_{\lambda}^h f(x_k)+\lambda h^2|\delta|^2+2\lambda h^2 r_{x_k}\cdot \delta 
				    -\omega_f(|h\delta|)\\[1.5ex]
				    & \geq M_{\lambda}^h f(x_k)+2\lambda h^2 r_{x_k}\cdot \delta 
				    -\omega_f(|h\delta|)\,,
	\end{split}
\end{equation}
where we have used the fact that 
\begin{equation}\label{Sec6.Theo.CnvrgMor.Eq.07}
	M_{\lambda}^h f(x_k)\leq \lambda  h^2|r_{x_k}|^2+f(x_k+hr_{x_k})\,.
\end{equation}
By the Cauchy-Schwartz inequality,
\begin{equation}\label{Sec6.Theo.CnvrgMor.Eq.08}
	r_{x_k}\cdot \delta \geq - |r_{x_k}||\delta| \geq -\sqrt{n}|r_{x_k}|
\end{equation}
and since $\omega_f$ is non-decreasing, we have that
\begin{equation}\label{Sec6.Theo.CnvrgMor.Eq.09}
	\omega_f(h|\delta|)\leq \omega_f(h\sqrt{n})\,.
\end{equation}
As a result, taking into account \eqref{Sec6.Theo.CnvrgMor.Eq.08}  and 
\eqref{Sec6.Theo.CnvrgMor.Eq.09} into \eqref{Sec6.Theo.CnvrgMor.Eq.06} yields
\begin{equation}\label{Sec6.Theo.CnvrgMor.Eq.10}
		M_{\lambda}^h f(x_k) \geq M_{\lambda}^h f(x_k)-2\lambda h^2\sqrt{n} |r_{x_k}| 
				    -\omega_f(h\sqrt{n})\,.
\end{equation}
From \eqref{Sec6.Theo.CnvrgMor.Eq.05}, 
\begin{equation}\label{Sec6.Theo.CnvrgMor.Eq.11}
	|z_{x_k}-x_k|=h|r_{x_k}+\delta|
\end{equation}
which, combined with \eqref{Sec6.Theo.CnvrgMor.Eq.04}, gives
\begin{equation}\label{Sec6.Theo.CnvrgMor.Eq.12}
	h|r_{x_k}+\delta|\leq\frac{d(\lambda)}{\sqrt{\lambda}}
\end{equation}
which implies that
\begin{equation}\label{Sec6.Theo.CnvrgMor.Eq.13}
	-h|r_{x_k}|\geq -h\sqrt{n} - \frac{d(\lambda)}{\sqrt{\lambda}}
\end{equation}
Using \eqref{Sec6.Theo.CnvrgMor.Eq.13} into \eqref{Sec6.Theo.CnvrgMor.Eq.10}, and by taking into account 
of \eqref{Sec6.Theo.CnvrgMor.Eq.01}, we get
\begin{equation}\label{Sec6.Theo.CnvrgMor.Eq.14}
		M_{\lambda}^h f(x_k) \geq M_{\lambda} f(x_k)  \geq M_{\lambda}^h f(x_k)-2\lambda {n}h^2-2h\sqrt{\lambda}d(\lambda)
					    -\omega_f(h\sqrt{n})\,,
\end{equation}
which concludes the proof.\hfill \qed\\

\noindent {\bf Proof of Corollary \ref{Sec4.Cor.CnvrgMorLips}:} The result follows from Theorem \ref{Sec4.Theo.CnvrgMor} by taking 
$\omega_f(t)=Lt$. \hfill \qed\\


\noindent {\bf Proof of Theorem \ref{Sec4.Theo.DisSch}:}
For the sake of clarity, in the following,  we present the proof only for $n=1$ and $n=2$, given that the proof of the 
results for $n>2$ build on the arguments required for $n-1$.

\noindent \textit{Case $n=1$:} We prove the statement by induction. Clearly, the statement 
is true for $m=0,1$ and for all $x_k=kh$, $k\in\mathbb{Z}$.
Suppose now $f_m=g_m$ for some $m\geq 1$. Then for $m+1$, we have
\begin{equation}
	\begin{split}
		f_{m+1}(x_k)&=\min\{f_m(x_k),\,f_m(x_k-h)+\lambda \tau_{m+1}h^2, 
		f_m(x_k+h)+\lambda \tau_{m+1}h^2\}\\[1.5ex]
		&=\min\{A,\, B,\,C\}\,
	\end{split}
\end{equation}
where we have let $A:=f_m(x_k)$, $B:=f_m(x_k-h)+\lambda \tau_{m+1}h^2$ and $C:=f_m(x_k+h)+\lambda \tau_{m+1}h^2$.
By the inductive assumption, we have then
\begin{equation}
	A=g_m(x_k)=\min\{ f_0(x_k+jh)+\lambda j^2h^2,\; |j|\leq m\},
\end{equation}
and
\begin{equation}
	\begin{split}
		B & =g_m(x_{k-1})+\lambda \tau_{m+1}h^2\\[1.5ex]
		  & =\min \{f_0(x_{k-1}+jh)+\lambda(j^2+\tau_{m+1})h^2,\; |j|\leq m\}\\[1.5ex]
		  & =\min \{f_0(x_{k-(m+1)}+\lambda (m^2+\tau_{m+1})h^2,\; B_1\}\\[1.5ex]
		  & =\min\{f_0(x_{k-(m+1)}+\lambda(m+1)^2h^2,\; B_1\}.
	\end{split}
\end{equation}
where we have let
\begin{equation}
	B_1:=\min\{ f_0(x_{k-1}+jh)+\lambda (j^2+\tau_{m+1})h^2,\;  -(m-1)\leq j\leq m\}\,.
\end{equation}
If we change index $j-1=r$ in $B_1$, we have
\begin{equation}
	B_1=\min\{ f_0(x_{k}+rh)+\lambda ((r+1)^2+\tau_{m+1})h^2,\; -m\leq r\leq m-1\}.
\end{equation}
Now we show that
\begin{equation}
	f_0(x_{k}+rh)+\lambda ((r+1)^2+\tau_{m+1})h^2\geq f_0(x_{k}+rh)+\lambda r^2h^2,\quad -m\leq r\leq m-1,
\end{equation}
as the right hand side of the above inequality is one of the terms to be minimised in $A$ above. 
This inequality is equivalent to
$(r+1)^2+\tau_{m+1}\geq r^2$ for $-m\leq r\leq m-1$, which is equivalent 
to $2r+1+2m+1\geq 0$ and thus to $2(m+r)+2\geq 0$.
As $-m\leq r\leq m-1 $, the last inequality above clearly holds.
Thus we have
\begin{equation}
	\min\{A,\, B\}=\min\{g_m(x_k),\, f_0(x_{k-(m+1)})+\lambda (m+1)^2h^2\}.
\end{equation}
Similarly, we can easily see that 
\begin{equation}
	\min\{A,\, C\}=\min\{g_m(x_k),\, f_0(x_{k+(m+1)})+\lambda (m+1)^2h^2\}\,,
\end{equation}
thus, 
\begin{equation}
\begin{split}
	f_{m+1}(x_k)&=\min\{A,\, B,\, C\}\\[1.5ex]
	&=\min\{g_m(x_k),\, f_0(x_{k-(m+1)})+\lambda (m+1)^2h^2,\, f_0(x_{k+(m+1)})+\lambda (m+1)^2h^2\}\\[1.5ex]
	&=g_{m+1}(x_{k}).
\end{split}
\end{equation}
It therefore follows that the statement holds for $n=m+1$ and by induction holds for all $n=1,2,\dots$.
\hfill\qed\\


\noindent \textit{Case $n=2$:}  
The proof of the two-dimensional case is very similar to the one-dimensional case but 
is slightly more involved. Again, we prove the statement $f_n=g_n$ by induction. 
For $m=0,1$, the statement is clearly true. Suppose for some $m\geq 1$, we have
$f_m=g_m$. Then for $m+1$ and each $(x_k, y_j)$, we have $f_{m+1}(x_k,y_j)=\min\{A,\, B\}$, where
\begin{equation}
	\begin{split}
		A &=f_m(x_k,y_j)\\[1.5ex]
		  &=g_m(x_k,y_j)\\[1.5ex]
		  &=\min\{f_0(x_k+rh,y_j+sh)+\lambda(r^2+s^2)h^2,\; |r|\leq m,\,|s|\leq m\},
	\end{split}
\end{equation}
and
\begin{equation}
	B=\min\{g_m(x_k+ph,y_j+qh)+\lambda (p^2+q^2)\tau_{m+1}h^2,\; |p|\leq 1,\, |q|\leq 1,\,|p|+|q|\neq 0\}.
\end{equation}
There are eight terms to be minimised in $B$. We consider first two typical cases.

Case $(a)$: $p=0$ and $q=1$. This is similar to the one-dimensional case as we have
\begin{equation}
	\begin{split}
& g_m(x_k+ph,y_j+qh)+\lambda \tau_{m+1}h^2=g_m(x_k,y_{j+1})+\lambda\tau_{m+1}h^2\\[1.5ex]
&\phantom{xxx}=\min\{\min\{f_0(x_k+rh,y_{j}+(m+1)h)+\lambda(r^2+m^2+\tau_{m+1})h^2,\; |r|\leq m\},\; C_1 \}\\[1.5ex]
&\phantom{xxx}=\min\{\min\{f_0(x_{k+r},y_{j+(m+1)})+\lambda(r^2+(m+1)^2)h^2,\; |r|\leq m\},\; C_1 \}\\[1.5ex]
&\phantom{xxx}=\min\{B_1, C_1\}\,
	\end{split}
\end{equation}
where we have let 
\begin{equation}
\begin{split}
	B_1&:=\min\{f_0(x_{k+r},y_{j+(m+1)})+\lambda(r^2+(m+1)^2)h^2,\; |r|\leq m\}\\[1.5ex]
	C_1&:=\min\{ f_0(x_k+rh,y_{j+1}+sh)+\lambda(r^2+s^2+\tau_{m+1})h^2,\, \, 
		|r|\leq m,\, -m\leq s\leq m-1\}
\end{split}
\end{equation}
Note that all terms to be minimised in $B_1$ are in the definition of $g_{m+1}(x_k,y_j)$ but not in
$g_{m}(x_k,y_j)$. As in the one-dimensional case, every term to be minimised in $C_1$
is greater than a corresponding term to be minimised in $A$ as we have seen for the Case $(i)$.
Therefore
\begin{equation}
	\min\{A, \, g_m(x_k,y_{j+1})+\lambda\tau_{m+1}h^2\}=\min\{A,\, B_1\}\,.
\end{equation}

Case $(b)$: $p=1$ and $q=-1$. In this case, we have
\begin{equation}
	\begin{split}
   &g_m(x_k+ph,y_j+qh)+\lambda (p^2+q^2)\tau_{m+1}h^2=g_m(x_{k+1},y_{j-1})+2\lambda\tau_{m+1}h^2\\[1.5ex]
&\phantom{xxx}=\min\{f_0(x_{k+1}+rh,y_{j-1}+sh)+\lambda(r^2+s^2+2\tau_{m+1})h^2,\;  |r|\leq m,\, |s|\leq m \}\\[1.5ex]
&\phantom{xxx}=\min\{f_0(x_{k+(m+1)},y_{j-(m+1})+\lambda(2m^2+2\tau_{m+1})h^2,\; D_1\}\\[1.5ex]
&\phantom{xxx}=\min\{f_0(x_{k+(m+1)},y_{j-(m+1)})+2\lambda (m+1)^2h^2,\; D_1\}.
	\end{split}
\end{equation}
In the above, $f_0(x_{k+(m+1)},y_{j-(m+1})+2\lambda (m+1)^2h^2$ is in the definition of $g_{m+1}(x_k,y_j)$
which is one of the terms to be minimised that is not in the definition of
$g_{m}(x_k,y_j)$. We also have
\begin{equation}
	\begin{split}
	D_1=\min\Big\{f_0(x_{k+1}+rh,y_{j-1}+sh)+\lambda(r^2+s^2+2\tau_{m+1})h^2,& -m\leq r\leq m-1,\\[1.5ex]
									     & -(m-1)\leq s\leq m \Big\}.
	\end{split}
\end{equation}
Again, as in the Case $(i)$, every term involved in the minimisation in $D_1$ is greater than a 
corresponding term in $A$. Therefore, we have
\begin{equation}
	\min\Big\{A,\, g_m(x_{k+1},y_{j-1})+2\lambda\tau_{m+1}h^2\Big\}=
	\min\Big\{A,\,f_0(x_{k+(m+1)},y_{j-(m+1})+2\lambda (m+1)^2h^2\Big\}\,.
\end{equation}
For the Case $(c)$: $p=0$, $q=-1$; Case $(d)$: $p=1$, $q=0$ and Case $(e)$: $p=-1$, $q=0$, the proofs 
are similar to that of Case $(a)$
while for the Case $(f)$ $p=1$, $q=1$; Case $(g)$: $p=-1$, $q=-1$ and Case $(h)$: $p=-1$, $q=1$, 
the proofs are similar to that of Case $(b)$. Therefore, we have
\begin{equation}
	f_{m+1}(x_k,y_j)=\min\{A,\, B\}=\min\{A,\, U,\, V\}
\end{equation}
where 
\begin{equation}
	\begin{split}
		U:=\min\Bigg\{\min\Big\{&f_0(x_{k+r},y_{j+(m+1)}),\, f_0(x_{k+r},y_{j-(m+1)}),
		  f_0(x_{k+(m+1)},y_{j+r}),\,\\[1.5ex]
		& f_0(x_{k-(m+1)},y_{j+r})\Big\}
		+\lambda ((m+1)^2+r^2)h^2,\; |r|\leq m\Bigg\},
	\end{split}
\end{equation}
and
\begin{equation}
	\begin{split}
		V:=\min\Big\{&f_0(x_{k+(m+1)},y_{j+(m+1)}),\, f_0(x_{k+(m+1)},y_{j-(m+1)}),\,
		f_0(x_{k-(m+1)},y_{j+(m+1)}),\\[1.5ex]
		&f_0(x_{k-(m+1)},y_{j-(m+1)})\Big\}+2\lambda (m+1)^2h^2.
	\end{split}
\end{equation}
By comparing the terms involved in $g_{m+1}(x_k,y_j)$ and those involved in $A$, $U$ and $V$, we see that
\begin{equation}
	f_{m+1}(x_k,y_j)=\min\{A,\, U,\, V\}=g_{m+1}(x_k,y_j).
\end{equation}
Thus the statement holds for $n=m+1$. By induction, we conclude then that 
$f_{n}(x_k,y_j)=g_{n}(x_k,y_j)$ holds for all $n=0,1,2,\dots$ and for all $x_k=kh,\, y_j=jh$ 
with $k,\,j\in\mathbb{Z}$. The proof is finished.
\hfill \qed\\


\noindent {\bf Proof of Proposition \ref{Sec4.Prop.Est}.} 
For any grid point $x_k$, we have that there exists a $r_{x_k}\in\mathbb{Z}^n$ such that
\begin{equation}\label{Sec6.Eq.Est.01}
	M_{\lambda}^{h}(f)(x_k)=\lambda|r_{x_k}h|^2+f(x_k+r_{x_k}h)\leq f(x_k)
\end{equation}
thus, 
\[
	\lambda|r_{x_k}h|^2\leq  f(x_k)-f(x_k+r_{x_k}h)\leq \osc(f)
\]
which implies that
\[
	|r_{x_k}|\leq \frac{1}{h}\sqrt{\frac{\osc(f)}{\lambda}}\,.
\]
As a result, by taking in \eqref{Sec4.Eq.02}
\[
	m\geq \floor{\frac{1}{h}\sqrt{\frac{\osc(f)}{\lambda}}}+1
\]
where $\floor{x}$ denotes the integer part of $x$, we are guaranteed that  
$r_{x_k}$, which enters \eqref{Sec6.Eq.Est.01}, belongs to the feasible set of \eqref{Sec4.Eq.02}
and is, therefore, attained by  Algorithm \ref{Algo:MoreauEnv}. This concludes the proof.
\hfill\qed\\




\begin{thebibliography}{44}
\bibitem{Zha08a} K. Zhang, 
		Compensated convexity and its applications, 
		Anal. Non-Lin. H. Poincar\'e Inst. 25 (2008) 743--771.

\bibitem{Zha08b} K. Zhang, 
		Convex analysis based smooth approximations of maximum functions and squared-distance functions,  
		J. Nonlinear Convex Anal.  9  (2008) 379--406.

\bibitem{ZOC15a} K. Zhang, A. Orlando, E. Crooks,
		Compensated convexity and Hausdorff stable geometric singularity extractions,
		Mathematical Models and Methods in Applied Sciences 25 (2015) 747--801,
		DOI: 10.1142/S0218202515500189.

\bibitem{ZOC15b} K. Zhang, A. Orlando, E. Crooks, 
		Compensated convexity and Hausdorff stable extraction of intersections for smooth manifolds, 
		Mathematical Models and Methods in Applied Sciences 25 (2015) 839--873,
		DOI: 10.1142/S0218202515500207.
\bibitem{ZCO16b} K. Zhang, E. Crooks, A. Orlando, 
		Compensated convex transforms and geometric singularity extraction from semiconvex functions 
		(in Chinese). Scientia Sinica Mathematica 46 (2016) 1--22. DOI: 10.1360/N012015-00339,
		(revised English version available at \url{https://arxiv.org/abs/1610.01451})

\bibitem{ZCO15c} K. Zhang, E. Crooks, A. Orlando,
		Compensated convexity, multiscale medial axis maps and sharp regularity of the squared distance function,
		SIAM Journal on Mathematical Analysis 47 (2015) 4289--4331.

\bibitem{ZCO16a} K. Zhang, E. Crooks, A. Orlando,  
		Compensated convexity methods for approximations and interpolations of sampled functions in 
		Euclidean spaces: Theoretical Foundations. 		
		SIAM Journal on Mathematical Analysis 48 (2016) 4126--4154.		

\bibitem{Obe08} {A.M.} Oberman, 
		Computing the convex envelope using a nonlinear partial differential equation.
		Math. Models Methods Appl. Sci. 18 (2008) 759--780.

\bibitem{Mor65}	{J.-J.} Moreau,
		Proximat\'e dualit\'e dans un espace Hilbertien,
		Bull. Soc. Math. Fr. 93 (1965) 273--299.

\bibitem{Mor66}	{J.-J.} Moreau,
		Fonctionnelles convexes.
		S\'{e}minaire "Sur les \'{e}quations aux d\'{e}riv\'{e}es partielles".
		Lecture Notes, Coll\'{e}ge de France, 1966.

\bibitem{LL86}	J.M. Lasry and P.L. Lions,
		A remark on regularization in Hilbert Spaces,
		Israel Math. J. 55 (1986) 257--266.

\bibitem{AA93}	H. Attouch, D. Aze,
		Approximations and regularizations of arbitrary functions in Hilbert spaces by the
		Lasry-Lions methods,
		Anal. Non-Lin. H. Poincar\'e Inst. 10 (1993) 289--312.

\bibitem{CS04}	P. Cannarsa,  C. Sinestrari,
		Semiconcave Functions, Hamilton-Jacobi Equations and Optimal Control,
		Birkh\"auser, Boston,  2004

\bibitem{Ser82}	J. Serra,
		Image Analysis and Mathematical Morphology.
		Volume 1, Academic Press, London, 1982.

\bibitem{Soi04}	P. Soille,
		Morphological Image Analysis: Principles and Applications,
		Springer, Berlin, 2nd edition, 2004.

\bibitem{Shi09}	F.Y. Shih,
		Image Processing and Mathematical Morphology,
		CRC Press, Boca Raton, USA, 2009.

\bibitem{BS92}	R. Van den Boomgaard, A.W.M. Smelders,
		The morphological structure of images,
		in Proceedings 11th IAPR International Conference
		on Pattern Recognition. The Hague, The Netherlands: IEEE Computer
		Society Press, Los Alamitos, CA, (1992),  268--271.

\bibitem{BDMS96} R. Van Den Boomgaard, L. Dorst, S. Makram--Ebeid, J. Schavemaker
		Quadratic structuring functions in mathematical morphology. 
		In: Maragos P., Schafer R.W., Butt M.A. (eds) 
		Mathematical Morphology and its Applications to Image and Signal Processing. 
		Computational Imaging and Vision, vol 5, 1996, Springer, Boston, MA.

\bibitem{Jac92}	P.T. Jackway,
		Morphological scale-space,
		in Proceedings 11th IAPR International Conference on Pattern Recognition.
		The Hague, The Netherlands: IEEE Computer Society Press, Los Alamitos, CA, (1992),
		252--255.

\bibitem{RW98}	{R.T.} Rockafellar, {R. J-B.} Wets,
		Variational Analysis, 
		Springer, Berlin, 1998. 

\bibitem{Car19}	M. Carlsson,
		On convex envelopes and regularization of non-convex functionals without moving global minima.
		Journal of Optimization Theory and Applications 183 (2019) 66--84.

\bibitem{ZCO18} K. Zhang, E. Crooks, A. Orlando,  
		Compensated convexity methods for approximations and interpolations of sampled functions in 
		Euclidean Spaces: Applications to contour lines, sparse data and inpainting.
		SIAM J. Imaging Sciences 11 (2018) 2368--2428.

\bibitem{OCZ20}	A. Orlando, E. Crooks, K. Zhang,
		Compensated convex based transforms for image processing and shape interrogation.
		\textit{Accepted for publication in: 
		Ke Chen, Sch\"onlieb C.-B., Xue-Cheng Tai, Younes L. (Eds),
		Handbook of Mathematical Models and Algorithms in Computer Vision and Imaging, Springer}

\bibitem{HL01}	{J.-B.} Hiriart-Urruty, C. Lemar\'echal,
		Fundamentals of Convex Analysis,
		Springer, New York, 2001.

\bibitem{Roc70} {R.T.} Rockafellar, 
		Convex Analysis, 
		Princeton Univ. Press, New Jersey, 1970.

\bibitem{CLSW98} F.H. Clarke, Yu.S. Ledyaev, R.J. Stern, P.R. Wolenski, 
		Nonsmooth Analysis and Control Theory,
		Springer-Verlag, New York, 1998.

\bibitem{DL93}	R.A. DeVore, G.G. Lorentz,
		Constructive Approximation,
		Springer-Verlag, Berlin, 1993.

\bibitem{Luc09b} Y. Lucet, 
		What shape is your conjugate? A survey of computational convex analysis and its applications,
		SIAM J. on Optimization 20 (2009) 216--250.

\bibitem{Luc06} Y. Lucet, 
		Fast Moreau envelope computation I: numerical algorithms,
		Numer. Algor. 43 (2006) 235--249.

\bibitem{Luc09a} Y. Lucet, 
		New sequential exact Euclidean distance transform algorithms based on convex analysis,
		Image and Vision Computing 27 (2009) 37--44.

\bibitem{Cor96} L. Corrias,
		Fast Legendre-Fenchel transform and applications to Hamilton-Jacobi equations and conservations laws.
		SIAM Journal on Numerical Analysis 33 (1996) 1534--1558.

\bibitem{FH12}	P.F. Felzenszwalb, D.P. Huttenlocher,
		Distance transforms of sampled functions,
		Theory of Computing 8 (2012) 415--428.

\bibitem{HM94}  C. Huang, O. Mitchell, 
		A Euclidean distance transform using grayscale morphology decomposition, 
		IEEE Transactions on Pattern Analysis and Machine Intelligence 16 (1994) 443--448.

\bibitem{SM92}	F.Y. Shih, O. Mitchell,  
		A mathematical morphology approach distance transformation, 
		IEEE Transactions on Image Processing 1 (1992) 197--204.

\bibitem{BDH96}	C.B. Barber, D.P. Dobkin, H. Huhdanpaa,
		The quickhull algorithm for convex hulls, 
		ACM Trans. Math. Software 22 (1996) 469--483.

\bibitem{LLS92} L. Lam, L. Seong-Whan Lee, Y.S. Ching, Y. Suen, 
		Thinning methodologies-A comprehensive survey,
		IEEE Transactions on Pattern Analysis and Machine Intelligence 14 (1992) 869--885.

\bibitem{CHN05} R.H. Chan, {C.-W.} Ho, M. Nikolova,
		Salt-and-pepper noise removal by median-type noise detectors and detail-preserving regularization,
		IEEE Trans. Image Processing 14 (2005) 1479--1485.

\bibitem{CCM07} {J.-F.} Cai, R. Chan, B. Morini,
		Minimization of an edge-preserving regularization functional by conjugate gradient type methods,
		in Image Processing Based on Partial Differential Equations,
		{X.-C.} Tai, {K.-A.} Lie, T. F. Chan, S. Osher (Eds), Springer (2005) 109--122.

\bibitem{Get12} P. Getreuer,
		Total variation inpainting using {S}plit {B}regman,
		Image Processing On Line 2 (2012) 147--157.
\end{thebibliography}
\end{document}